\documentclass[english]{svjour3}
\usepackage[T1]{fontenc}
\usepackage[latin9]{inputenc}
\usepackage{algorithm2e}
\usepackage{amsmath}
\usepackage{amssymb}
\usepackage{hyperref}

\makeatletter

\providecommand{\tabularnewline}{\\}

\newenvironment{svmultproof}{\begin{proof}}{\qed\end{proof}}

\newcommand{\pPi}{\mathnormal{\Pi}}
\newcommand{\im}{\operatorname{im}}

\newcommand{\rank}{\operatorname{rank}}

\newcommand{\R}{\mathbb{R}}
\newcommand{\Natu}{\mathbb{N}}
\newcommand{\Scan}{S_{\textrm{can}}}
\newcommand{\Ncan}{N_{\textrm{can}}}
\newcommand{\dt}{{\rm d}}
\newcommand{\dist}{\operatorname{dist}}

\newcommand{\T}{\ast}

\LinesNumbered

\makeatother

\usepackage{babel}

\begin{document}
\title{On the computation of accurate initial conditions for linear higher-index
differential-algebraic equations and its application in initial value
solvers}
\titlerunning{Accurate initial conditions for DAEs}
\author{Michael Hanke \and Roswitha März}
\institute{
Michael Hanke\at KTH Royal Institute of Technology, Stockholm, Sweden\\
\email{hanke@kth.se}
\and
Roswitha M\"arz\at Humboldt University, Berlin, Germany\\
\email{Roswitha.Maerz@hu-berlin.de}
}
\maketitle
\begin{abstract}
In contrast to regular ordinary differential equations, the problem of accurately setting initial conditions just emerges in the context of differential-algebraic equations where the dynamic degree of freedom of the system is smaller than the absolute dimension of the described process, and the actual lower-dimensional configuration space of the system is deeply implicit. For linear higher-index   differential-algebraic equations,  we develop an appropriate  numerical method based on properties of canonical subspaces and on the so-called geometric reduction. Taking into account the fact that higher-index differential-algebraic equations lead to ill-posed problems in  naturally given norms, we modify this  approach to serve as transfer conditions from one time-window to the next  in a time stepping procedure and  combine it with window-wise overdetermined least-squares collocation  to construct the first fully numerical solvers for higher-index initial-value problems.

\keywords{accurate initial conditions \and higher-index differential-algebraic equations \and canonical subspaces \and geometric reduction \and least-squares collocation}

Classification: 65L80 \and 65L05 \and 65L08 \and 65L60

\end{abstract}

\section{Introduction}

The numerical solution of initial value problems for higher-index
differential-algebraic equations (DAEs) is known to be a delicate
problem. Classical numerical methods like finite difference or collocation
methods loose their convergence properties known from their application
to explicit ordinary differential equations (ODEs). While appropriate
modifications of classical methods are well established for index-1
and special classes of low-index DAEs featuring usable special structures, their use will fail for
more general higher-index DAEs (e.g. \cite{BCP96,HaiWaII,KuhMeh06,LMT,TMS14}).

Known solvers for general higher-index DAEs as in \cite{RaRh94,KuMeRaWe97,EstLam21} are based on the generation and evaluation of so-called derivative array functions, which is extremely costly. The arrays are provided manually or by automatic differentiation. Then the index of the given DAE is reduced \emph{a-priori} or \emph{a-posteriori} and the necessary derivatives are determined using the derivative  array functions.
By the way,  these solvers need so-called \emph{consistent initial values} to start, but the computation of consistent initial values is by itself a hard problem \cite{LePeGe91,EstLam21b}.

Higher-index DAEs are known to lead to ill-posed problems in natural function spaces.
In order to overcome the problems related to the ill-posedness of
higher-index DAEs, in \cite{HaMaTiWeWu17,HaMaTi19,HaMa21C} we proposed
a family of least-squares methods for solving them. These methods do not need
any preprocessing nor higher derivatives of the problem data as those
proposed earlier \cite{RaRh94,KuMeRaWe97,EstLam21}. However, this
method needs \emph{accurate initial conditions} in the sense of \cite[Definition 2.3]{LMW}
and \cite[Theorem 2.52]{LMT}.

Accurate initial conditions for DAEs are something else than consistent initial values.
Roughly speaking, they consist of exactly $l$ conditions according to the degree of freedom $l$ in order to capture exactly one solution. To formulate reasonable $l$ conditions requires, of course, the knowledge of $l$ itself and, additionally information on canonical subspaces, and hence the determination of accurate initial
conditions is also a nontrivial problem. In a recent paper \cite{HaMa231},
we were able to substantiate a constructive way  for just this problem.  

The present paper is devoted to the construction of a numerical algorithm
based on the results of \cite{HaMa231} and its analysis for linear
DAEs. Then this algorithm will be combined with the time-stepping
method proposed in \cite{HaMaPad121} in order to establish a fully numerical initial
value solver for linear higher-index DAEs. The latter requires accurately
stated transfer conditions from one step to the next. We will use
the algorithm developed below for this purpose.

The paper is structured as follows: In Section~\ref{sec:AIC} we provide basic facts concerning the flow defined by the DAE, the canonical subspaces,  accurate initial and transfer condition.
In Section~\ref{s.Algoritm} we first summarize the method proposed in \cite{HaMa231} and develop an
implementable algorithm for determining a matrix capable of providing
accurate initial conditions and discuss properties of the associated tools.
Section~\ref{s.PropAICnum} is devoted
to a deeper analysis of the algorithm's properties. The algorithm
has been implemented with the goal of robustness and efficiency. Some
numerical examples of its performance are provided.
Finally, in Section~\ref{sec:An-initial-value} the convergence of a time stepper based on the least-squares
method of \cite{HaMaPad121} combined with the algorithm developed
in the present paper is shown and some numerical results are presented. Some technical details are
collected in the appendix.

Throughout the paper we will use the Euclidean norm $\lvert x\rvert$ for vectors $x\in\R^m$ and the spectral norm
$\lvert M\rvert$ for matrices $M\in\R^{m\times l}$.

\section{Formulating accurate initial and transfer conditions\label{sec:AIC}}

\subsection{The flow structure of a regular DAE}\label{subs.regular}
In this paper, we are interested in solving initial value problems for
linear DAEs
\begin{equation}
A(t)(Dx(t))'+B(t)x(t)=q(t),\quad t\in[a,b].\label{eq:DAE}
\end{equation}
Here, $x:[a,b]\rightarrow\R^{m}$ is the unknown vector-valued function
defined on the finite interval $[a,b]$. In order to apply collocation methods later we suppose continuously differentiable solutions satisfying the equation \eqref{eq:DAE} pointwise.
The right-hand side $q:[a,b]\rightarrow\R^{m}$
shall be sufficiently smooth. Moreover, sufficient smoothness is also
assumed of the coefficient matrix functions $A:[a,b]\rightarrow\R^{m\times k}$
and $B:[a,b]\rightarrow\R^{m\times m}$. With the intention of creating an appropriate collocation approach we assume that the DAE features
an explicit splitting of the components of $x$ into differentiated
and nondifferentiated (algebraic) ones. Without restricting generality
we will assume that $x_{1},\ldots,x_{k}$ are the differentiated components
while $x_{k+1},\ldots,x_{m}$ are the algebraic ones. Therefore, we
set
\[
D\in\R^{k\times m},\quad D=[I_{k},0]
\]
with the identity matrix $I_{k}\in\R^{k}$. Moreover, we assume $\rank A(t)\equiv k$
for all $t\in[a,b]$.
It is convenient
to rewrite (\ref{eq:DAE}) in the so-called \emph{standard form}
\begin{equation}
E(t)x'(t)+F(t)x(t)=q(t),\quad t\in[a,b]\label{eq:EF}
\end{equation}
where $E=AD$ and $F=B$.
Both formulations are equivalent if a solution of (\ref{eq:DAE})
has also differentiable algebraic components. This is what we assume
in the following.
\medskip

In the present paper we presume regularity in the sense of \cite[Definition 3.3]{HaMa231} as basic term, which is defined via the level-wise reduction from \cite[Sections 12 and 13]{RaRh02}\footnote{For  details we refer also to Section \ref{subs.Alg} below.}. If the DAE is regular with index $\mu\in \Natu$, then owing to \cite[Theorem 4.3]{HaMa231} also the  regular strangeness index $\mu^S$ \cite[Definition 3.15]{KuhMeh06} and the regular tractability index $\mu^T$ \cite[Definition 2.25]{LMT} are well-defined, and $\mu^{S}=\mu-1$, $\mu^{T}=\mu$ . In particular, this allows to combine tools and results of the related concepts.
\medskip

As pointed out in \cite{HaMa231} there are two particular continuously time-varying subspaces which shape the flow structure of a regular DAE, namely the so-called \emph{flow subspace $\Scan$} and its \emph{canonical complement $\Ncan$}. 
More precisely,
the subspace $\Scan(\bar{t})$, $\bar{t}\in[a,b]$, is the set of all
possible function values of solutions to the homogeneous DAE $A(Dx)'+Bx=0$,
that is
\begin{align*}
\Scan(\bar{t})= & \{\bar{x}\in\R^{m}:\text{there is a solution }x:(\bar{t}-\delta,\bar{t}+\delta)\cap[a,b]\rightarrow\R^{m}\\
 & \text{ of the homogeneous DAE such that }x(\bar{t})=\bar{x}\}.
\end{align*}
Thus, $\Scan(\bar t)$ is just the set of consistent initial values for the homogeneous DAE at time  $\bar{t}$. 
While $\Scan$ can be seen as a generalization to the state space  or configuration space for regular ODEs, the second subspace does not occur at all in the ODE theory and first appears  in the DAE context.
The second subspace $\Ncan(\bar{t})$ is a special pointwise complement to the flow subspace 
$\Scan(\bar{t})$,
\[
\R^{m}=\Ncan(\bar{t})\oplus\Scan(\bar{t}),\quad \bar t\in[a,b]
\]
such that each initial value problem
\[
A(t)(Dx(t))'+B(t)x(t)=q(t),\quad t\in(\bar{t}-\delta,\bar{t}+\delta)\cap[a,b],\quad x(\bar{t})-\bar{x}\in\Ncan(\bar{t}),
\]
with arbitrary  sufficiently smooth $q$ and $\bar{x}\in\R^{m}$ is uniquely
solvable without requiring any consistency condition for $q$ and its derivatives at $\bar t$ to be satisfied. These conditions define $\Ncan(\bar{t})$ uniquely.
It is called the \emph{canonical complement} to the flow subspace at $\bar{t}$. 

It should be underlined that both time-varying canonical subspaces of a regular DAE have constant dimension on $[a,b]$. Furthermore, the inclusion $\ker D\subseteq \Ncan(t)$ is valid.

The continuous  projector-valued function $\pPi_{\textrm{can}}:[a,b]\rightarrow\R^{m\times m}$ associated with the above  decomposition of $\R^m$ having the property 
\begin{align*}
 \Scan(t)=\im\pPi_{\textrm{can}}(t),\quad \Ncan(t)=\ker\pPi_{\textrm{can}}(t), \quad t\in [a,b],
\end{align*}
is called \emph{canonical projector function associated with the DAE}. In the context of the tractability index it plays an important role as completely decoupling tool \cite[Chapter 2]{LMT}. 
The general solution of the regular DAE \eqref{eq:DAE} is described as
\begin{align}
 x(t)=X(t,a)x_a+x_q(t),\quad t\in [a,b]   \label{eq:A}
\end{align}
in which $X(\cdot,a):[a,b]\rightarrow \R^{m\times m}$ is the maximal-size fundamental solution normalized at $a$, that means, it is the uniquely determined solution of the matrix IVP
\begin{align*}
 A(DX)'+BX=0, \quad X(a)=\pPi_{\textrm{can}}(a),
\end{align*}
the function  $x_q:[a,b]\rightarrow \R^{m}$ is the unique solution of the IVP
\begin{align*}
 A(Dx)'+Bx=q, \quad x(a)\in N_{\textrm{can}}(a),
\end{align*}
and $x_a\in \R^m$ is the arbitrary free constant that must be fixed in order to capture a specific solution in the flow. Regarding the favorable characteristics
\begin{align}
 \im X(t,a)=\im \pPi_{\textrm{can}}(t)= \Scan(t), \;  \ker X(t,a)=\ker  \pPi_{\textrm{can}}(a)= \Ncan(a), \; t\in [a,b],
 \label{eq:B}
\end{align}
we know that
\begin{align*}
 \pPi_{\textrm{can}}(t)x(t)&=X(t,a)x_a+ \pPi_{\textrm{can}}(t)x_q(t),\\
 (I-\pPi_{\textrm{can}}(t))x(t)&=(I-\pPi_{\textrm{can}}(t))x_q(t)=:v_q(t),\quad t\in [a,b],
\end{align*}
and hence
\begin{align}\label{eq:ICa}
  \pPi_{\textrm{can}}(a)x(a)=\pPi_{\textrm{can}}(a)x_a,\quad 
  (I-\pPi_{\textrm{can}}(a))x(a)=v_q(a).
\end{align}
As shown in \cite[Section 2.6.2]{LMT}, the term $\pPi_{\textrm{can}}x_q$ depends on $q$ itself, but not at all on its derivatives, while certain components of $q$ and all involved derivatives of $q$ are gathered into  the function $v_q$.
We emphasize that the function $v_q$ is invariant of the initial data $a$ and $x_a$.

The component $(I-\pPi_{\textrm{can}}(a))x_a$ does not have any impact here, only the component $\pPi_{\textrm{can}}(a)x_a\in \Scan(a)$ is significant.
Let $l$ denote the constant dimension of the flow subspace $\Scan(t)$. It is evident that $l$ is the number of linearly independent solutions of the homogeneous DAE, that is the dynamical degree of freedom of the DAE.

By construction, one has
\begin{align*}
 \Ncan(t)&\supseteq\ker D=\ker E(t),\\
 \Scan(t)&\subseteq\{z\in\R^m: B(t)z\in\im A(t)D\}=\{z\in\R^m: F(t)z\in\im E(t)\},\\
 l&\leq k=\rank D=\rank E(t),
\end{align*}
but $l$ itself is unknown, and neither the canonical subspaces nor the canonical projector function are given 
\emph{a-priori}.

We close this section by recalling the well-known fact \cite{GM86} that in the case of an index-one DAE it holds that 
\begin{align*}
 \Ncan(t)&=\ker D=\ker E(t),\\
 \Scan(t)&=\{z\in\R^m: B(t)z\in\im A(t)D\}=\{z\in\R^m: F(t)z\in\im E(t)\},\\
 l&=k=\rank D=\rank E(t).
\end{align*}

\subsection{Stating initial conditions and transfer conditions}

As in the previous section we will assume that (\ref{eq:DAE}) is regular with
index $\mu$. 

A value $x_{a}\in\R^{m}$ is called \emph{consistent} if there exists
a solution of (\ref{eq:DAE}) which satisfies the DAE pointwise on $[a,b]$ subject to the initial condition 
\begin{equation}
x(a)=x_{a}.\label{eq:IC}
\end{equation}
As a straightforward
generalization of the situation for explicit ordinary differential
equations one is tempted to consider the initial value problem (\ref{eq:DAE}), (\ref{eq:IC}) and, regarding \eqref{eq:ICa}, a consistent value $x_a$ must have the precise form
\begin{align*}
 x_a=x_a^S+ v_q(a), \quad x_a^S\in \Scan(a),
\end{align*}
otherwise the resulting IVP fails to be solvable.\footnote{We have no intention at all to determine consistent initial values or to develop initialization procedures and refer to the relevant literature, for example:\cite{Pant88,LePeGe91,MaSoe93,Pryce98,ReMaBa00,EstLam21b}.}

In the higher-index case, so-called \emph{hidden constraints} appear, which is accompanied by the appearance of higher derivatives of $q$  in the expression $v_q(a)$ in \eqref{eq:ICa}
 It is just this dependence which makes it
almost always impossible to provide consistent initial values \emph{a-priori}. Furthermore, 
for every given right-hand side $q$,
so-called initialization procedures to compute a consistent value $x_a$ must be executed additionally.

Moreover, even for an index-1 DAE, 
slight perturbations of a consistent initial value $x_{a}\in R^m$ 
may lead to an inconsistent one and the related IVP is no longer solvable. This becomes most obvious for systems
with a special structure, in particular semi-explicit systems. However,
in the latter case it seems easy to cure the situation: Provide only
initial conditions for the differentiated components. In the slightly
more general situation (\ref{eq:DAE}), this corresponds to providing
only initial values for the components $x_{1},\ldots,x_{k}$: For
given $x_{a}^{D}\in\R^{k}$, replace (\ref{eq:IC}) by $Dx(a)=x_{a}^{D}$.
Equivalently, this can be formulated as
\begin{equation}
D(x(a)-x_{a})=0\label{eq:ICP}
\end{equation}
where $x_{a}=D^{+}x_{a}^D$ with the Moore-Penrose inverse $D^{+}$. This works well in the index-one case with $l=k$.
Note that $x_{a}$ in (\ref{eq:ICP}) is only unique up to perturbations
in $\ker D$. The initial value problem (\ref{eq:DAE}), (\ref{eq:ICP})
is well-posed for index-1 DAEs in appropriate  norms. In particular, any small perturbation
of $x_{a}$ implies a small perturbation of the solution of the DAE,
only.\footnote{Such a construction was proposed
in \cite{GM86} in the general index-1 case, also \cite{LMT,LMW}.}

The initial condition (\ref{eq:ICP}) comprises $k=\rank D$ linearly independent equations. This is inappropriate for higher index DAEs featuring merely  $l<k$ free unknowns.
In order to generalize (\ref{eq:ICP}) in an appropriate manner for higher-index problems, it
seems reasonable to use a matrix $G^{a}\in\R^{l\times m}$ with $\ker G^{a}= \Ncan(a)$
and require 
\begin{equation}
G^{a}(x(a)-x_{a})=0.\label{eq:ICG}.
\end{equation}
We have
\begin{proposition}
\label{prop:accuratelyStated}Let $l$ be the degree of freedom of the DAE \eqref{eq:DAE} and the matrix $G^{a}\in\R^{l\times m}$ be such that
$\ker G^{a}=\Ncan(a)=\ker\pPi_{\textrm{can}}(a)$. Then, $G^{a}$ has full
row rank. It holds:
\begin{enumerate}
\item Let $g=G^{a}x_{a}$. Then, the conditions (\ref{eq:ICG}) and
$G^{a}x(a)=g$ are equivalent.
\item The initial value problem (\ref{eq:DAE}), (\ref{eq:ICG}) is uniquely
solvable for each $x_{a}\in\R^{m}$.
\item The solutions of the initial value problem (\ref{eq:DAE}), (\ref{eq:ICG})
depend continuously on $g$ and, thus, also on $x_{a}$.
\end{enumerate}
\end{proposition}
This proposition is a consequence of \cite[Theorems 2.44 and 2.59]{LMT}.
At this place it should be emphasized again that, for indices $\mu>1$, we cannot expect a smooth dependence
of the solutions on perturbations of the right-hand side $q$. Higher-index DAEs are ill-posed in the usual norms.
\footnote{and their 
dynamical degree of freedom $l$
is less than $k$.}
\begin{remark}
\label{rem:accu}
\begin{itemize}
\item Let us be given two matrices $G^a_{1}$ and $G^a_{2}$ that fulfill
the assumptions of Proposition \ref{prop:accuratelyStated}. Let the
solutions $x^{[i]}$ of (\ref{eq:DAE}) fulfill $G^a_{i}(x^{[i]}(a)-x_{a})=0$,
$i=1,2$. Then $x^{[1]}=x^{[2]}$.
\item In \cite{LMW}, an initial condition $G^a x(a)=g$ having the
property 3 of Proposition \ref{prop:accuratelyStated} is called \emph{accurately
stated}.\hfill\qed
\end{itemize}
\end{remark}
\begin{remark}\label{r.derivbound}
If $G^a$ fulfills the requirements of Proposition \ref{prop:accuratelyStated} then the generalized inverse $G^{a-}$ with $G^{a -}G^a=\pPi_{\textrm{can}}(a)$, $G^aG^{a -}=I_l$ is uniquely determined and the formula 
\begin{align*}
 x(t)&=X(t,a)x_a+x_q(t)=X(t,a)\pPi_{\textrm{can}}(a)x_a+x_q(t)\\
 &=X(t,a)G^{a -}G^ax_a+x_q(t)
\end{align*}
records all possible DAE solutions corresponding to the given right-hand side $q$. We fix the reference solution $x_{*}$. Hence, for any solution $x$, it holds
\begin{align*}
 x(t)=x_*(t)+X(t,a)G^{a -}(G^ax_a-G^ax_*(a)), \quad t\in[a,b].
\end{align*}
It becomes evident that, if $X(\cdot,a)$ and $q$ have continuous derivatives, then so are all solutions, and  as far as the derivatives exist it results that
\begin{align*}
 x^{(j)}(t)=x^{(j)}_*(t)+X^{(j)}(t,a)G^{a -}(G^ax_a-G^ax_*(a)), \quad t\in[a,b], \quad j\geq 0.
\end{align*}
Therefore, for the set $U_r$ of all solutions satisfying $|G^ax_a-G^ax_*(a)|\leq r$ there are uniform bounds $C_{j,r}$, such that
\begin{align*}
 \max_{a\leq t\leq b}|x^{(j)}(t)|\leq C_{j,r}, \quad t\in [a,b],\quad x\in U_r, \quad j\geq 0.
\end{align*}
\hfill\qed
\end{remark}
\smallskip

With this background it seems obvious that the most convincing kind
of initial condition is given by (\ref{eq:ICG}) instead of (\ref{eq:ICP})
or even (\ref{eq:IC}). However, for a practical application, the
problem of making a given $x_{a}$ consistent in (\ref{eq:IC}) or
(\ref{eq:ICP}) has been replaced by the problem of computing a matrix
$G^{a}\in\R^{l\times m}$ with $\ker G^{a}=\Ncan(a)=\ker\pPi_{can}(a)$. An advantageous
property of $G^{a}$ is that it depends on the coefficients $A,B$,
only, while a consistent initial value definitely depends also on the right-hand
side $q$ and the provided initialization $x_{a}$. Nevertheless, an accurate
$G^{a}$ is also not easily accessible. 
\medskip

The new least-squares collocation methods for higher-index DAEs that have been under development for several years 
 \cite{HaMaTiWeWu17,HaMaTi19,HaMa21C,HaMaPad121} work with accurate initial and transfer conditions.
To the best of the authors' knowledge, all other procedures for DAEs are tied to consistent values and initializations, respectively.
\bigskip

In order to also state \emph{transfer conditions} that are required in the time-stepping procedure below we now turn to the more general condition for an arbitrary  fixed $\bar t\in [a,b]$,
\begin{align}\label{G.1}
 G(\bar t)x(\bar t)=g,
\end{align}
given by a matrix function $G:[a,b]\rightarrow\R^{l\times m}$ that has full row-rank $l$ as well as the property 
\begin{align}\label{G.2}
\ker  G(t)=\Ncan(t), \quad t\in [a,b].
\end{align}
Such a function is composed in \cite[Lemma 3.1]{HaMaPad121} using  the projector-based analysis and tractability index framework. We are taking a different approach here, which we can better put into practice. Owing to the following theorem, an appropriate function $G$ 
 can be provided, for instance, by means of any basis of the flow subspace associated with the adjoint DAE
\begin{align}\label{G.3}
 -E(t)^*y'(t)+(F(t)^*-E'(t)^*)y(t)=0,\quad t\in[a,b],
\end{align}
and 
\begin{align*}
 -D^*(A^*y)'(t)+B(t)^*y(t)=0,\quad t\in[a,b],
\end{align*}
respectively, by means of the following assertion.

\begin{theorem}
\label{th-Ncan}\cite[Theorem 4]{HaMa231} Let the DAE \eqref{eq:EF}
be\emph{ }regular with index $\mu$ and degree of freedom $l$.
Then, the adjoint DAE 
is also regular with the same characteristics and 
\[
\Ncan(t)=\ker C_{\Scan^{*}}^{\T}(t)E(t),\quad t\in[a,b].
\]
in which  $C_{\Scan^{*}}(t)$ denotes any basis of the flow subspaces $\Scan^{*}(t)$ associated with the adjoint DAE.
\end{theorem}
\begin{corollary}
Owing to Theorem~\ref{th-Ncan}, $G(\bar{t})=C_{\Scan^{*}}^{\T}(\bar{t})E(\bar{t})$
gives rise to accurately stated initial conditions.
\end{corollary}
We underline that there are continuous basis functions  $C_{\Scan^{*}}^{\T}$  such that the resulting function $G$ is continuous, too. The numerical implementation introduced later is defined locally, only, which is why the continuity may be lost there. 

In order to deal with condition \eqref{G.1} it is convenient to slightly generalize the representation \eqref{eq:A}
and use \cite{LMT,LMW}
\begin{align}\label{G.4}
 x(t)=X(t,\bar t)\bar x+\bar x_q(t),\quad t\in[a,b],
\end{align}
in which $X(\cdot,\bar t):[a,b]\rightarrow \R^{m\times m}$ is the maximal-size fundamental solution \emph{normalized at $\bar t$} now, that is, it is the uniquely determined solution of the matrix IVP
\begin{align*}
 A(DX)'+BX=0, \quad X(\bar t)=\pPi_{\textrm{can}}(\bar t),
\end{align*}
the function  $\bar x_q:[a,b]\rightarrow \R^{m}$ is the unique solution of the IVP
\begin{align*}
 A(Dx)'+Bx=q, \quad x(\bar t)\in N_{\textrm{can}}(\bar t).
\end{align*}
$\bar x\in \R^m$ is the arbitrary free constant that must be fixed in order to capture a specific solution in the flow. Similarly to \eqref{eq:B} it holds
\begin{align*}
 \im X(t,\bar t)=\im \pPi_{\textrm{can}}(t)= \Scan(t), \;  \ker X(t,\bar t)=\ker  \pPi_{\textrm{can}}(\bar t)= \Ncan(\bar t), \; t\in [a,b],
\end{align*}
such that, for $t\in[a,b]$,
\begin{align*}
 \pPi_{\textrm{can}}(t)x(t)&=X(t,\bar t)\bar x+ \pPi_{\textrm{can}}(t)\bar x_q(t),\\
 (I-\pPi_{\textrm{can}}(t))x(t)&=(I-\pPi_{\textrm{can}}(t))\bar x_q(t)=v_q(t).
\end{align*}
Recall that the function $v_q$ is not affected by condition \eqref{G.1}.

The basic properties of the function $G$ ensure the existence of a function $G^-:[a,b]\rightarrow \R^{m\times l}$ being a pointwise reflexive generalized inverse such that
\begin{align}
 G(t)G^-(t)=I_l, \quad G^-(t)G(t)=\pPi_{\textrm{can}}(t),\quad t\in[a,b]. \label{Ggi}
\end{align}
The condition \eqref{G.1} can be written equivalently as 
\begin{align}\label{G.5}
 \pPi_{\textrm{can}}(\bar t)x(\bar t)= G(\bar t)^-g.
\end{align}
If $x_*, \tilde x_*:[a,b]\rightarrow \R^m$ are two solutions of the DAE \eqref{eq:EF} corresponding to the initial conditions $G(\bar t)x(\bar t)=g_*$ and $G(\bar t)x(\bar t)=\tilde g_*$, respectively, then the resulting relation 
\begin{align*}
 \tilde x_*(t)-x_*(t)=X(t,\bar t)\pPi_{\textrm{can}}(\bar t)(\tilde x_*(\bar t)-x_*(\bar t))=X(t,\bar t)G^-(\bar t)(\tilde g_*-g_*),\quad t\in[a,b],
\end{align*}
shows the continuity of the solution with respect to $g$.

As mentioned before, in practice the  matrix $G(\bar t)$ is usually not available  \emph{a-priori}\footnote{It should be added that the canonical projector function $\pPi_{\textrm{can}}$ is also described constructively in the context of projector-based analysis \cite{LMT}, but this is generally not yet practically feasible. Even the much simpler explicit representation of the canonical subspace $\Ncan$ in \cite{LMT} is, in general, too complicated due to the derivatives of projector functions involved.}, but it must first be obtained by numerical procedures, and hence, in practice, we are confronted with the perturbed initial condition
\begin{align}\label{G.6}
 \tilde Gx(\bar t)=g,
\end{align}
instead of \eqref{G.1}. How does this affect the solution of the DAE? Recall that \eqref{G.1} yields \eqref{G.5}. In contrast, now we have
\begin{align*}
 G(\bar t)x(\bar t)&=g +(G(\bar t)-\tilde G)x(\bar t)\\
 &=g +(G(\bar t)-\tilde G)\pPi_{\textrm{can}}(\bar t)x(\bar t)   +(G(\bar t)-\tilde G)(I-\pPi_{\textrm{can}}(\bar t))x(\bar t)\\
 &=g +(G(\bar t)-\tilde G)G^-(\bar t)G(\bar t)x(\bar t)   +(G(\bar t)-\tilde G)v_q(\bar t),
\end{align*}
and therefore, with 
\begin{gather}
\mathcal F:= I-(G(\bar t)-\tilde G)G^-(\bar t)= \tilde GG^-(\bar t),\label{G.F}\\
\intertext{it holds}
\mathcal F G(\bar t)x(\bar t)=g + (G(\bar t)-\tilde G)v_q(\bar t).\nonumber
\end{gather}
The question if the initial condition \eqref{G.6} uniquely fixes a DAE solution is now traced back to the question if the matrix $\mathcal F$ is nonsingular. If so, then
\begin{align}\label{G.7}
G(\bar t)x(\bar t)=
\mathcal F^{-1}  g + \mathcal F^{-1}(G(\bar t)-\tilde G)v_q(\bar t)=:\tilde g,
\end{align}
and hence the perturbed initial condition \eqref{G.6} determines a unique DAE solution, too.
The next theorem gives exact information about how the solution of the perturbed problem changes. It turns out that a close approximation $\tilde G$ of $G(\bar t)$ is quite reasonable.
\begin{lemma}\label{l.F}
 Let $\tilde G\in\R^{l\times m}$.
 Then the matrix $\mathcal F$ given in \eqref{G.F} is nonsingular, if and only if the condition
 \begin{align}\label{G.8}
  \ker \tilde G\cap \Scan(\bar t)=\{0\}
 \end{align}
 is satisfied. Moreover, if $\mathcal F$ is nonsingular, $\tilde G$ has full row rank.
\end{lemma}
\begin{svmultproof}
 From $\mathcal Fz=0 $, that is $\tilde GG^-(\bar t)z=0$ and $w=G^-(\bar t)z$ it follows that $w\in \ker \tilde G\cap \Scan(\bar t)$. Now \eqref{G.8} implies $w=0$, thus $z=Gw=0$.
 In the opposite direction, we assume that $\mathcal F= \tilde GG^-(\bar t)$ is nonsingular.
 Consequently, $\tilde G$ has full row rank since obviously $\rank \tilde G \leq l$ and
 \begin{align*}
 l = \rank \mathcal F \leq \min(\rank \tilde G,\rank G^-(\bar t)) \leq \rank\tilde G.
 \end{align*}
 Then,
 $U=G^-(\bar t)(\tilde GG^-(\bar t))^{-1}\tilde G$ represents a projector matrix such that $U^2=U$, $\im U= \Scan(\bar t)$ and $\ker U=\ker \tilde G$, and consequently $\ker \tilde G\oplus \Scan(\bar t)=\R^m$.
\end{svmultproof}
\begin{theorem}\label{t.tildeG}
 Let the DAE \eqref{eq:DAE} be regular with index $\mu\in\Natu$ and dynamical degree of freedom $l>0$, $q:[a,b]\rightarrow\R^m$ be sufficiently smooth.
 
 Let $G:[a,b]\rightarrow\R^{l\times m}$ have the property \eqref{G.2}, $G^-$ be defined by \eqref{Ggi}, and $\tilde G\in\R^{l\times m}$ be an approximation of $G(\bar t)$ such that 
 \begin{align*}
  \rho((G(\bar t)-\tilde G)G^-(\bar t))<1
 \end{align*}
where $\rho(\cdot)$ denotes the spectral radius.
For arbitrarily given $g_*\in\R^l$, let $x_*:[a,b]\rightarrow\R^m$ denote the unique solution of the DAE \eqref{eq:DAE} corresponding to the  initial condition $G(\bar t)x(\bar t)=g_*$.

Then, the matrix $\mathcal F$ given by \eqref{G.F} is nonsingular and the DAE \eqref{eq:DAE} possesses a unique solution $\tilde x_*:[a,b]\rightarrow\R^m$ subject to the perturbed initial condition  $\tilde Gx(\bar t)=g_*$.

Furthermore, with
\begin{align*}
 \tilde g_*:= G(\bar t)\tilde x_*(\bar t)=\mathcal F^{-1}g_* - \mathcal F^{-1}(\tilde G-G(\bar t))v_q(\bar t).
\end{align*}
it holds that
\begin{align*}
 \tilde x_*(t)-x_*(t)&=X(t,\bar t)\pPi_{\textrm{can}}(\bar t)(\tilde x_*(\bar t)-x_*(\bar t))=X(t,\bar t)G^-(\bar t)(\tilde g_*-g_*)\\
 &= X(t,\bar t)G^-(\bar t)\mathcal F^{-1}(G(\bar t)-\tilde G)(G^{-}(\bar t)g_{*}+v_q(\bar t))   ,\quad t\in[a,b].
\end{align*}
In particular, if even $\ker \tilde G=\Ncan(\bar t)$ then $\tilde g_*=\mathcal F^{-1}g_*$ and

\begin{align*}
 \tilde x_*(t)-x_*(t)=X(t,\bar t)(\tilde G^- -G^-(\bar t))g_*,
 \quad t\in[a,b].
\end{align*}
\end{theorem}
\begin{svmultproof}
 First of all, consider the equation $\mathcal Fz=0$, rewritten as $z=(G(\bar t)-\tilde G)G^-(\bar t)z$.
 Owing to the condition on the spectral radius this immediately implies $z=0$, thus $\mathcal F$ is nonsingular. By Lemma \ref{l.F} we may use the projector $U=G^-(\bar t)(\tilde GG^-(\bar t))^{-1}\tilde G$ onto $ \Scan(\bar t)$ along $\ker \tilde G$. Next we choose the generalized inverse $\tilde G^-$ of $\tilde G$ such that $\tilde G^-\tilde G=U$,  $\tilde G\tilde G^-=I_l$. Next, since $\mathcal F G(\bar t)\tilde G^-=\tilde GG^-(\bar t)G(\bar t)\tilde G^-=\tilde GG^-(\bar t)G(\bar t)U\tilde G^-=\tilde GU\tilde G^-=\tilde G\tilde G^-=I_l$ it results that $\mathcal F^{-1}=(\tilde GG^-(\bar t))^{-1}=G(\bar t)\tilde G^-$.
 As we have pointed out before, the solution $\tilde x_*$  is uniquely determined since $\mathcal F$ is nonsingular.
 
From $\tilde G\tilde x_*(\bar t)=g_*$ we obtain 
\begin{align*}
  \underbrace{\tilde G^-\tilde G\pPi_{\textrm{can}}(\bar t)}_{=\pPi_{\textrm{can}}(\bar t)}\tilde x_*(\bar t)&=\tilde G^-g_*-\tilde G^-\tilde G (I-\pPi_{\textrm{can}}(\bar t)) x_*(\bar t)=\tilde G^-g_*-\tilde G^-\tilde G v_q(\bar t),
\end{align*}
and in turn
\begin{align*}
 \tilde g_*&=G(\bar t)\tilde x_*(\bar t)=G(\bar t)\tilde G^-g_*-G(\bar t)\tilde G^-\tilde Gv_q(\bar t)=\mathcal F^{-1}g_*-\mathcal F^{-1}\tilde Gv_q(\bar t)\\
 &=\mathcal F^{-1}g_*-\mathcal F^{-1}(\tilde G-G(\bar t))v_q(\bar t).
\end{align*}
Finally, if even $\ker \tilde G=\Ncan(\bar t)$ then  $(\tilde G-G(\bar t))v_q(\bar t)=0$,  $\tilde g_*=\mathcal F^{-1}g_*=G(\bar t)\tilde G^-g_*$, and $U=\tilde G^-\tilde G$ coincides with $\pPi_{\textrm{can}}(\bar t)$ which leads to
\begin{align*}
 X(t,\bar t)G^-(\bar t)(\tilde g_*-g_*)&=X(t,\bar t)G^-(\bar t)(G(\bar t)\tilde G^- g_*-g_*)\\
 &=X(t,\bar t)(G^-(\bar t)G(\bar t)\tilde G^- g_*-G^-(\bar t)g_*)\\
 &=X(t,\bar t)(\tilde G^- g_*-G^-(\bar t)g_*)=X(t,\bar t)(\tilde G^- -G^-(\bar t))g_*.
\end{align*}
\end{svmultproof}
For further use, we collect some relevant properties in the following lemma.
\begin{lemma}\label{l.FG}
 Let there be given the matrix functions
 \begin{align*}
  G&:[a,b]\rightarrow \R^{l\times m},\; \rank G(t)=l,\; \ker G(t)=\Ncan(t),\; t\in [a,b], \\
  \tilde G&:[a,b]\rightarrow \R^{l\times m}.
  \end{align*}
  Furthermore, let for each $t\in[a,b]$, 
  $G^-(t)\in\R^{m\times l}$ be the reflexive generalized inverse of $G(t)$ satisfying
 \begin{align*}
  G^-(t)G(t)=\pPi_{\textrm{can}}(t),\; G(t)G^-(t)=I_l
  \end{align*}
  and $G^+(t)$ the Moore-Penrose inverse of $G(t)$.
Finally, set
\begin{align*}
  \mathcal F(t)&:=I_l-(G(t)-\tilde G(t))G^-(t), \; t\in [a,b].
 \end{align*}
Let  $G$ be bounded,
$
 |G(t)|\leq c_G \text{ for } t\in [a,b],
$
and satisfy for each $\bar t\in [a,b]$  the condition\footnote{In this context, if the function $G$ is piecewise continuous, then condition \eqref{G.11} is valid.}
\begin{align}\label{G.11}
 \lim_{t\rightarrow\bar t} G(t)w=0\quad \Longrightarrow \; w\in \ker G(\bar t)=\Ncan(\bar t).
\end{align}

\begin{description}
 \item[\rm (1)] Then there are bounds $c_{G^-}$, $c_{G^+}$  such that 
 \begin{align*}
  |G^-(t)|\leq c_{G^-}, \quad  |G^+(t)|\leq c_{G^+}, \quad t\in [a,b].
 \end{align*}
\item[\rm (2)] If $\tilde G$ is close to $G$ in the sense that there is a constant $\chi<1$ with
\begin{align}\label{F}
| G(t)-\tilde G(t)| c_{G^-}\leq \chi<1,\quad\text{ uniformly for all } t\in[a,b],
\end{align}
then $\mathcal F(t)$ remains nonsingular on $[a,b]$, $\rank\tilde G(t) = l$, the decomposition $\ker \tilde G(t)\oplus \Scan(t)=\R^m$ is valid, and $\tilde G^-(t):=G^-(t)\mathcal F(t)^{-1}$ defines  the reflexive generalized inverse of $\tilde G(t)$ featuring $\im \tilde G^-(t)=\Scan(t)$, $\ker \tilde G^-(t)=\{0\}$.
Furthermore,  there is a bound $c_{\tilde G^-}$ such that
\begin{align*}
 |\tilde G^-(t)|&= |G^-(t)\mathcal F(t)^{-1}|\leq c_{\tilde G^-}, \quad t\in [a,b], \\
 |G^-(t)-\tilde G^-(t)|&\leq  c_{G^-} c_{\tilde G^-}|G(t)-\tilde G(t)|,\quad t\in [a,b].
 \end{align*}
\end{description}
\end{lemma}
\begin{svmultproof}
{\rm (1):} The first assertion is trivial if $G$ is continuous, but here $G$ is not necessarily continuous. Recall that $\im G^-(t)=\im\pPi_{\textrm{can}}(t)=\Scan(t)$.
We have to show the existence of a bound $c_{G^-}$ such that
\begin{align*}
 |G^-(t)z|\leq c_{G^-} \quad \text{ for all }\quad z\in\R^l, |z|=1,\; t\in[a,b].
\end{align*}
Suppose, by contradiction, that, to each $n\in\Natu$ there exist $t_n\in[a,b]$ and $z_n\in\R^l$, $|z_n|=1$, such that $|G(t_n)z_n|\geq n$. The resulting sequences reside in compact sets and have accumulation points. For simplicity we assume $t_n\rightarrow t_*\in[a,b]$ and $z_n\rightarrow z_*\in\R^l$.

Denote $w_n=G^-(t_n)z_n\in \Scan(t_n)$ and $\bar w_n=w_n/|w_n|$. Then
\begin{align*}
 |w_n|&\geq n,\\
 G(t_n)w_n&= G(t_n)G^-(t_n)z_n=z_n,\quad |G(t_n)w_n|=1,\\
 |G(t_n)\bar w_n|&\leq 1/n.
\end{align*}
In turn, the sequence of the $\bar w_n\in\Scan(t_n)$, $|\bar w_n|=1$, has an accumulation point and we may suppose $\bar w_n\rightarrow \bar w_*$, $G(t_n)\bar w_n\rightarrow0$.
Regarding the relation
\begin{align*}
 (I-\pPi_{\textrm{can}}(t_*))\bar w_*=\bar w_*-\pPi_{\textrm{can}}(t_*)(\bar w_*-\bar w_n)+(\pPi_{\textrm{can}}(t_*)-\pPi_{\textrm{can}}(t_n))\bar w_n+\bar w_n\rightarrow 0,
\end{align*}
we obtain $\bar w_*\in\Scan(t_*)$. Here we used the continuity of $\pPi_{\textrm{can}}(t)$.
Furthermore, we have also
\begin{align*}
 |G(t_n)\bar w_*|&=|G(t_n)(\bar w_*-\bar w_n)+G(t_n)\bar w_n|\\
 &\leq c_G |\bar w_*-\bar w_n|+|G(t_n)\bar w_n|\rightarrow 0,
\end{align*}
which leads to $\bar w_*\in \Ncan(t_*)$, in turn to  $\bar w_*\in \Ncan(t_*)\cap\Scan(t_*)=\{0\}$. However, this contradicts the property $|\bar w_*|=1$, and hence $G^-(t)$ is uniformly bounded.

By \cite[Theorem 1.20]{NaVo76}, $\lvert G^+(t)\rvert\leq \lvert G^-(t)\rvert \leq c_{G^-}$ and, thus,
$c_{G^+}=c_{G^-}$ is a bound for $\lvert G^+(t)\rvert$.

 {\rm (2):}
 Condition \eqref{F} ensures the nonsingularity of  $\mathcal F(t)=\tilde G(t)G^-(t)$ and $|\mathcal F(t)^{-1}|\leq (1-\chi)^{-1}$, \; $t\in[a,b]$.
 Owing to Lemma \ref{l.F}, $\rank\tilde G(t) = l$ and the decomposition $\ker \tilde G(t)\oplus \Scan(t)=\R^m$ is valid. We have obviously
 \begin{align*}
  |\tilde G^-(t)|=|G^-(t)\mathcal F(t)^{-1}|\leq c_{G^-}(1-\chi)^{-1}=:c_{\tilde G^-}, \quad t\in [a,b].
 \end{align*}
 and the remaining inequality results immediately from
 \begin{align*}
  \tilde G^-(t)-G^-(t)&=G^-(t)\mathcal F(t)^{-1}-G^-(t)=G^-(t)\mathcal F(t)^{-1}(I-\mathcal F(t))\\&=G^-(t)\mathcal F(t)^{-1}(G(t)-\tilde G(t))G^-(t),\quad t\in [a,b].
 \end{align*}
\end{svmultproof}
 Owing to Lemma \ref{l.F} and its proof, we may use the   projector 
 \[
 U(t)=G^-(t)(\tilde G(t)G^-(t))^{-1}\tilde G(t)
 \]
 and apply the generalized inverse $\tilde G(t)^{-}$ so that $\tilde G(t)^{-}\tilde G(t)=U(t)$, $\tilde G(t)\tilde G(t)^{-}=I_l$. This gives $\mathcal F(t)^{-1}=G(t)\tilde G(t)^-$ and 
\begin{align*}
  G^-(t)\mathcal F(t)^{-1}=G^-(t)G(t)\tilde G(t)^-=\pPi_{\textrm{can}}(t)\tilde G(t)^- , \quad t\in [a,b],
 \end{align*}
The resulting inequality
\begin{align}\label{G.10}
 |\tilde x_*(t)-x_*(t)|\leq |X(t,\bar t)| c_{\tilde G^-} |G(\bar t)-\tilde G(\bar t)||\pPi_{\textrm{can}}(\bar t)x_*(\bar t)+v_q(\bar t)|,\quad t, \bar t\in [a,b],
\end{align}
guarantees that, regardless of the point $\bar t$ at which the transfer condition is required, $\tilde x_*$  converges uniformly to $x_*$ as soon as $\tilde G(\bar t)$ tends to $G(\bar t)$.

\section{An algorithm to provide
accurate initial and transfer conditions}\label{s.Algoritm}

\subsection{Basic reduction tools}\label{subs.Alg}
For the problem of computing a suitable matrix $G^{a}$, a solution
has been advocated in \cite{HaMa231} which will be considered subsequently.
Further below, we will also consider its application in an initial
value solver for higher-index DAEs proposed in \cite{HaMaPad121}, in which  transfer conditions from window to window at several times are realized by means of a matrix function $G:[a,b]\rightarrow \R^{l\times m}$  with $\ker G(t)= \Ncan(t)$. 

The key for all subsequent considerations is the reduction procedure
which has been proposed in \cite[Section 12]{RaRh02} and is investigated
in more detail in \cite{HaMa231}. For its description, it is convenient
to write (\ref{eq:DAE}) in standard form and use the coefficient pair $\{E,F\}=\{AD,B\}$.
The reduction procedure consists of the following
steps:
\begin{enumerate}
\item Assume that $E$ has constant rank $r$ on $[a,b]$. Set $E_{0}=E$, $F_{0}=F$,
$m_{0}=m$, and $r_{0}=r$. Let $Z_{0}:[a,b]\rightarrow\R^{m_{0}\times(m_{0}-r_{0})}$
be a function such that $Z_{0}(t)$ is a basis of $(\im E_{0}(t))^{\perp}$
and $Y_{0}:[a,b]\rightarrow\R^{m_{0}\times r_{0}}$ a function such
that $Y_{0}(t)$ is a basis of $\im E_{0}(t)$ for every $t\in[a,b].$
The homogeneous DAE is then equivalent to the system
\[
Y_{0}^{\T}E_{0}x'+Y_{0}^{\T}F_{0}x=0,\quad Z_{0}^{\T}F_{0}x=0.
\]
Here and in the following we omit the argument $t$ for the sake of
simplifying the notation.
\item If $[E_{0},F_{0}]$ has full (row) rank, then $\rank Z_{0}^{\T}F_{0}=m_{0}-r_{0}$.
Hence, the subspace $S_{0}=\ker Z_{0}^{\T}F_{0}$ has dimension $r_{0}$
in this case. So each of the solutions of the homogeneous DAE (\ref{eq:EF})
must stay in $S_{0}$. Assume that $C_{0}:[a,b]\rightarrow\R^{m_{0}\times r_{0}}$
is a differentiable function such that $C_{0}(t)$ is a basis of $S_{0}(t)$,
$t\in[a,b]$. Then $x$ has a representation of the form $x=C_{0}x_{[1]}$
with $x_{[1]}:[a,b]\rightarrow\R^{r_{0}}$ being a smooth function.
Denote $E_{1}=Y_{0}^{\T}E_{0}C_{0}$ and $F_{1}=Y_{0}^{\T}(F_{0}C_{0}+E_{0}C_{0}')$.
Then it holds
\[
E_{1}x_{[1]}'+F_{1}x_{[1]}=0
\]
on $[a,b]$.
\item Note that $E_{1}(t),F_{1}(t)\in\R^{r_{0}\times r_{0}}$. If $r_{0}=m$,
then $Y_{0}(t)\equiv I_{m}$ can be chosen. Moreover, $S_{0}(t)\equiv\R^{m_{0}}$
and $C_{0}(t)=I_{m}$ is appropriate. In turn, $E_{1}=E_{0}$ and
$F_{1}=F_{0}$ such that the process becomes stationary.
\item Let now $0<r_{0}<m$. The steps 1-3 can be repeated with $E_{0},F_{0}$
replaced by $E_{1},F_{1}$ provided that the different assumptions
are fulfilled.
\item If $0=r_{0}<m$, then $E_{0}=0$, and the DAE does not have any dynamical
degrees of freedom. In that case, we can define formally a matrix
$C_{0}(t)\in\R^{0\times m}$which is an empty matrix in the sense
of Matlab.\footnote{Matlab is a registered trademark of The MathWorks Inc.}
\cite{MATLAB}
\end{enumerate}
This process must stop since either one of the necessary assumptions
is not fulfilled or at a stage where $E_{\mu}$ has full rank. In
the latter case, let $\mu$ be the smallest integer with this property.

Recall that by
\cite[Definition 2]{HaMa231} the pair $\{E,F\}$ is defined to be \emph{regular with index $\mu$}
if the just described process is well-defined and stops in step 3 or 5.
The number $\mu$ is  the index of $\{E,F\}$ on $[a,b]$.

The following theorem is of importance. Note that $C_{\mu-1}=I_{r_{[\mu-1]}}$.\footnote{If $r_{\mu-1}=0$, $I_{r_{[\mu-1]}}$ is, by convention, an empty
matrix.}
\begin{theorem}
\label{th-Scan}\cite[Theorem 3.5 and 5.1]{HaMa231} Let the pair
$\{E,F\}$ in (\ref{eq:EF}) (and, equivalently, (\ref{eq:DAE}))
be regular with index $\mu$. 
\begin{enumerate}
\item Then the flow subspace $\Scan(t)$ is an $r_{\mu-1}$-dimensional
subspace of $\R^{m}$ and the dynamical degree of freedom of the DAE is $l=r_{\mu-1}$.
Let $C_{0},C_{1},\ldots,C_{\mu-1}$ be the
matrix-valued functions arising in step 2 of the reduction procedure.
Let $C:[a,b]\rightarrow\R^{m\times r_{\mu-1}}$ be given by $C=C_{0}C_{1}\cdots C_{\mu-1}$.
Then $C(t)$ is a basis of $\Scan(t)$ for all $t\in[a,b].$
\item For each $t\in[a,b]$, there exists the canonical complement $\Ncan(t)$
such that $\R^{m}=\Scan(t)\oplus\Ncan(t)$.
\item The initial condition 
\begin{equation}
G(\bar{t})x(\bar{t})=g\label{eq:AICimp}
\end{equation}
is accurately stated at $\bar{t}\in[a,b]$ for any $g\in\R^{r_{\mu-1}}$
if $G(\bar{t})$ has full row rank and $\ker G(\bar{t})=\Ncan(\bar{t})$.
\end{enumerate}
\end{theorem}
Let us note that (\ref{eq:AICimp}) can be replaced by $G(\bar{t})(x(\bar{t})-x_{\bar{t}})=0$,
$x_{\bar{t}}\in\R^{m}$.

For later discussions, the reduction procedure is sketched in a more
concise form in Algorithm~\ref{alg:Reduction-procedure}.

\begin{algorithm}

\caption{Reduction procedure\label{alg:Reduction-procedure}}

\textbf{global} $\mu$\;
\SetKwProg{Fn}{function}{}{end}
\SetKwFunction{Cbasis}{Cbasis}
\Fn{\Cbasis{$E,F$}} {
    $m = \dim E$\;
    $r = \rank E$\;
    \If{$m = r$}{
        \textbf{return} $I_m$\;}
    $\mu = \mu+1$\;
    Let $Y$ be a basis of $\im E$\;
    Let $Z$ be a basis of $(\im E)^\perp$\;
    Let $C$ be a basis of $\ker Z^\T F$\;
    $E_{\rm new} = Y^\T EC$\;
    $F_{\rm new} = Y^\T(FC+EC')$\;
    \textbf{return} $C\cdot\Cbasis{$E_{\rm new},F_{\rm new}$}$\;
}
\end{algorithm}

In order to determine a matrix $G(\bar{t})$ giving rise to an accurately
stated initial condition one will need related information concerning $\Ncan(\bar{t})$.
Following Theorem \ref{th-Ncan},   this can be accomplished by using the adjoint DAE to (\ref{eq:DAE})
and (\ref{eq:EF}), respectively.

Consequently, Algorithm~\ref{alg:Ncan} is suitable to compute a
matrix appropriate for stating accurate initial conditions.

\begin{algorithm}
\caption{Computation of a matrix for stating accurate initial conditions\label{alg:Ncan}}

\SetKwProg{Fn}{function}{}{end}
\SetKwFunction{Accurate}{AccurateInitialConditions}
\SetKwFunction{Cbasis}{Cbasis}
\Fn{\Accurate{$E,F,\bar{t}$}} {
    $C = \Cbasis{$-E^{\T},F^{\T}-(E^{\T})'$}$\;
    \textbf{return} $C(\bar{t})^{\T}E(\bar{t})$\;
}
\end{algorithm}

\begin{remark}
Algorithm~\ref{alg:Reduction-procedure} allows the bases to be chosen
arbitrarily. A rather natural choice seems to be the use of orthonormal
bases $Y(t)$, $Z(t)$, $C(t)$. We will later see (Section~\ref{sec:Properties})
that this choice is most appropriate in the time stepping approach.

It is worth to mention that Chistyakov \cite{Chist}
and Jansen \cite{Jansen14} proposed their own versions of the reduction
algorithm. In particular, they use special constructions for the basis
function $C$ \cite[Appendix A.3]{HaMa231}, see also Appendix \ref{subsec:Choices}.\hfill\qed
\end{remark}
\hspace{0pt}
\begin{remark}\label{r.array}
The key ingredient
for the computation of $G(\bar t)$ is the recursive determination of
bases $Y,$ $Z$, $C$ and the derivation of $C$. In each step of
the recursion, the dimension $m$ of the matrices $E$ and
$F$ is reduced such that the final dimension is just the number
$l$ of dynamical degrees of freedom of the given DAE (\ref{eq:DAE}).

There is an alternative way to generate a basis of the flow subspace $\Scan$,  if so-called  derivative array functions 
obtained from (\ref{eq:EF})
by $\kappa$-fold differentiation are available. More precisely, $\kappa$-fold differentiation of the DAE \eqref{eq:EF} yields the inflated system
\begin{align*}
 \mathcal E_{[\kappa]}(t) x'_{[\kappa]}(t)+ \mathcal F_{[\kappa]}x(t)=q_{[\kappa]}(t),\quad t\in[a,b],
\end{align*}
 with  coefficient functions $\mathcal{E}_{[\kappa]}:[a,b]\rightarrow\R^{m(\kappa+1)\times m(\kappa+1)}$, $\mathcal{F}_{[\kappa]}:[a,b]\rightarrow\R^{m(\kappa+1)\times m}$,
\begin{align*}
\mathcal{E}_{[\kappa]} & =\left[\begin{array}{cccccc}
E & 0 & 0 & 0 & \cdots & 0\\
E'+F & E & 0 & 0 & \ldots & 0\\
E''+2F' & 2E'+F & E & 0 & \ldots & 0\\
 \vdots &  &  & \ddots &  & \vdots\\
 &  &  &  & E & 0\\
E^{(\kappa)}+\kappa F^{(\kappa-1)} & \cdots &  &  & \kappa E'+F &  E
\end{array}\right],\;
\mathcal{F}_{[\kappa]}  =\begin{bmatrix}
F\\F'\\\vdots\\F^{(\kappa)}\end{bmatrix},
\end{align*}
enlarged right-hand side $q_{[\kappa]}:[a,b]\rightarrow\R^{m(\kappa+1)}$ and unknown $x_{[\kappa]}:[a,b]\rightarrow\R^{m(\kappa+1)}$,
\begin{align*}
 q_{[\kappa]}  =\begin{bmatrix}
q\\q'\\\vdots\\q^{(\kappa)}\end{bmatrix}, \quad
x_{[\kappa]}  =\begin{bmatrix}
x\\x'\\\vdots\\x^{(\kappa)}\end{bmatrix}.
\end{align*}
The so-called  \emph{differentiation index} $\nu$ is, if it exists, the  smallest integer $\kappa$ such that $\mathcal{E}_{\kappa}$ has constant rank, say $r_{[\kappa]}$, and 
is smoothly one-full (cp. \cite{Camp87}). Of course, this requires rank
determinations of $\mathcal{E}_{\kappa}$ the size of which becomes larger with increasing $\kappa$.
If the index is not known,
there is no natural cheep procedure to find it whereas it is a natural
byproduct of the above reduction procedure. 

If the DAE \eqref{eq:EF} is regular with index $\mu\in\Natu$ then the differentiation index $\nu$ is well defined and $\nu=\mu$. Furthermore, then the matrix functions $\mathcal E_{[\kappa]}$ feature constant ranks and the matrix functions $[\mathcal E_{[\kappa]}, \mathcal F_{[\kappa]}]$ have full row-rank. 

With a continuous basis function  $Z_{[\kappa]}:[a,b]\rightarrow \R^{m(\kappa+1)\times(m(\kappa+1)-r_{[\kappa]})}$ of $(\im \mathcal E_{[\kappa]})^{\perp}$ we form the sets
\begin{align*}
 S_{[\kappa]}(t)&=\ker Z^*_{[\kappa]}(t)\mathcal F_{[\kappa]}(t),\\
 C_{[\kappa]}(t)&=\{x\in\R^m:x=z+(Z^*_{[\kappa]}(t)\mathcal F_{[\kappa]}(t))^{+}Z^*_{[\kappa]}(t)q_{[\kappa]}(t), z\in S_{[\kappa]}(t) \},\quad t\in[a,b].
\end{align*}
By construction, the values $x(t)$ of DAE solutions belong to the sets  $C_{[\kappa]}(t)$. It has been pointed out already in \cite[Section 2.4.3]{BCP96}  that  $C_{[\nu]}(t)$ is precisely the set of all consistent initial values at time $t$ associated with the DAE. Moreover one has $S_{[\nu]}=S_{[\nu-1]}=\Scan$. Hence,
\begin{align}\label{G.9}
 \Scan(t)=\ker Z^*_{[\nu-1]}(t)\mathcal F_{[\nu-1]}(t), \quad t\in[a,b].
\end{align}
Owing to Theorem~\ref{th-Ncan}, providing a continuous basis $C_{\Scan^*}$ of the flow subspace $\Scan^*$ associated with the adjoint of the given DAE \eqref{eq:EF} we arrive at
\begin{align*}
 \Ncan(t)=\ker C_{\Scan^*}^*(t)E(t), \quad t\in[a,b].
\end{align*}
The fact that this procedure consists of only one, albeit usually high-dimensional, step can be seen as an advantage because the numerical differentiation of the bases is no longer required. However, this assumes that the coefficient functions $\mathcal E_{[\nu-1]}$ and $\mathcal F_{[\nu-1]}$ are really already available, and that they are without errors.\footnote{In other contexts, the derivatives of the coefficients
$E,F$ and the right-hand side $q$ are be provided by numerical differentiations
(e.g., \cite{LePeGe91}) or automatic differentiation (e.g., \cite{EstLam21b}), but without  error estimates that could be adopted now.}
\hfill\qed
\end{remark}
\begin{remark}\label{r.Hyp}
 In the strangeness index framework, according to \cite[Hypothesis 3.48]{KuhMeh06}, the given regular DAE \eqref{eq:EF} with index $\mu$ is remodeled by means of the (slightly modified) inflated system
 \begin{align*}
 \mathcal E_{[\mu-1]}(t) x'_{[\mu-1]}(t)+ \mathcal F_{[\mu-1]}x(t)=q_{[\mu-1]}(t),\quad t\in[a,b],
\end{align*}
to the regular index-one DAE,
\begin{align*}
 Y^*Ex'+Y^*Fx&=Y^*q,\\
 Z^*_{[\mu-1]}\mathcal F_{[\mu-1]}x&=Z^*_{[\mu-1]}q_{[\mu-1]}.
\end{align*}
For this aim, a continuous basis $C$ of the flow-subspace $\Scan$ is generated which  yields $Z^*_{[\mu-1]}\mathcal F_{[\mu-1]}C=0$. The property $ker E\cap \Scan=\{0\}$ is a consequence of regularity and hence the matrix function $EC$ has full column-rank $l=\dim \Scan$. Then there is a continuous basis $Y$ of $\im EC$ leading to the above new DAE. The index-one criterion $\ker Y^*E\cap \Scan=\{0\}$ is actually satisfied, since $Y^*Ez=0, z=Cw$ implies $Y^*ECw=0$, thus $w=0$, in turn $z=0$. By construction, $\tilde N:=\ker Y^*E$ has dimension $m-l$. A comparison with equation \eqref{G.8} and Theorem \ref{t.tildeG} raises the question of whether this does not also shift the flow. In contrast, if $Y$ is chosen to form a basis of the flow subspace associated to the adjoint DAE, then by Theorem \ref{th-Ncan}, $\ker Y^*E=\ker C^*_{\Scan^*}E=\Ncan$ results.
\hfill\qed
\end{remark}

The key ingredients for the implementation of Algorithms~\ref{alg:Reduction-procedure}
and \ref{alg:Ncan} are the rank determination, the computation of
the bases, and the differentiation. These components will be investigated
in more detail subsequently.

\subsection{Numerical differentiation\label{subsec:Numerical-differentiation}}

Let $\bar{t}\in[a,b]$ be given. In line 13 of Algorithm~\ref{alg:Reduction-procedure}
and line 2 of Algorithm~\ref{alg:Ncan} the differentiation of matrix
valued functions at $\bar{t}$ is required. Hence, they must be known
in a neighborhood of $\bar{t}$. The recursive nature of Algorithm~\ref{alg:Reduction-procedure}
requires that even the (approximation of the) derivatives must be
known in a neighborhood of $\bar{t}$. These requirements suggest the
application of a collocation approach.

Let $\tau>0$ be given with $\tau<b-a$. Choose an interval $[c,c+\tau]\subseteq[a,b]$
such that $\bar{t}\in[c,c+\tau]$. The choice $c=\bar{t}-\tau/2$
will lead to central approximations while $c=\bar{t}$ or $c=\bar{t}-\tau$
will give one-sided approximations. Let $N$ be a positive integer
and let, for any nonnegative integer $K$, $\mathfrak{P}_{K}$ denote
the space of all polynomials of degree less than or equal to $K$.
Let $M=N+1$ and select collocation points
\begin{equation}
c\leq\sigma_{1}<\sigma_{2}<\cdots<\sigma_{M}\leq c+\tau.\label{eq:inodes}
\end{equation}

\begin{definition}
Let $f\in C^{1}([a,b],\R)$. Let $p\in\mathfrak{P}_{N}$ be the polynomial
interpolating $f$ at the nodes $\sigma_{i}$, $i=1,\ldots,M$. The
\emph{spectral derivative} $\Delta_{\tau}^{N}f$ to $f$ on $[c,c+\tau]$
is given by $\Delta_{\tau}^{N}f=p'$. $\Delta_{\tau}^{N}$ is the
\emph{spectral differentiation operator}.
\end{definition}
Spectral differentiation is a well-known construction \cite{Bal00,BalTru03}. 

The interpolation error of polynomial interpolation has the form \cite[Theorem 2.1.4.1]{StBu93}
\[
f(t)-p(t)=w(t)f[\sigma_{1},\ldots,\sigma_{M},t]
\]
where $w(t)=\prod_{i=1}^{M}(t-\sigma_{i})$, and $f[\sigma_{1},\ldots,\sigma_{M},t]$
denotes the $M$-th divided difference. Then we have
\[
f'(t)-\Delta_{\tau}^{N}f(t)=w'(t)f[\sigma_{1},\ldots,\sigma_{M},t]+w(t)f[\sigma_{1},\ldots,\sigma_{M},t,t].
\]
The following property follows immediately.
\begin{proposition}
\label{prop.gen}Let $f\in C^{N+1}([a,b],\R)$. Then it holds
\[
\lVert\Delta_{\tau}^{N}f-f'\rVert_{C([c,c+\tau],\R)}\leq C\tau^{N}.
\]
\end{proposition}
The constant $C$ in the latter estimate depends on the size of the
$(N+1)$-st derivative of $f$ and the node distribution. In particular,
the Chebyshev nodes of the second kind on $[-1,1]$,
\begin{equation}
\tilde{\sigma}_{i}=\cos\left(\frac{M-i}{M-1}\pi\right),\quad i=1,\ldots,M\label{eq:Cheb2}
\end{equation}
shifted to the interval $[c,c+\tau]$ have exceptional properties.
\begin{theorem}
\label{th.cheb}Let the nodes (\ref{eq:inodes}) be the scaled and
shifted Chebyshev nodes of the second kind (\ref{eq:Cheb2}), that
is, $\sigma_{i}=\frac{1}{2}\tau(\tilde{\sigma_{i}}+1)+c$. Let $f\in C^{N+1}([a,b],\R)$.
Then it holds
\[
\lVert\Delta_{\tau}^{N}f-f'\rVert_{C([c,c+\tau],\R)}\leq2^{-2N+1}(2+2\log N)\frac{\|f^{(N+1)}\rVert_{C([a,b],\R)}}{N!}\tau^{N}.
\]
\end{theorem}
\begin{svmultproof}
Let $p$ be the interpolation polynomial of $f$ at the nodes $\sigma_{i}$,
$i=1,\ldots,M$. Define the functions $\text{\ensuremath{\tilde{p}(s)=p(\frac{\tau}{2}(s+1)+c})}$
and $\tilde{f}(s)=f(\frac{\tau}{2}(s+1)+c)$ for $s\in[-1,1]$. Then,
\cite[Theorem 2.1]{Hav80} provides the estimate
\[
\lVert\frac{\dt}{\dt s}(\tilde{p}-\tilde{f})\rVert_{C([-1,1],\R)}\leq(2+2\log N)\dist_{[-1,1]}(\frac{\dt}{\dt s}\tilde{f},\mathfrak{P}_{N-1})
\]
where
\[
\dist_{[-1,1]}(\frac{\dt}{\dt s}\tilde{f},\mathfrak{P}_{N-1})=\inf\{\lVert\frac{\dt}{\dt s}\tilde{f}-u\rVert_{C([-1,1],\R)}:u\in\mathfrak{P}_{N-1}\}.
\]
Since $\frac{\dt}{\dt s}\tilde{f}(s)=\frac{\tau}{2}f'(t)$, 
\[
\dist_{[-1,1]}(\frac{\dt}{\dt s}\tilde{f},\mathfrak{P}_{N-1})=\frac{\tau}{2}\dist_{[c,c+\tau]}(f',\mathfrak{P}_{N-1}).
\]
Let $q\in\mathfrak{P}_{N-1}$ denote the interpolation polynomial of
$f'$ where the nodes $\rho_{i}$ are the Chebyshev nodes of the first
kind, $\tilde{\rho}_{i}=-\cos\left(\frac{2i-1}{2M-2}\pi\right)$,
$i=1,\ldots,M-1$, scaled and shifted to the interval $[c,c+\tau]$,
that is, $\rho_{i}=\frac{\tau}{2}(\tilde{\rho}_{i}+1)+c$. The error
can be estimated by
\[
\lVert f'-q\rVert_{C([c,c+\tau],\R)}\leq\frac{\|f^{(N+1)}\rVert_{C([a,b],\R)}}{N!}\lVert w\|_{C([c,c+\tau],\R)}
\]
where $w(t)=\prod_{i=1}^{M-1}(t-\rho_{i})$. By the special choice
of the interpolation nodes, $w(t)=2^{-N+1}\left(\frac{\tau}{2}\right)^{N}T_{N}(\frac{2}{\tau}(t-c)-1)$
with the Chebyshev polynomial $T_{N}$. Since $\lVert T_{N}\rVert_{C([-1,1],\R)}=1$,
we obtain
\begin{align*}
\frac{\tau}{2}\lVert\Delta_{\tau}^{N}f-f'\rVert_{C([c,c+\tau],\R)} & =\lVert\frac{\dt}{\dt s}(\tilde{p}-\tilde{f})\rVert_{C([-1,1],\R)}\\
 & \leq2^{-N+1}(2+2\log N)\frac{\|f^{(N+1)}\rVert_{C([a,b],\R)}}{N!}\left(\frac{\tau}{2}\right)^{N+1}.
\end{align*}
\end{svmultproof}

\begin{remark}
A further discussion about optimal node distributions for the spectral
differentiation is provided in \cite{Hoa15}. In particular, for odd
$N$, a distribution on $[-1,1]$ is determined such that both boundaries
belong to the set of nodes and the error for the approximate derivative
at the boundaries is smaller than those of the Chebyshev nodes of
the second kind while only being moderately larger inside the interval.
Hence, they are better suited for one-sided approximations, that is
$\bar{t}=c$ or $\bar{t}=c+\tau$.\hfill\qed
\end{remark}
Algorithmically, for given values $f(\sigma_{1}),\ldots,f(\sigma_{M})$,
the values $\Delta_{\tau}^{N}f(\sigma_{1}),\ldots,\Delta_{\tau}^{N}f(\sigma_{M})$
must be computed. Once the latter are known, the derivative $p'=\Delta_{\tau}^{N}f$
can be evaluated by polynomial interpolation.

Let $\mathbf{f}=[f(\sigma_{1}),\ldots,f(\sigma_{M})]^{\T}$
and $\mathbf{\mathbf{f}_{\tau}}=[\Delta_{\tau}^{N}f(\sigma_{1}),\ldots,\Delta_{\tau}^{N}f(\sigma_{M})]^{\T}$.
Since the mapping $\mathbf{f}\mapsto\mathbf{f}_{\tau}$ is linear,
there exists a matrix $\hat{D}_{\tau}^{N}$ (the so-called \emph{spectral
differentiation matrix}) such that $\mathbf{f}_{\tau}=\hat{D}_{\tau}^{N}\mathbf{f}$
\cite{BerTre04}. The investigations in \cite{Bal00,BalTru03} (mostly
for the case of Chebyshev nodes of the second kind) show that a more
accurate version is provided by an expression derived in \cite{SchWer86}\footnote{In the context of spectral differentiation, this is known as the \emph{nullsum
 trick}.}\footnote{The index-value $2$ in $D_{2,ij}^{N}$ stands for the length of the interval
 $[-1,1]$.}
\begin{multline}
\Delta_{\tau}^{N}f(\sigma_{i})=-\frac{2}{\tau}\sum^M_{j=1,j\neq i}\frac{w_{j}}{w_{i}}\frac{f(\sigma_{i})-f(\sigma_{j})}{\tilde{\sigma}_{i}-\tilde{\sigma}_{j}}=-\frac{2}{\tau}\sum_{j=1,j\neq i}^{M}D_{2,ij}^{N}(f(\sigma_{i})-f(\sigma_{j})),\\
i=1,\ldots,M\label{eq:specDiff}
\end{multline}
where $\tilde{\sigma}_{i}$ denotes the scaled nodes on $[-1,1]$,
that is $\tilde{\sigma}_{i}=\frac{2}{\tau}(\sigma_{i}-c)-1$. The
coefficients $w_{i}$ are the weights of the barycentric polynomial
interpolation formula \cite{BerTre04}
\[
w_{i}=\frac{1}{\prod_{j=1,j\neq i}^{M}(\tilde{\sigma}_{i}-\tilde{\sigma}_{j})}.
\]
Note that the weights $w_{i}$, and hence the coefficients $D_{2,ij}^{N}$,
are independent of $[c,c+\tau]$ and can, thus, be precomputed. For
a number of often used nodes these weights are known explicitly \cite{BerTre04,WaHuVa14}.
In particular, for the Chebyshev nodes of the second kind (\ref{eq:Cheb2}),
these weights have a very simple form. Up to a common scaling factor
it holds
\[
w_{i}=(-1)^{M-i}\delta_{i},\quad\delta_{i}=\begin{cases}
1/2, & j=1\text{{\,or\,}}j=M,\\
1, & \text{otherwise}
\end{cases}
\]
Next, assume that perturbations of the data are present, that is, instead
of $f$ just a function $f_{\epsilon}\in C([c,c+\tau],\R)$ with $\lVert f_{\epsilon}-f\rVert_{C([c,c+\tau],\R)}\leq\epsilon$
is available. In the presence of perturbations, we can estimate 
\begin{align}
\lVert\Delta_{\tau}^{N}f_{\epsilon}-f'\rVert_{C([c,c+\tau],\R)} & \leq\lVert\Delta_{\tau}^{N}f_{\epsilon}-\Delta_{\tau}^{N}f\rVert_{C([c,c+\tau],\R)}+\lVert\Delta_{\tau}^{N}f-f'\rVert_{C([c,c+\tau],\R)}\nonumber \\
 & \leq(1+\Lambda_{M})\lvert\mathbf{f}_{\epsilon,\tau}-\mathbf{f}_{\tau}\rvert_{\R^{\infty}}+\lVert\Delta_{\tau}^{N}f-f'\rVert_{C([c,c+\tau],\R)}\nonumber \\
 & \leq(1+\Lambda_{M})\lvert D_{2}^{N}\rvert_{\infty}\lvert\mathbf{f}_{\epsilon}-\mathbf{f}\rvert_{\R^{\infty}}\frac{4}{\tau}+\lVert\Delta_{\tau}^{N}f-f'\rVert_{C([c,c+\tau],\R)}\nonumber \\
 & \leq(1+\Lambda_{M})\lvert D_{2}^{N}\rvert_{\infty}\frac{4\epsilon}{\tau}+\lVert\Delta_{\tau}^{N}f-f'\rVert_{C([c,c+\tau],\R)}.\label{eq:perturbedDiff}
\end{align}
Here, $D_{2}^{N}=(D_{2,ij}^{N})_{i=1,\ldots,M,j=1,\ldots,M}$ with
$D_{2,ii}^{N}=0$. $\lvert\cdot\rvert_{\infty}$ is the row sum norm.
Moreover, $\Lambda_{M}$ is the Lebesgue constant belonging to the
node sequence (\ref{eq:inodes}).\footnote{The size of the Lebesgue constant does not matter later on.} For the nodes (\ref{eq:Cheb2}), an estimate $\Lambda_{M}\leq\alpha\log M$
\cite{Fornberg1994} holds true.
\begin{proposition}
\label{prop:ChebCond}For the Chebyshev nodes of the second kind (\ref{eq:Cheb2})
it holds
\[
\lvert D_{2}^{N}\rvert_{\infty}\leq16N^{2}.
\]
\end{proposition}
The proof is rather technical without contributing much to the topic
of this paper. Therefore, it is moved to Appendix~\ref{subsec:ProofA2}.
\begin{remark}
\label{rem:CondEquidistant}It is well known that equidistant nodes
are inappropriate for interpolation. This is emphasized once more
in the case of numerical differentiation. For equidistant nodes, $\lvert D_{2}^{N}\rvert_{\infty}$
is exponentially growing. For a proof, consult the appendix.\hfill\qed
\end{remark}
Generalizing the idea of spectral collocation one can consider least-squares
approximations by polynomials instead of interpolation. More precisely,
assume $M>N$ nodes (\ref{eq:inodes}) be given. Let $p\in\mathfrak{P}_{N}$
be the polynomial minimizing the functional
\begin{equation}
\Phi^{C}(q)=\sum_{i=1}^{M}(q(\tau_{i})-f(\tau_{i}))^{2}.\label{eq:LSQ}
\end{equation}
In the case of $M=N+1$, this polynomial coincides with the interpolating
polynomial. Setting $\Delta_{\tau}^{N}f=p'$ we obtain an approximation
to the derivative of $f$ also in the more general case. A convergence
proof can be provided by noting that the index-2 DAE
\begin{equation}
\begin{bmatrix}
0 & 0\\
-1 & 0
\end{bmatrix}\begin{bmatrix}
x_{1}'\\
x_{2}'
\end{bmatrix}+\begin{bmatrix}
1 & 0\\
0 & 1
\end{bmatrix}\begin{bmatrix}
x_{1}\\
x_{2}
\end{bmatrix}=\begin{bmatrix}
f\\
0
\end{bmatrix},\quad t\in[c,d].\label{eq:KCF2}
\end{equation}
has the unique solution $x_{1}=f$, $x_{2}=f'$. An application of
the least-squares collocation method of \cite{HaMaT122} amounts to
minimizing the functional
\[
\Phi(q)=\int_{c}^{c+\tau}(q(t)-f(t))^{2}\dt t.
\]
We observe that (\ref{eq:LSQ}) is just a discrete version of this
procedure. Convergence properties for certain modifications of $\Phi^{C}$
are given in \cite{HaMaT122}. So far however, we were only able to
show the convergence of order $O(\tau^{N-1})$ for those modifications. However,
we believe that  $O(\tau^{N})$ convergence can be obtained.
\begin{remark}
Implementation issues
\begin{itemize}
\item For given $\bar{t}\in[a,b]$, a surrounding interval $[c,c+\tau]$
must be chosen. If $\bar{t}=a$ or $\bar{t}=b$, the one-sided choices
$c=a$ and $c=b-\tau$, respectively, are necessary. The error estimates
do not show any preference about the location of $\bar{t}$ inside
$[c,c+\tau]$. However, symmetric formulae show often more accuracy
than unsymmetric ones. This suggests to use $c=\bar{t}-\tau/2$ whenever
possible.
\item In the end of the day, the derivative is only required to be computed
at $\bar{t}$. Therefore, it is most appropriate if $\bar{t}=\sigma_{i}$
for some $i$ in (\ref{eq:inodes}). This will avoid the final interpolation
step. We observed even a much higher accuracy of the derivative approximation
at the nodes. Note that this cannot be explained by the presence of
the Lebesgue constant in (\ref{eq:perturbedDiff}).
\item The theory is most complete for Chebyshev nodes of the second kind.
Since $\pm1$ are nodes, they can be used both for one-sided approximations
and central differentiations (the latter only for odd $M$). In
the one-sided case, a slightly better accuracy seems to be possible
using the nodes proposed in \cite{Hoa15}. However, many other nodes
are applicable, too (Chebyshev nodes of the first kind, Gauss nodes,
Lobatto nodes, Radau nodes, and many more).
\item The case $M>N+1$ does not seem to be suitable with respect to the
approximation quality. However, even in that case we obtained reasonable
results. 
\item The mapping $\mathbf{f}\mapsto\Delta_{\tau}^{N}\mathbf{f}$ is linear
also for $M>N+1$. A representation similar to (\ref{eq:specDiff})
can be derived as follows. Let $T_{0},\ldots,T_{N}$ be the Chebyshev
polynomials and let $\tilde{\sigma}_{i}=\frac{2}{\tau}(\sigma_{i}-c)+1$,
$i=1,\ldots,M$ be the collocation nodes transformed to $[-1,1]$.
Then we form the generalized Vandermonde matrices $V=[T_{j}(\tilde{\sigma}_{i})]_{i=1,\ldots,M,j=0,\ldots,N}$
and $V_{D}=[T_{j}'(\tilde{\sigma}_{i})]_{i=1,\ldots,M,j=0,\ldots,N}$.
For
reasonable collocation nodes, the condition number of $V$ is small. Let $V^+$ denote the
Moore-Penrose inverse of $V$.
Then it holds $\Delta_{\tau}^{N}\mathbf{f}=\frac{2}{\tau}V_{D}V^+\mathbf{f}$.
For a more robust implementation, the null sum trick \cite{BalTru03}
can be applied here, too. Denote $D_{2}^{N}=V_{D}V^+$.
Since the function $f(t)\equiv1$ will be differentiated exactly by
this method, it holds $\sum_{j=1}^{M}D_{2,ij}^{N}=0$ for $i=1,\ldots,M$.
Therefore, 
\[
\Delta_{\tau}^{N}f(\tau_{i})=\frac{2}{\tau}\sum_{j=1,j\neq i}^{M}D_{2,ij}^{N}(f(\sigma_{j})-f(\sigma_{i}))
\]
which closely resembles (\ref{eq:specDiff}).\hfill\qed
\end{itemize}
\end{remark}

\subsection{Rank determination and computation of differentiable bases}

A critical point is the rank determination and the computation of
a smooth basis $C(t)$ in line 11 of Algorithm~\ref{alg:Reduction-procedure}
in a neighborhood of the point $\bar{t}\in[a,b]$. Following the setting
of Section~\ref{subsec:Numerical-differentiation}, let $[c,c+\tau]\subseteq[a,b]$
with $\bar{t}\in[c,c+\tau]$. In order to be able to compute the derivatives,
$C$ must be known at the collocation nodes (\ref{eq:inodes}). The
method of choice for the determination of the ranks is the singular
value decomposition. In uncritical cases, a $QR$-decomposition with
column pivoting might be possible \cite[Chapter 2]{Bjo15}. Unfortunately,
these algorithms in the forms implemented in standard libraries are
unsuitable for our purposes. In the case of the singular value decomposition,
the singular values and, therefore, the bases of the corresponding
images and nullspaces, are ordered in decreasing order. Moreover,
the bases of the nullspaces are not uniquely defined. A similar property
holds in the case of the $QR$-decomposition: The diagonal elements
of the $R$-factor are ordered in decreasing order. The following example
exhibits the behavior.
\begin{example}
\label{Exa:svdQR}
Consider the matrix
\[
H(t)=\begin{bmatrix}
\sin t & 0& 0\\
0 & \cos t& 0
\end{bmatrix},\quad t\in(0,\pi/2).
\]
Let $\hat{t}\in(0,\pi/2)$ denote the unique point at which $\sin\hat{t}=\cos\hat{t}$.
Then it holds $\sin t<\cos t$ for $t\in(0,\hat{t})$ and $\sin t>\cos t$
for $t\in(\hat{t},\pi/2)$. The singular value decompositions and
the $QR$-factorizations have been computed both in Matlab, version
R2022b \cite{MATLAB}, and GNU Octave, version 6.4.0 \cite{Octave}.
Both results are identical. In particular, the singular value decomposition
$H=USV^{\T}$ is
\[
H(t)=\begin{cases}
\begin{bmatrix}
1 & 0\\
0 & 1
\end{bmatrix}\begin{bmatrix}
\sin t & 0 & 0\\
0 & \cos t & 0
\end{bmatrix}\begin{bmatrix}
1 & 0 &0\\
0 & 1 &0\\ 0 & 0& 1
\end{bmatrix}^{\T}, & t\in(\hat{t},\pi/2),\\
\begin{bmatrix}
0 & 1\\
1 & 0
\end{bmatrix}\begin{bmatrix}
\cos t & 0 & 0\\
0 & \sin t & 0
\end{bmatrix}\begin{bmatrix}
0 & 1 & 0\\
1 & 0 & 0\\ 0& 0 & 1
\end{bmatrix}^{\T}, & t\in(0,\hat{t}).
\end{cases}
\]
This leads to 

\[
\ker H(t)=\begin{cases}
\ker \begin{bmatrix}
\sin t & 0 & 0\\
0 & \cos t & 0
\end{bmatrix}\begin{bmatrix}
1 & 0 &0\\
0 & 1 &0\\ 0 & 0& 1
\end{bmatrix}^{\T}= \im \begin{bmatrix}
0 \\ 0\\
1 
\end{bmatrix} , & t\in(\hat{t},\pi/2),\\
\ker \begin{bmatrix}
\cos t & 0 & 0\\
0 & \sin t & 0
\end{bmatrix}\begin{bmatrix}
0 & 1 & 0\\
1 & 0 & 0\\ 0& 0 & 1
\end{bmatrix}^{\T}= \im \begin{bmatrix}
0 \\ 0\\
1 
\end{bmatrix} , & t\in(0,\hat{t}).
\end{cases}
\]
The piecewise continuous function
\[
G(t)=\begin{cases}
\begin{bmatrix}
\sin t & 0 & 0\\
0 & \cos t & 0
\end{bmatrix}\begin{bmatrix}
1 & 0 &0\\
0 & 1 &0\\ 0 & 0& 1
\end{bmatrix}^{\T}, & t\in(\hat{t},\pi/2),\\
\begin{bmatrix}
\cos t & 0 & 0\\
0 & \sin t & 0
\end{bmatrix}\begin{bmatrix}
0 & 1 & 0\\
1 & 0 & 0\\ 0& 0 & 1
\end{bmatrix}^{\T}, & t\in(0,\hat{t})
\end{cases}
\]
stands as an example for the function $G$ in Lemma \ref{l.FG}.

Similarly, the $QR$-factorization with pivoting $HP=QR$ and, equivalently,
$H=QRP^{\T}$ provides
\[
H=\begin{cases}
\begin{bmatrix}
1 & 0\\
0 & 1
\end{bmatrix}\begin{bmatrix}
\sin t & 0 & 0\\
0 & \cos t & 0
\end{bmatrix}\begin{bmatrix}
1 & 0 & 0\\
0 & 1 & 0\\
0 & 0 & 1
\end{bmatrix}^{\T}, & t\in(\hat{t},\pi/2],\\
\begin{bmatrix}
0 & -1\\
-1 & 0
\end{bmatrix}\begin{bmatrix}
-\cos t & 0 & 0\\
0 & -\sin t & 0
\end{bmatrix}\begin{bmatrix}
0 & 1 & 0\\
1 & 0 & 0\\
0 & 0 & 1
\end{bmatrix}^{\T}, & t\in[0,\hat{t}).
\end{cases}
\]
The minus sign in the latter decompositions are a result of the stabilization
procedures in standard implementations as used in Linpack \cite{LINPACK}
and Lapack \cite{LAPACK}. Both Matlab and GNU Octave use different
implementations of Lapack (Intel\footnote{Intel is a registered trademark of Intel Corporation.}
Math Kernel Library \cite{MKL} and OpenBLAS \cite{OpenBLAS}, respectively). Similarly as in the case of the singular value decomposition, a piecewise continuous $G$ is obtained.
\qed
\end{example}
Let us remark in passing that both types of decompositions deliver
bases that are useful in Algorithm~\ref{alg:Reduction-procedure}.
Moreover, these bases are orthonormal bases such that the following
matrix multiplications $Y^{\T}EC$ and $Y^{\T}FC$ as well as $Y^{\T}EC'$
are optimally conditioned. Below we describe strategies for computing
smooth basis functions.

\subsubsection{Singular value decomposition}
\begin{lemma}
\label{L.Pode}Let $S(t)\subseteq\R^{m}$, $t\in[c,c+\tau]$, be a
time varying subspace which has a continuously differentiable basis.
Let $P(t)$ denote the orthoprojection of $\R^{m}$onto $S(t)$, $t\in[c,c+\tau]$
and let, for a fixed $\hat{t}\in[c,c+\tau]$, $\hat{C}$ be an orthonormal
basis of $S(\hat{t})$.

Then, the problem
\begin{equation}
C'=(P'P-PP')C,\quad C(\hat{t})=\hat{C},\label{eq:Code}
\end{equation}
is uniquely solvable and, for each $t\in[c,c+\tau]$, the solution
value $C(t)$ represents an orthonormal basis of $S(t)$.
\end{lemma}
\begin{svmultproof}
With $r=\rank\hat{C}$, we have
\begin{align*}
(C^{\T}C)' & =(C^{\T})'C+C^{\T}C'=C^{\T}(PP'-P'P)C+C^{\T}(P'P-PP')C\equiv0,\\
C^{\T}(\hat{t})C(\hat{t}) & =\hat{C}^{\T}\hat{C}=I_{r}.
\end{align*}
Hence, $C^{\T}(t)C(t)\equiv I_{r}$ for all $t\in[c,c+\tau]$. Moreover,
from
\begin{align*}
((I_{m}-P)C)' & =(I_{m}-P)'C+(I_{m}-P)C'=(I_{m}-P)'+(I_{m}-P)P'PC\\
 & =(I_{m}-P)'C-(I_{m}-P)'PC=(I_{m}-P)'(I_{m}-P)C,\\
(I_{m}-P(\hat{t}))C(\hat{t}) & =(I_{m}-P(\hat{t}))\hat{C}=0,
\end{align*}
 we derive that $C=PC$, and in turn that $\im C(t)=S(t)$, $t\in[c,c+\tau]$.
Therefore, $C(t)$ is an orthonormal basis of $S(t)$.
\end{svmultproof}

A simple consequence of Lemma~\ref{L.Pode} is
\begin{corollary}
\label{C.Code}Under the conditions of Lemma~\ref{L.Pode}, $C'=P'C$,
in particular $C'(\hat{t})=P'(\hat{t})\hat{C}$.
\end{corollary}
\begin{svmultproof}
From the differential equation in the assumptions of Lemma~\ref{L.Pode},
we obtain
\[
C'=P'C+PC'=P'C+P(P'P-PP')PC=P'C.
\]
\end{svmultproof}

Lemma~\ref{L.Pode} and Corollary~\ref{C.Code} suggest the following
construction for obtaining a differentiable basis function as well
as its derivative:
\begin{enumerate}
\item Let a pointwise orthonormal basis $U(t)$ of $S(t)$ be given. Then,
$P(t)=U^{\T}(t)U(t)$ is the orthogonal projection onto $S(t)$.
\item By using the methods of Section~\ref{subsec:Numerical-differentiation},
determine $P'$.
\item Solve the differential equation (\ref{C.Code}). Here, an (overdetermined)
collocation method based on the approximations of Section~\ref{subsec:Numerical-differentiation}
can be used. In the case of $N=M+1$, this amounts to using a classical
spectral collocation method \cite{Fornberg1994}. Note that the choice
$N>M+1$ for solving ordinary differential equations has been proposed
earlier for Chebyshev nodes of the first kind \cite{Albasiny78}.
\end{enumerate}
We use this
procedure for computing a differentiable basis $C(t)$ in Algorithm~\ref{alg:Reduction-procedure}, which is, 
obviously, quite expensive. However, in the time-stepping method below it will be applied for accurate transitions from one window to the next in low frequency only. In this context, we can do without the continuity of  $Z$ while $Y$ is required to be continuous to ensure a continuous subsequent pair  $E_{\textrm{new}}$ and $F_{\textrm{new}}$.\footnote{If, on the other hand, solutions of the inhomogeneous DAE are to be calculated using the reduction method, $Z$ must also be smooth.}

The idea of approximating smooth singular value decompositions by
differential equations is already used, for instance, in \cite{BunByeMehNic91,KunMeh91}.
The procedure of \cite{KunMeh91} assumes the derivative of the matrix
whose smooth basis shall be computed to be explicitly available.
Such a strategy is not applicable in the present case since the reduction
procedure proceeds recursively including new matrices which already
depend on, possibly quite a few, earlier derivatives. The procedure of \cite{BunByeMehNic91}
is much more expensive than the present one and is hard to implement
in the recursive context necessary here.
\begin{remark}
It may well happen that, at some iteration of the recursive procedure,
$E(t)$ has rank 0 in line 4 of Algorithm~\ref{alg:Reduction-procedure}.
This means that the DAE does not have any dynamical degrees of freedom.
Therefore, a strategy must be implemented for a decision if a given
matrix is numerically zero or not.\hfill\qed
\end{remark}

\subsubsection{$QR$-decomposition}

The critical component of the $QR$-decomposition with pivoting is
the choice of the permutation matrix (when pivoting) and the selection
of the sign of the Householder reflection (in stabilizing the computations),
cf. \cite[Section 2.4]{Bjo15} for details. Both problems are of a
discrete nature and, hence, nondifferentiable. A differentiable implementation
can be obtained if these components are fixed. Under the assumption
that $\tau$ is sufficiently small we can assume that a change in
the strategies for pivoting and stabilizations have only a minor influence
on the stability and, thus, on the overall accuracy of the decomposition.
Therefore, the following procedure seems reasonable:
\begin{enumerate}
\item Assume that $E:[c,c+\tau]\rightarrow\R^{m\times m}$ is a matrix valued
function such that $S(t)=\ker E(t)$ is a subspace having a continuously
differentiable basis $C(t)$.
\item For a fixed $\hat{t}\in[c,c+\tau]$, compute the $QR$-decomposition
of $E(\hat{t})$ with pivoting according to a standard algorithm,
e.g. using the strategy of \cite{LINPACK}. Record the pivoting order
and the sign selection when determining the Householder reflection.
\item Find the rank $r$ of $E(\hat{t})$ and determine $C(\hat{t})$ as
an orthonormal basis of $S(\hat{t})$.
\item For $t\neq\hat{t}$, compute a $QR$-decomposition of $E(t)$ by applying
exactly $r$ Householder reflections using the pivoting and stabilization
strategy recorded in step 2.
\end{enumerate}
\begin{remark}
Implementation issues
\begin{itemize}
\item The rank determination in step 3 of the algorithm is critical. In
order to make this decision more reliable, a cheap estimator for the
singular values described in \cite{Bis90} can be used.
\item Step 4 may fail if the rank of $E(t)$ is less than $r$. In that
case, the assumptions of the algorithm are not fulfilled. Similarly
to step 3, an estimator for the singular values is reasonable to use
in order to ensure that $\rank E(t)=r$.
\item In line 4 of Algorithm~\ref{alg:Reduction-procedure} it may well
happen that $r=0$. This corresponds to a DAE which does not have
any dynamical degrees of freedom. Similarly to the SVD approach, a
test for checking if $E(t)$ is numerically the zero matrix must be
included in step 2.\hfill\qed
\end{itemize}
\end{remark}
More reliable versions of $QR$-decompositions are the so-called rank-revealing
$QR$-decompositions \cite{ChaIps94}. They are only slightly more
expensive than the standard $QR$-decompositions with pivoting but
much more reliable in the rank determination and, consequently, in
the computation of a basis $C$. However, these algorithms compute
iteratively a series of permutations such that they are not differentiable
as well. A way out in this case may be the application of the procedure
suggested by Lemma~\ref{L.Pode} or Corollary~\ref{C.Code}.

\section{Properties of the implementations, examples and Conjecture AIC}\label{s.PropAICnum}
\subsection{Properties of the implementations and Conjecture AIC}\label{sec:Properties}
Algorithm~\ref{alg:Reduction-procedure} is implemented using the
previously discussed ingredients. In order to be more specific let,
for a given integer $N>1$ and a sequence of collocation points $-1\leq\tilde{\sigma}_{1}<\cdots<\tilde{\sigma}_{M}\leq1$,
a discrete differentiation operator $\Delta_{\tau}^{N}$ on $[c,c+\tau]\subseteq[a,b]$
with $\tau>0$, $\sigma_{i}=\frac{1}{2}\tau(\tilde{\sigma_{i}}+1)+c$
and $\bar{t}\in[c,c+\tau]$ be defined. The implementation of the
reduction procedure is then sketched in Algorithm~\ref{alg:Reduction-procedure-1}.

\begin{algorithm}
\caption{Discrete reduction procedure\label{alg:Reduction-procedure-1}}

\SetKwProg{Fn}{function}{}{end}
\SetKwFunction{Accurate}{AccurateInitialConditions}
\SetKwFunction{Cbasis}{Cbasis}
\Fn{\Accurate{$E,F,\bar{t},N,\tau$}} {
    $C = \Cbasis{$-E^{\T},F^{\T}-(E^{\T})',\bar{t},N,\tau$}$\;
    \textbf{return} $C^{\T}E(\bar{t})$\;
}
\vspace{2em}
\textbf{global} $\mu$\;
\SetKwProg{Fn}{function}{}{end}
\SetKwFunction{Cbasis}{Cbasis}
\Fn{\Cbasis{$E,F,\bar{t},N,\tau$}} {
    $m = \dim E$\;
    $r = \rank E$\;
    \If{$m = r$}{
        \textbf{return} $I_m$\;}
    $\mu = \mu+1$\;
    Let $Y$ be an \emph{orthonormal} basis of $\im E$\;
    Let $Z$ be an \emph{orthonormal} basis of $(\im E)^\perp$\;
    Let $C$ be an \emph{orthonormal} basis of $\ker Z^TF$\;
    $E_{\rm new} = Y^\T EC$\;
    $F_{\rm new} = Y^\T(FC+E\Delta_\tau^NC)$\;
    \textbf{return} $C(\bar{t})\cdot\Cbasis{$E_{\rm new},F_{\rm new},\bar{t},N,\tau$}$\;
}
\end{algorithm}

Let in the following $G_{\tau}(t)$ be computed by Algorithm~\ref{alg:Reduction-procedure-1}
for $t\in[a,b]$. Moreover, let $G_{0}(t)$ be the matrix computed
by that algorithm but using exact differentiation instead of spectral
differentiation. In any case, we assume that the purely algebraic routines work without errors, i.e. only approximation errors of the numerical differentiation play a role. Consequently, $\ker G_{0}(t)=\ker G(t)=\Ncan(t)$ applies. We emphasize,
however, that $G_0(t)$ is only locally defined by Algorithm~\ref{alg:Reduction-procedure-1} and thus not necessarily continuous on $[a,b]$.

\begin{lemma}
\label{Cor:Gbound}There exists a constant $c_G$ such that $\lvert G_{\tau}(t)\rvert\leq c_{G}$
for all $\tau\geq 0$ and all $t\in[a,b]$.
\end{lemma}
\begin{svmultproof}
We observe that $E$ is continuous and, thus, bounded. Moreover, 
\[
     C(t)=C_{\mu-1}^{\T}(t)C_{\mu-2}^{\T}(t)\cdots C_{1}^{\T}(t)                                                            
\]
with matrices $C_i(t)$ having orthonormal columns. Hence, $\lvert C(t)\rvert=1$. This provides us with $\lvert G_\tau(t)\rvert = \lvert C(t)E(t)\rvert \leq \lvert E(t)\rvert\leq \max_{a\leq t\leq b}|E(t)|=: c_G$.
\end{svmultproof}

\begin{remark}
The bound for the norms of $G_{\tau}(t)$ 
is crucial for the convergence proofs of the stepwise procedure later
on. The important ingredient for the result to hold is the orthonormality
of the bases $C_{i}(t)$. For the algorithms by Chistyakov and Jansen,
such properties cannot be guaranteed as the basis representations
they used shows (cf. Appendix~\ref{subsec:Choices}). \hfill\qed
\end{remark}
The output of Algorithm~\ref{alg:Reduction-procedure-1} is not necessary
a matrix $G_{\tau}(t)$ with $\ker G_{\tau}(t)=\Ncan(t)$ such that
it is not suitable to formulate accurate initial conditions. An appropriate
measure for determining the inaccuracy of $G_{\tau}(t)$ is the \emph{opening}
(or, \emph{gap}) between $\ker G_{\tau}(t)$ and $\Ncan(t)$.

The opening $\Omega(\mathcal{U},\mathcal{V})$ between two subspaces
$\mathcal{U},\mathcal{V}\subseteq\R^{m}$ is defined by
\begin{align*}
\Omega(\mathcal{U},\mathcal{V}) & =\max\{\omega(\mathcal{U},\mathcal{V}),\omega(\mathcal{V},\mathcal{U})\},\\
\omega(\mathcal{U},\mathcal{V}) & =\sup\{d(u,\mathcal{V}):u\in\mathcal{U},\lvert u\rvert=1\}
\end{align*}
whereby $d(u,\mathcal{V})=\inf\{\lvert u-v\rvert:v\in\mathcal{V}\}$
\cite{Kato58}, \cite[Section 39]{AchGlas77}, see also \cite[Section 2.6.3]{GolubVanLoan1989}.
Note that, in general, $\omega(\mathcal{U},\mathcal{V})\neq\omega(\mathcal{V},\mathcal{U})$.
For example, if $\mathcal{U}\subset\mathcal{V}$ and $\dim\mathcal{U}<\dim\mathcal{V}$,
it holds $\omega(\mathcal{U},\mathcal{V})=0$ while $\omega(\mathcal{V},\mathcal{U})=1$.
The definition of the opening is equivalent to the more symmetric
expression
\[
\Omega(\mathcal{U},\mathcal{V})=\rvert P_{\mathcal{U}}-P_{\mathcal{V}}\lvert
\]
where $P_{\mathcal{U}}$ and $P_{\mathcal{V}}$ denote the orthogonal
projections onto $\mathcal{U}$ and $\mathcal{V}$, respectively \cite[Section 39]{AchGlas77}.

The importance of the opening comes from the following consideration:
Let $x_{0}\in\R^{m}$ be given. Moreover, let $x_{\tau}\in\R^{m}$
be such that $G_{\tau}(t)(x_{\tau}-x_{0})=0$ and $P$ the orthogonal
projection of $\R^{m}$ onto $\Ncan(t)$. Then, $P(x_{\tau}-x_{0})$
is the best approximation to $x_{\tau}-x_{0}$ in $\Ncan(t)$. Hence,
the error becomes
\begin{align*}
\lvert(I_{m}-P)(x_{\tau}-x_{0})\rvert & =d(x_{\tau}-x_{0},\Ncan)=d\left(\frac{x_{\tau}-x_{0}}{\lvert x_{\tau}-x_{0}\rvert},\Ncan(t)\right)\lvert x_{\tau}-x_{0}\rvert\\
 & \leq\Omega(\ker G_{\tau}(t),\Ncan(t))\lvert x_{\tau}-x_{0}\rvert.
\end{align*}
So the gap $\Omega(\ker G_{\tau}(t),\Ncan(t))$ provides a measure
for the inaccuracy of the initial value determined by using $G_{\tau}(t)$
instead of accurate initial conditions.
\begin{property}
In finite dimensional spaces, $\Omega(\mathcal{U},\mathcal{V})$ has
the following properties:
\begin{enumerate}
\item $0\leq\Omega(\mathcal{U},\mathcal{V})\leq1$.
\item If $\Omega(\mathcal{U},\mathcal{V})<1$, then $\dim\mathcal{U}=\dim\mathcal{V}$.
\item If $\dim\mathcal{U}=\dim\mathcal{V}$, then $\omega(\mathcal{U},\mathcal{V})=\omega(\mathcal{V},\mathcal{U})$.
\item If $\dim\mathcal{U}\neq\dim\mathcal{V}$, then $\Omega(\mathcal{U},\mathcal{V})=1$.
\end{enumerate}
\end{property}
A practical computable method for the opening of two subspaces is
provided in the next proposition.
\begin{proposition}
\label{prop:Opening}Let $\mathcal{U},\mathcal{V}\subseteq\R^{m}$
be two subspaces with $\dim\mathcal{U}=\dim\mathcal{V}=r$. Furthermore,
let $V\in\R^{m\times(m-r)}$ be an orthonormal basis of $\mathcal{V}^{\perp}$
and $U\in\R^{m\times r}$ be an orthonormal basis of $\mathcal{U}$.
Then it holds
\[
\Omega(\mathcal{U},\mathcal{V})=\lvert V^{\T}U\rvert,
\]
that is, $\Omega(\mathcal{U},\mathcal{V})$ is the largest singular
value of $V^{\T}U$.\footnote{In \cite{GolubVanLoan1989}, a more symmetric version for the computation
of the opening than the one provided here is given. Let $\hat{U}$
and $\hat{V}$ be orthonormal bases of $\im\mathcal{U}$ and $\im\mathcal{V}$,
respectively, then $\Omega(\mathcal{U},\mathcal{V})=\sqrt{1-\rho_{\textrm{min}}^{2}(\hat{U}^{\ast}\hat{V})}$,
where $\rho_{\textrm{min}}(\hat{U}^{\ast}\hat{V})$ denotes the smallest
singular value of $\hat{U}^{\ast}\hat{V}$. However, in practice it
turned out that this formula is more prone to rounding errors than
the one provided.}
\end{proposition}
\begin{svmultproof}
Since $\dim\mathcal{U}=\dim\mathcal{V}$, it holds $\Omega(\mathcal{U},\mathcal{V})=\omega(\mathcal{U},\mathcal{V})$.
Let $P$ be the orthogonal projector onto $\mathcal{V}$. Then it
holds $d(x,\mathcal{V})=\lvert(I-P)x\rvert$. Furthermore,
\[
(I-P)x=\sum_{i=1}^{m-r}(x,v_{i})v_{i}=\sum_{i=1}^{m-r}v_{i}v_{i}^{\T}x=VV^{\T}x
\]
where $V=(v_{1},\ldots,v_{m-r})$. For $x\in\mathcal{U}$, there exists
a $c\in\R^{r}$ such that $x=Uc$, and, for each $c\in\R^{r}$, $x=Uc\in\mathcal{U}$.
It holds $\lvert x\rvert=\lvert c\rvert$. Hence,
\begin{align*}
d(\mathcal{U},\mathcal{V}) & =\sup_{x\in\mathcal{U},\lvert x\rvert=1}d(x,\mathcal{V})=\sup_{\lvert c\rvert=1}\lvert VV^{\T}Uc\rvert=\left(\sup_{\lvert c\rvert=1}c^{\T}U^{\T}VV^{\T}VV^{\T}Uc\right)^{1/2}\\
 & =\left(\sup_{\lvert c\rvert=1}c^{\T}U^{\T}VV^{\T}Uc\right)^{1/2}=\lvert V^{\T}U\rvert.
\end{align*}
The latter expression equals the largest singular value of $\lvert V^{\T}U\rvert$.
\end{svmultproof}

So far, we were unable to establish rigorous convergence bounds for
the computed approximation $G_{\tau}(t)$ in Algorithm~\ref{alg:Reduction-procedure-1}
However, for further convergence estimates in an initial value solver
in Section~\ref{sec:An-initial-value}, some estimations will be
needed. A certain heuristics can be based on \cite[Theorem 1.3]{Hav84}
which provides an estimation of the kind
\[
\lVert\tilde{p}^{(\kappa)}-\tilde{f}^{(\kappa)}\rVert_{C([-1,1],\R)}\leq c_{\kappa}N^{\kappa-1}\dist_{[-1,1]}(\tilde{f}^{(\kappa)},\mathfrak{P}_{N-\kappa}),
\]
 for $2\leq \kappa\leq N$ and $\tilde{f}\in C^{\kappa}([-1,1],\R)$, where
$\tilde{p}$ is the interpolation polynomial of $\tilde{f}$ using
Chebyshev nodes of the second kind on $[-1,1]$. Similarly to the
proof of Theorem~\ref{th.cheb} this provides us with an estimate
\[
\lVert\left(\Delta_{\tau}^{N}\right)^{\kappa}f-f^{(\kappa)}\rVert_{C([c,c+\tau],\R)}\leq c(\kappa,N)\lVert f^{(N+1)}\rVert_{C([a,b],\R)}\tau^{N+1-\kappa}
\]
for $f\in C^{N+1}([a,b],\R)$ and $1\leq \kappa\leq N$ with a constant $c(\kappa,N)$
depending only on $\kappa$ and $N$. The execution of the algorithm includes
$\mu-1$ differentiations where $\mu$ denotes the index of $\{E,F\}$.
So we will expect an error of the kind $\lvert G_{\tau}(t)-G_{0}(t)\rvert\leq C_d\tau^{N+2-\mu}$.
However, in most experiments so far, the subspace $S_{0}(t)$ in
the first iteration of the algorithm is constant such that $C_{0}'\equiv0$
is computed exactly. This way, we expect to gain another power of
$\tau$,
\begin{equation}
\lvert G_{\tau}(t)-G_{0}(t)\rvert\leq C_d\tau^{N+3-\mu}.\label{eq:simplified-error}
\end{equation}
This is what we observed in practice.

\begin{description}
\item [{\textbf{Conjecture}}] \textbf{AIC}\hspace{1cm}
\begin{enumerate}
\item There exists an $\tau_{0}>0$ uniformly in $t\in[a,b]$ such that
$\rank G_{\tau}(t)=l$ equals the number of dynamical degrees of freedom
of the DAE (\ref{eq:EF}) for $\tau\leq\tau_{0}$. In particular,
if there are no dynamical degrees of freedom, $G_{\tau}(t)$ becomes
a matrix of dimension $0\times m$, that means, it vanishes. Hence,
it is exact.
\item Let the number of dynamical degrees of freedom $l$ be positive. Beyond
rounding errors (and its propagation by the stability properties of
the algorithms of linear algebra), the dominating error is the discretization
error of the numerical differentiation. This and the numerical results
presented later suggest that an error estimation of the kind
\[
\lvert G_{\tau}(t)-G_{0}(t)\rvert\leq C_d\tau^{N+2-\mu}
\]
holds true for all $\tau\leq\tau_{0}$ where $N$ is the order of
the polynomials used for the numerical differentiation operator $\Delta_{\tau}^{N}$.
$G_{0}(t)$ is determined by Algorithm \ref{alg:Reduction-procedure-1},
but using exact differentiation instead of the discrete approximation.
$\mu$ is the index of $\{E,F\}$. The constant $C_d$ depends only
on the matrix pair $\{E,F\}$ and its derivatives. In particular,
it does not depend on $t\in[a,b]$. Motivated by Example~\ref{Exa:svdQR}, we assume that
$G_0$ is piecewise continuous.
\end{enumerate}
\end{description}

\begin{corollary}
Let Assumption AIC hold and the number of dynamical degrees of freedom $l$ be positive. Then there exists a constant $c$ such that, for every $t\in[a,b]$,
\[
 \Omega(\ker G_{\tau}(t),\Ncan(t)) \leq c\tau^{N+2-\mu}.
\]
\end{corollary}
\begin{svmultproof}
$P=I-G_{0}(t)^{+}G_{0}(t)$ is the orthogonal projector onto $\ker G_{0}(t)=\Ncan(t)$.
Similarly, $Q_{\tau}=I-G_{\tau}(t)^{+}G_{\tau}(t)$ is the orthogonal
projector onto $\ker G_{\tau}(t)$. Then we have
\begin{align*}
Q_{\tau}-P&=G_{0}(t)^{+}G_{0}(t)-G_{\tau}(t)^{+}G_{\tau}(t)\\
&=G_{0}(t)^{+}(G_{0}(t)-G_{\tau}(t))+(G_{0}(t)^{+}-G_{\tau}(t)^{+})G_{\tau}(t).
\end{align*}
By Lemma~\ref{l.FG}, there is a uniform bound $c_{G^{+}}$ of $G_0(t)^+$.  
As a consequence
of \cite[Satz 8.2.5]{KiSch88} we obtain
\[
\lvert G_{0}(t)^{+}-G_{\tau}(t)^{+}\rvert\leq 2C_d\tau^{N+2-\mu}
\]
for all sufficiently small $\tau$ independent of $t$. Hence,
\begin{align*}
\Omega(\ker G_{\tau}(t),\Ncan(t))&=\lvert Q_{\tau}-P\rvert \\
&\leq c_{G^{+}}C_d\tau^{N+2-\mu}+2c_{G}C_d\tau^{N+2-\mu}\\
&=C_d(c_{G^{+}}+2c_{G})\tau^{N+2-\mu}
\end{align*}
holds true.
\end{svmultproof}

\subsection{Numerical examples\label{sec:Numerical-examples}}

The methods and algorithms have been implemented in C++. This implementation
used the GNU C++ compiler (version 7.5.0) \cite{GCC}, the Eigen library
(version 3.4.0) \cite{eigenweb}, the Intel(R) Math Kernel Library
(version 2023.0.0) , the GNU Scientific Library (version 2.7) \cite{GSL},
and SuiteSparse (version 5.10.1) \cite{SPQR}. The implementation
has been tuned for speed.

We tested the method at a number of examples. In all these tests,
the index as well as the number of dynamical degrees of freedom were
correctly determined. In particular, if the DAE does not have any
dynamical degrees of freedom ($l=0$), the error measure is $\Omega(\ker G_{\tau},\Ncan)=0$.
So it seems that the implementation is rather reliable. Below, we
present some examples which exhibit the typical behavior of the method
and its implementation.

\subsubsection{\label{subsec:LinCaMo}The linearized example by Campbell and Moore}

We consider the linear DAE (\ref{eq:DAE}) with the coefficients
\begin{align*}
A(t) & =\begin{bmatrix}
1\\
 & 1\\
 &  & 1\\
 &  &  & 1\\
 &  &  &  & 1\\
 &  &  &  &  & 1\\
 &  &  &  &  & 0
\end{bmatrix},\quad D(t)=\begin{bmatrix}
1\\
 & 1\\
 &  & 1\\
 &  &  & 1\\
 &  &  &  & 1\\
 &  &  &  &  & 1 & 0
\end{bmatrix},\\
B(t) & =\begin{bmatrix}
0 & 0 & 0 & -1 & 0 & 0 & 0\\
0 & 0 & 0 & 0 & -1 & 0 & 0\\
0 & 0 & 0 & 0 & 0 & -1 & 0\\
0 & 0 & \sin t & 0 & 1 & -\cos t & -2\rho\cos^{2}t\\
0 & 0 & -\cos t & -1 & 0 & -\sin t & -2\rho\sin t\cos t\\
0 & 0 & 1 & 0 & 0 & 0 & 2\rho\sin t\\
2\rho\cos^{2}t & 2\rho\sin t\cos t & -2\rho\sin t & 0 & 0 & 0 & 0
\end{bmatrix},\quad\rho\neq0.
\end{align*}
The structure of this example has been investigated in detail in \cite{HaMa231}
and used for numerical tests in \cite{HaMaPad121,HaMaTi19}. It is
the linearized version of a test problem proposed in \cite{CaMo95}.
Accurately stated initial conditions are provided by
\[
G(t)=\begin{bmatrix}
H(t) & 0_{2\times3} & 0_{2\times1}\\
H(t)(\mathfrak{A(t)}+\Omega'(t))\Omega(t) & H(t) & 0_{2\times1}
\end{bmatrix}
\]
for any $t\in\R$, where
\begin{align*}
H(t) & =\begin{bmatrix}
\sin t & -\cos t & 0\\
0 & 1 & \cos t
\end{bmatrix},\quad\mathfrak{A}(t)=\begin{bmatrix}
0 & 1 & -\cos t\\
-1 & 0 & -\sin t\\
0 & 0 & 0
\end{bmatrix},\\
\Omega(t) & =\begin{bmatrix}
\cos^{4}t & \sin t\cos^{3}t & -\sin t\cos^{2}\\
\sin t\cos^{3}t & \sin^{2}t\cos^{2}t & -\sin^{3}t\cos t\\
-\sin t\cos^{2}t & -\sin^{3}t\cos t & \sin^{2}t
\end{bmatrix}.
\end{align*}
Here, $0_{i\times j}$ denotes a zero matrix of dimension $i\times j$.
In the following computations we used $\rho=5$.

In all subsequent experiments, we use $\bar{t}=0$ for evaluating
$G_{\tau}(\bar{t})$. In a central approximation, the chosen interval
becomes $[c,c+\tau]=[-\tau/2,\tau/2]$ while for the one-sided approximation,
$[c,c+\tau]=[0,\tau]$ is used. In order to have $\bar{t}=0$ as a
collocation point in the central approximation, the number of collocation
points $M$ must be odd. So the order of the approximating polynomial
becomes $N=M-1$ (spectral differentiation) and $N=M-2$ (least-squares
approximation), respectively. Hence, we expect the convergence orders
$\tau^{M-1}$ and $\tau^{M-2}$, respectively.

Table~\ref{tab:LCMCheb} shows the errors by using Chebyshev nodes
of the second kind (\ref{eq:Cheb2}) and $QR$-factorization. The
first observation is that with rather small $M$ and not too small
$\tau$ the error is in the order of magnitude of the rounding unit.
This indicates that rounding errors are negligible in this setting.
The order of convergence for the spectral differentiation is as expected
$\tau^{M-1}$. Rather surprising is the fact that even for the least
squares approximation the convergence order is $\tau^{M-1}$ which
is higher than the theoretically expected order. Moreover, for low
$M$ and larger $\tau$, the errors between both methods are identical.
We do not have an explanation for this behavior.\footnote{In fact, the numerical computations for this example made us suspicious
that the matrix used for accurate initial condition in \cite[Equation (6.2)]{HaMaPad121}
may be wrong. The considerations in \cite{HaMa231} showed indeed
that that matrix was erroneous.}

\begin{table}
\caption{Opening of the nullspace of $G_{\tau}(0)$ and $\Ncan(0)$ using Chebyshev
nodes of the second kind and using the central approximation on $[-\tau/2,\tau/2]$.
The table shows the results for the polynomial degree $N=M-1$ (spectral
differentiation, above) and $N=M-2$ (differentiation by least squares
collocation, below). Equal data are marked in boldface\label{tab:LCMCheb}}

\begin{centering}
\begin{tabular}{|c|c|c|c|c|c|}
\hline 
$M$\textbackslash{} $\tau$ & 0.1 & 0.05 & 0.25 & 0.125 & 0.00625\tabularnewline
\hline 
\hline 
3 & \textbf{3.29e-03} & \textbf{8.22e-04} & \textbf{2.05e-04} & \textbf{5.14e-05} & \textbf{1.28e-05}\tabularnewline
\hline 
5 & \textbf{2.62e-06} & \textbf{1.64e-07} & \textbf{1.03e-08} & \textbf{6.41e-10} & \textbf{4.01e-11}\tabularnewline
\hline 
7 & \textbf{8.69e-10} & \textbf{1.36e-11} & \textbf{2.12e-13} & 3.23e-15 & 2.97e-16\tabularnewline
\hline 
9 & 1.57e-13 & 1.29e-15 & 5.40e-16 & 9.10e-16 & 3.04e-16\tabularnewline
\hline 
11 & 6.90e-16 & 3.27e-16 & 1.09e-15 & 2.25e-16 & 5.12e-16\tabularnewline
\hline 
\end{tabular}
\par\end{centering}
\bigskip
\centering{}%
\begin{tabular}{|c|c|c|c|c|c|}
\hline 
$M$\textbackslash{} $\tau$ & 0.1 & 0.05 & 0.25 & 0.125 & 0.00625\tabularnewline
\hline 
\hline 
3 & \textbf{3.29e-03} & \textbf{8.22e-04} & \textbf{2.05e-04} & \textbf{5.14e-05} & \textbf{1.28e-05}\tabularnewline
\hline 
5 & \textbf{2.62e-06} & \textbf{1.64e-07} & \textbf{1.03e-08} & \textbf{6.41e-10} & \textbf{4.01e-11}\tabularnewline
\hline 
7 & \textbf{8.69e-10} & \textbf{1.36e-11} & \textbf{2.12e-13} & 3.18e-15 & 2.71e-16\tabularnewline
\hline 
9 & 1.56e-13 & 2.54e-15 & 1.31e-15 & 1.073e-15 & 1.06e-15\tabularnewline
\hline 
11 & 1.60e-15 & 1.16e-15 & 2.26e-15 & 1.12e-15 & 1.26e-15\tabularnewline
\hline 
\end{tabular}
\end{table}

Table~\ref{tab:LCMCheb-1} presents the error for a one-sided approximation
using Chebyshev nodes of the second kind and $QR$-factorization.
In contrast to the first experiment, the convergence orders are as
expected. Rounding errors start to dominate the accuracy at errors
in the order of $10^{3}\epsilon_{\textrm{mach}}$ with the rounding
unit $\epsilon_{\textrm{mach}}$ of the machine used. In the present
example, it seems to be better to use higher order polynomials with
larger stepsizes than lower order polynomials with small stepsizes.

\begin{table}
\caption{Experiments under the same conditions as in Table \ref{tab:LCMCheb}.
However, the interval for the approximation is $[0,\tau]$\label{tab:LCMCheb-1}}

\begin{centering}
\begin{tabular}{|c|c|c|c|c|c|}
\hline 
$M$\textbackslash{} $\tau$ & 0.1 & 0.05 & 0.25 & 0.125 & 0.00625\tabularnewline
\hline 
\hline 
3 & 6.79e-03 & 1.67e-03 & 4.15e-04 & 1.03e-04 & 2.57e-05\tabularnewline
\hline 
5 & 5.39e-06 & 3.33e-07 & 2.07e-08 & 1.29e-09 & 8.04e-11\tabularnewline
\hline 
7 & 1.76e-09 & 2.74e-11 & 4.17e-13 & 8.51e-14 & 1.83e-13\tabularnewline
\hline 
9 & 3.08e-13 & 1.06e-14 & 1.27e-13 & 9.52e-14 & 5.61e-13\tabularnewline
\hline 
11 & 4.62e-14 & 2.36e-14 & 5.21e-14 & 3.93-13 & 6.560e-13\tabularnewline
\hline 
\end{tabular}
\par\end{centering}
\bigskip
\centering{}%
\begin{tabular}{|c|c|c|c|c|c|}
\hline 
$M$\textbackslash{} $\tau$ & 0.1 & 0.05 & 0.25 & 0.125 & 0.00625\tabularnewline
\hline 
\hline 
3 & 1.05e-01 & 5.11e-02 & 2.53e-02 & 1.26e-02 & 6.27e-03\tabularnewline
\hline 
5 & 1.82e-04 & 2.34e-05 & 2.97e-06 & 3.75e-07 & 4.70e-08\tabularnewline
\hline 
7 & 1.16e-07 & 3.66-09 & 1.15e-10 & 3.47e-12 & 1.34e-13\tabularnewline
\hline 
9 & 3.51e-11 & 2.10e-13 & 1.84e-14 & 3.71e-14 & 6.31e-13\tabularnewline
\hline 
11 & 8.79e-14 & 4.84e-14 & 1.23e-13 & 7.80e-13 & 9.00e-13\tabularnewline
\hline 
\end{tabular}
\end{table}

In the next experiment (Table~\ref{tab:LCMRad}), the setting of
Table~\ref{tab:LCMCheb-1} is used with the Chebyshev nodes replaced
by the Radau nodes. The results of both experiments are comparable
while those with Radau nodes provide slightly better approximations.

\begin{table}
\caption{Experiments under the same conditions as in Table \ref{tab:LCMCheb}.
However, the interval for the approximation is $[0,\tau]$ and Radau
nodes are used\label{tab:LCMRad}}

\begin{centering}
\begin{tabular}{|c|c|c|c|c|c|}
\hline 
$M$\textbackslash{} $\tau$ & 0.1 & 0.05 & 0.25 & 0.125 & 0.00625\tabularnewline
\hline 
\hline 
3 & 4.05e-03 & 1.00e-03 & 2.49e-04 & 6.19e-05 & 1.54e-05\tabularnewline
\hline 
5 & 3.41e-06 & 2.11e-07 & 1.31e-08 & 8.18e-10 & 5.10e-11\tabularnewline
\hline 
7 & 1.22e-09 & 1.91e-11 & 2.91e-13 & 1.74e-13 & 6.69e-13\tabularnewline
\hline 
9 & 2.48e-13 & 4.41e-14 & 5.04e-14 & 2.82e-14 & 2.94e-13\tabularnewline
\hline 
11 & 2.36e-14 & 3.49e-14 & 3.01e-14 & 1.24e-12 & 6.64e-13\tabularnewline
\hline 
\end{tabular}
\par\end{centering}
\bigskip
\centering{}%
\begin{tabular}{|c|c|c|c|c|c|}
\hline 
$M$\textbackslash{} $\tau$ & 0.1 & 0.05 & 0.25 & 0.125 & 0.00625\tabularnewline
\hline 
\hline 
3 & 9.01e-02 & 4.42e-02 & 2.19e-02 & 1.09e-02 & 5.43e-03\tabularnewline
\hline 
5 & 1.51e-04 & 1.93e-05 & 2.45e-06 & 3.09e-07 & 3.87e-08\tabularnewline
\hline 
7 & 9.15e-08 & 2.89e-09 & 9.12e-11 & 3.010e-12 & 1.15e-13\tabularnewline
\hline 
9 & 2.70e-11 & 2.14e-13 & 1.35e-14 & 1.38e-14 & 3.08e-13\tabularnewline
\hline 
11 & 3.08e-14 & 5.45e-14 & 1.43e-13 & 1.21e-12 & 8.51e-13\tabularnewline
\hline 
\end{tabular}
\end{table}

Finally, we tested the sensitivity with respect to the choice of the
orthogonalization procedure. Under the conditions of Table~\ref{tab:LCMRad},
the $QR$-factorization has been replaced by an SVD and the procedure
suggested by Corollary~\ref{C.Code}. Thus, this procedure is much
more expensive. The results are presented in Table~\ref{tab:LCMRad-1}.
In the case of the spectral differentiation, the results are in general
slightly better than before. However, for the least-squares approximation,
we see that $M=3$ (meaning $N=1$) is a too low order of approximation.
More surprisingly, the order of convergence for $M>3$ is one higher
than before and larger than expected!

\begin{table}
\caption{This experiment is identical to that in Table~\ref{tab:LCMRad} with the
exception of the $QR$-decomposition replaced by the singular value
decomposition\label{tab:LCMRad-1}}

\begin{centering}
\begin{tabular}{|c|c|c|c|c|c|}
\hline 
$M$\textbackslash{} $\tau$ & 0.1 & 0.05 & 0.25 & 0.125 & 0.00625\tabularnewline
\hline 
\hline 
3 & 3.94e-03 & 9.86e-04 & 2.47e-04 & 1.64e-05 & 1.54e-05\tabularnewline
\hline 
5 & 1.04e-05 & 6.53e-07 & 4.09e-08 & 2.56e-09 & 1.60e-10\tabularnewline
\hline 
7 & 6.85e-09 & 1.08e-10 & 1.80e-12 & 3.35e-13 & 5.63e-13\tabularnewline
\hline 
9 & 2.17e-12 & 9.00e-14 & 1.42e-13 & 1.88e-14 & 7.98e-14\tabularnewline
\hline 
11 & 7.58e-14 & 8.37e-14 & 3.14e-13 & 7.28e-13 & 9.25e-13\tabularnewline
\hline 
\end{tabular}
\par\end{centering}
\bigskip
\centering{}%
\begin{tabular}{|c|c|c|c|c|c|}
\hline 
$M$\textbackslash{} $\tau$ & 0.1 & 0.05 & 0.25 & 0.125 & 0.00625\tabularnewline
\hline 
\hline 
3 & 5.15e-01 & 5.04e-01 & 5.01e-01 & 5.00e-01 & 5.00e-01\tabularnewline
\hline 
5 & 4.20e-06 & 2.64e-07 & 1.65e-08 & 1.03e-09 & 6.45e-11\tabularnewline
\hline 
7 & 1.28e-08 & 2.00e-10 & 3.175e-12 & 9.67e-14 & 7.69e-15\tabularnewline
\hline 
9 & 8.28e-12 & 5.43e-14 & 4.78e-14 & 9.22e-14 & 2.54e-13\tabularnewline
\hline 
11 & 4.26e-14 & 5.21e-14 & 1.85e-13 & 2.78e-13 & 1.63e-12\tabularnewline
\hline 
\end{tabular}
\end{table}

\subsubsection{An example by Chua and Riaza}

In \cite[Section 2.2.6]{Riaza08}, a simple circuit with current controlled
resistors proposed in \cite{ChWuHuZh93} is thoroughly investigated
in order to demonstrate the index analysis for electrical networks.
This example has the nice property that, depending on the behavior of the incorporated
circuit elements, the modeling DAE exhibits different indices and dynamical degrees of freedom.

Consider
\begin{align*}
A(t) & =\begin{bmatrix}
C_{1}(t) & 0 & 0\\
0 & C_{2}(t) & 0\\
0 & 0 & L(t)\\
0 & 0 & 0\\
0 & 0 & 0
\end{bmatrix},\quad D(t)=\begin{bmatrix}
1 & 0 & 0 & 0 & 0\\
0 & 1 & 0 & 0 & 0\\
0 & 0 & 1 & 0 & 0
\end{bmatrix},\\
B(t) & =\begin{bmatrix}
C_{1}'(t) & 0 & 0 & -1 & 1\\
0 & C_{2}'(t) & 1 & 1 & 0\\
0 & -1 & L'(t) & 0 & 0\\
-1 & 1 & 0 & -R_{1}(t) & 0\\
1 & 0 & 0 & 0 & -R_{2}(t)
\end{bmatrix},\quad t\in[a,b].
\end{align*}
 In contrast to the previous example, the derivative of $E^{\T}$
in step 2 of Algorithm~\ref{alg:Ncan} does not vanish identically.
This system has
\begin{itemize}
\item index 1 and $l=3$ if none of the function $C_{1},C_{2},R_{1},R_{2},L$ have
a zero on $[a,b]$;
\item index 2 and $l=2$ if $R_{1}(t)\equiv0$ on $[a,b]$ and none of the functions
$C_{1},C_{2},R_{2},L$ has a zero on $[a,b]$;
\item index 3 and $l=1$ if $R_{1}(t)\equiv0$, $C_{1}(t)+C_{2}(t)\equiv0$ on $[a,b]$
and none of the functions $C_{1},C_{2},R_{2},L$ has a zero on $[a,b]$.
\end{itemize}
Appropriate accurate
initial conditions at $t\in[a,b]$ are given by
\begin{itemize}
\item in the index 1 case: $G(t)=D(t)$;
\item in the index 2 case: $G(t)=\begin{bmatrix}
C_{1}(t)/C_{2}(t) & 1 & 0 & 0 & 0\\
0 & 0 & 1 & 0 & 0
\end{bmatrix}$;
\item in the index 3 case: $G(t)=\begin{bmatrix}
-1, & 1, & -L(t)/(R_{2}(t)C_{1}(t), & 0, & 0\end{bmatrix}$.
\end{itemize}
For the following experiments, we choose
\begin{itemize}
\item in the index 1 case: $C_{1}(t)=\sin t+2$, $C_{2}(t)=\cos t+2$, $L(t)=t^{2}+1$,
$R_{1}(t)=\frac{1}{2}\sin(2t)+1$, $R_{2}(t)=\sin t+\cos t+2$.
\item in the index 2 case: $C_{1}(t)=\sin t+2$, $C_{2}(t)=\cos t+2$, $L(t)=t^{2}+1$,
$R_{1}(t)\equiv0$, $R_{2}(t)=\sin t+\cos t+2$.
\item in the index 3 case: $C_{1}(t)=\sin t+2$, $C_{2}(t)=-C_{1}(t)$,
$L(t)=t^{2}+1$, $R_{1}(t)\equiv0$, $R_{2}(t)=\sin t+\cos t+2$.
\end{itemize}
This example is well-suited to test if the rank decision strategies
work as expected. Perturbations of the arising matrices may lead to
a rank reduction for the index 2 and 3 cases.\footnote{The structural index has been determined to be 2 in both the index 2 and the index 3 cases by the DAESA software \cite{DAESA}.}
However, even for the coarsest grids ($h=0.5$) and lowest interpolation
orders ($N=1$), the gap between $\Ncan(0)$ and $G_{\tau}(0)$ was
in the order of magnitude of the rounding unit for all index cases.

\section{An initial value solver based on accurate initial conditions\label{sec:An-initial-value}}

In \cite{HaMaPad121}, we developed a numerical method for the solution
of initial value problems for higher-index DAEs based on a least-squares
approach. This method is pleasantly simple since it does not need
any preprocessing of the given DAE nor any index reduction procedures.
In the present context it is important to note that this method is
based solely on accurate initial conditions and does not need nor
computes consistent initial values. This method works surprisingly
well. However, in \cite{HaMaPad121} we were able to show convergence
only under the condition that, for any $\bar{t}\in[a,b]$, a matrix
$G(\bar{t})$ is already available which can be used to precisely provide an accurate
initial condition $G(\bar{t})x(\bar{t})=g$. The algorithm now proposed
in Section~\ref{sec:AIC} is just capable to provide reasonable approximations
to such matrices.

To be more precise, consider the initial value problem for (\ref{eq:DAE})
subject to the accurate initial condition
\begin{equation}
G^{a}x(a)=g.\label{eq:Gar}
\end{equation}
It will be assumed that $G^{a}$ belongs to the data of the problem.
Algorithm~(\ref{alg:Reduction-procedure-1}) can be used to check
if $G^{a}$ gives indeed rise to an accurate initial condition.
\medskip

Despite the fact that we are looking for smooth solutions below, the convergence statements for the least-squares method are naturally
formulated in terms of Hilbert space norms that will be derived from Sobolev space norms. Let $\bar{t}\in[a,b)$
and an $H>0$ be fixed such that $\mathcal{I}=[\bar{t},\bar{t}+H]\subseteq[a,b]$.
On the subinterval $[\bar{t},\bar{t}+H]$, we are searching for solutions $x\in H_{D}^{1}(\mathcal{I})$, whereby
\[
H_{D}^{1}(\mathcal{I})=\{x\in L^{2}(\mathcal{I},\R^{m}):Dx\in H^{1}(\mathcal{I},\R^{m})\},
\]
such that (\ref{eq:DAE}) holds\footnote{Later on actually pointwise, since then the right-hand side function $q$ is supposed to be at least continuous. } on $\mathcal{I}$ for given $q\in L^{2}(\mathcal{I},\R^{m})$
and 
\begin{equation}
G_{0}(\bar{t})x(\bar{t})=g.\label{eq:G0bar}
\end{equation}
Here, $G_{0}(t)$ is computed by Algorithm~\ref{alg:Reduction-procedure-1}
but using exact differentiation instead of a discrete approximation.
In particular, this means that (\ref{eq:G0bar}) is an accurate initial
condition.

The norm in $H_{D}^{1}(\mathcal{I})$ is given by
\[
\lVert x\rVert_{H_{D}^{1}(\mathcal{I})}=\left(\lVert x\rVert_{L^{2}(\mathcal{I},\R^{m})}^{2}+\lVert Dx\rVert_{L^{2}(\mathcal{I},\R^{k})}^{2}\right)^{1/2},x\in H_{D}^{1}.
\]
Note that below the following inequality ensured by \cite[Lemma 3.2]{HaMaPad121} plays its role:
\begin{align}\label{RelDx(t)}
 |Dx(t)|\leq C_{H} \lVert x\rVert_{H_{D}^{1}(\mathcal{I})}, \quad t\in[\bar t,\bar t+H], \quad x\in H_{D}^{1}(\mathcal{I}),
\end{align}
with the constant $C_{H}= (2\max\{H,H^{-1}\})^{1/2}\geq \sqrt{2}$ depending on the interval length $H$.

We formulate some general assumptions to be used throughout  the following discussion.
\begin{description}
\item [{\textbf{General}}] \textbf{Assumption A}: 

Let the DAE \eqref{eq:DAE} be regular with index $\mu\geq 2$ and 
dynamical degree of freedom $l=r_{[\mu-1]}>0$.
Let $x_{\ast}\in H_{D}^{1}(a,b)$ be the solution of (\ref{eq:DAE}), (\ref{eq:Gar})
for given $q\in L^{2}((a,b),\R^{m})$ and $g\in\R^{l}$.  Let all given data as well as the solution $x_*$  be sufficiently smooth.\footnote{The detailed smoothness assumptions result from the requirements of the numerical methods and analysis below.}
Set 
\begin{align*}
 C_{H}= \sqrt{2\max\{H,\frac{1}{H}\}}\geq \sqrt{2}.
\end{align*}
Agree upon that
$G_{\tau}(\bar{t})$ for $\tau>0$ is computed by Algorithm~\ref{alg:Reduction-procedure-1}
while $G_{0}(\bar{t})$ is computed by the same algorithm but with
exact differentiation. Moreover, let there exists an $\tau_{0}>0$ and
a $C_{d}>0$, both independent of $\bar{t}\in[a,b]$, such that (cf.
Conjecture AIC)
\begin{align*}
\lvert G_{\tau}(\bar{t})-G_{0}(\bar{t})\rvert&\leq C_{d}\tau^{N+2-\mu},\\
\rank G_{\tau}(\bar t)&=l, \quad 0\leq\tau\leq \tau_0,\quad  \bar t\in [a,b].
\end{align*}
\qed
\end{description}

\subsection{Least-squares collocation on a subinterval $[\bar{t},\bar{t}+H]$}

In the least-squares collocation procedure to be introduce below,
stepsizes $h$, polynomial degrees $N$, and collocation points $\theta_{i}$,
$i=1,\ldots,M$ are used, too. They are independent of the corresponding
quantities used in the discrete differentiation procedure in Algorithm~\ref{alg:Reduction-procedure-1}.
In order to distinguish them from each other, we will add subscripts
$c$ (collocation) and $d$ (differentiation), respectively: for example,
$M_{c}$ denotes the number of collocation points used in the collocation
approach while $M_{d}$ is the number of collocation points in the
discrete differentiation procedure. Moreover, it is appropriate to
distinguish between the sets of collocation points in both algorithms.
For the discrete differentiation, Chebyshev nodes of the second kind
are most appropriate (Theorem~\ref{th.cheb}) while it is wise to
use Gauss-Legendre or Gauss-Radau nodes for the least-squares collocation
\cite[Theorem 2]{HaMaT122}. For the least-squares collocation, also
Chebyshev nodes of the first kind are reasonable choices.

The interval $\mathcal{I}=[\bar{t},\bar{t}+H]$ will be subdivided
by the grid
\[
\pi:\bar{t}=t_{0}<\cdots<t_{n}=\bar{t}+H
\]
where $t_{j}=\bar{t}+jh$ and $h=H/n$.\footnote{The restriction to constant stepsizes is done for simplifying the
notation. Results for quasi-uniform grids are provided in \cite{HaMaTi19}.} Let $C_{\pi}(\mathcal{I},\R^{m})$ be the space of all function $x:\mathcal{I}\rightarrow\R^{m}$
such that $x$ restricted onto each subinterval $(t_{j-1},t_{j})$
is continuous with continuous extensions onto the boundary points $t_{j-1}$
and $t_{j}$.

Let an $N_{c}>0$ be given. We will approximate the solutions of (\ref{eq:DAE})
on $\mathcal{I}$ subject to the approximate initial condition
\begin{equation}
G_{\tau}(\bar{t})x(\bar{t})=g\label{eq:AccIC}
\end{equation}
by piecewise polynomials of degree $N_{c}$. In general, one has to expect $\ker G_{\tau}(\bar{t})\neq\Ncan(\bar{t})$, which means that 
this initial condition fails to be an accurate one.  This is in contrast to the study in \cite{HaMaPad121}, where it was based exclusively on accurate initial conditions.

We define the
ansatz space for the collocation approach by
\begin{align*}
\mathcal X_{\pi}=\{p\in C_{\pi}(\mathcal{I},\R^{m}): & Dp\in C(\mathcal{I},\R^{k})\\
 & p_{\kappa}\vert_{(t_{j-1},t_{j})}\in\mathfrak{P}_{N},\kappa=1,\ldots,k,\\
 & p_{\kappa}\vert_{(t_{j-1},t_{j})}\in\mathfrak{P}_{N-1},\kappa=k+1,\ldots,m,\\
 & j=1,\ldots,n\}.
\end{align*}
Obviously, $\mathcal X_{\pi}$ is a finite dimensional subspace of $H_{D}^{1}(\mathcal{I})$.
Then, the least-squares functional is given by
\begin{equation}
\Phi_{\tau,\pi}(x)=\int_{\mathcal{I}}\lvert A(t)(Dx(t))'+B(t)x(t)-q(t)\rvert^{2}\dt t+\lvert G_{\tau}(\bar{t})x(\bar{t})-g\rvert^{2},\quad x\in\mathcal X_{\pi}.\label{eq:Phi}
\end{equation}
This functional cannot be evaluated directly. Instead, it must be
approximated by a discrete process. A straightforward discretization
may not be appropriate since the problem (\ref{eq:DAE}), (\ref{eq:AccIC})
is ill-posed \cite{HaMaTiWeWu17} and, thus, extremely sensitive to perturbations.
However, the following approximation is applicable: Let $M_{c}$ points
\[
0<\theta_{1}<\theta_{2}<\cdots<\theta_{M_{c}}<1
\]
 be given. The set of collocation points is
\[
S_{\pi,M_{c}}=\{t_{ji}=t_{j-1}+\theta_{i}h:j=1,\ldots,n,\quad i=1,\ldots,M_{c}\}.
\]
Then we define the interpolation operator $R_{\pi,M_{c}}:C_{\pi}(\mathcal{I},\R^{m})\rightarrow C_{\pi}(\mathcal{I},\R^{m})$
by assigning, to each $u\in C_{\pi}(\mathcal{I},\R^{m})$, the piecewise
polynomial $R_{\pi,M_{c}}u$ with
\[
R_{\pi,M_{c}}u\vert_{(t_{j-1}t_{j})}\in\mathfrak{P}_{M_{c}-1},j=1,\ldots,n,\quad R_{\pi,M_{c}}u(t)=u(t),t\in S_{\pi,M_{c}}.
\]
The discrete functional is given by
\begin{equation}
\Phi_{\tau,\pi,M_{c}}(x)=\int_{\mathcal{I}}\lvert R_{\pi,M_{c}}(A(t)(Dx(t))'+B(t)x(t)-q(t))\rvert^{2}\dt t+\lvert G_{\tau}(\bar{t})x(\bar{t})-g\rvert^{2},\quad x\in\mathcal X_{\pi}.\label{eq:PhiM}
\end{equation}
This functional has a representation \cite{HaMaT122}
\[
\Phi_{\tau,\pi,M_{c}}(x)=U^{\T}\mathcal{L}U+\lvert G_{\tau}(\bar{t})x(\bar{t})-g\rvert^{2},\quad x\in\mathcal X_{\pi},
\]
with a positive definite symmetric matrix $\mathcal{L}$ depending
only on $\theta_{1},\ldots,\theta_{M_{c}}$ and 
\[
U=\begin{bmatrix}
U_{1}\\
\vdots\\
U_{n}
\end{bmatrix}\in\R^{mM_{c}n},\quad U_{j}=\left(\frac{h}{M_{c}}\right)^{1/2}\begin{bmatrix}
A(t_{j1})(Dx(t_{j1}))'+B(t_{j1})x(t_{j1})-q(t_{j1})\\
\vdots\\
A(t_{jM_c})(Dx(t_{jM_c}))'+B(t_{jM_c})x(t_{jM_c})-q(t_{jM_c})
\end{bmatrix}.
\]
The convergence of this method has been proven in \cite{HaMaTi19,HaMa21}
if $\ker G_{\tau}(\bar{t})=\Ncan(\bar{t})$, that is, (\ref{eq:AccIC})
is accurate. As mentioned before, the latter condition is usually not fulfilled. However,
it holds true for $G_{0}(\bar{t})$. We will generalize the corresponding
theorems of \cite{HaMaTi19,HaMa21} for the functionals (\ref{eq:Phi})
and (\ref{eq:PhiM}) now. In order to do so it is convenient to reformulate
the problem in an operator setting. Supposing the General Assumption to be satisfied and $\tau$ to be sufficiently small, we define
\[
\mathcal{T}:\mathcal X\rightarrow\mathcal Y,\quad\mathcal{T}_{\tau}:\mathcal X\rightarrow\mathcal Y,
\quad\mathcal X=H_{D}^{1}(\mathcal{I}),\quad\mathcal Y=L^{2}(\mathcal{I},\R^{m})\times\R^{l}
\]
by
\[
\mathcal{T}x(t)=\begin{bmatrix}
A(t)(Dx(t))'+B(t)x(t)\\
G_{0}(\bar{t})x(\bar{t})
\end{bmatrix},\quad\mathcal{T}_{\tau}x(t)=\begin{bmatrix}
A(t)(Dx(t))'+B(t)x(t)\\
G_{\tau}(\bar{t})x(\bar{t})
\end{bmatrix},\quad x\in\mathcal X.
\]
Note that, by construction, $G_{\tau}(t) = C(t)E(t) = C(t)A(t)D$ for some matrix $C(t)$. Hence, $\ker G_{\tau}(t)\supseteqq\ker D$ for $\tau\geq0$ and all $t\in[a,b]$ such that both $\mathcal{T}$ and  $\mathcal{T}_{\tau}$ are well-defined for $x\in\mathcal X$.
Both operators are bounded and injective. 
Regarding Theorem \ref{t.tildeG} we find that $\mathcal T$ and $\mathcal T_{\tau}$ share their range. Since $\im \mathcal T=\im \mathcal T_{\tau}\subset\mathcal Y$ is a nonclosed, both inverse operators $\mathcal T^{-1}$ and $\mathcal T_{\tau}^{-1}$ are unbounded, and hence we are confronted with ill-posed operator equations.

Let $\mathcal{P}_{\pi}:\mathcal X\rightarrow\mathcal X_{\pi}$ be the orthogonal
projector. Finally, define
\[
\mathcal{R}_{\pi,M_{c}}:C_{\pi}(\mathcal{I},\R^{m})\times\R^{l}\rightarrow\mathcal Y,\quad\mathcal{R}_{\pi,M_{c}}\begin{bmatrix}
w\\
r
\end{bmatrix}=\begin{bmatrix}
R_{\pi,M_{c}} & 0\\
0 & I_{l}
\end{bmatrix}\begin{bmatrix}
w\\
r
\end{bmatrix}.
\]
The initial value problem (\ref{eq:DAE}), (\ref{eq:G0bar}) reads
$\mathcal{T}x=\begin{bmatrix}
\bar{q}\\
g
\end{bmatrix}$ while $\mathcal{T}_{\tau}x=\begin{bmatrix}
\bar{q}\\
g
\end{bmatrix}$ refers to the approximate condition (\ref{eq:AccIC}). Here, $\bar{q}$
is the restriction of $q$ onto $\mathcal{I}$. 
The
minimum norm least-squares solutions of (\ref{eq:Phi}) and (\ref{eq:PhiM})
exist and are given by
\[
(\mathcal{T}_{\tau}\mathcal{P}_{\pi})^{+}\begin{bmatrix}
\bar{q}\\
g
\end{bmatrix}\text{ and }(\mathcal{R}_{\pi,M_{c}}\mathcal{T}_{\tau}\mathcal{P}_{\pi})^{+}\begin{bmatrix}
\bar{q}\\
g
\end{bmatrix},
\]
respectively. Here, $(\mathcal{T}_{\tau}\mathcal{P}_{\pi})^{+}$ and
$(\mathcal{R}_{\pi,M_{c}}\mathcal{T}_{\tau}\mathcal{P}_{\pi})^{+}$
denote the Moore-Penrose inverses of the corresponding operators.

In order to motivate the following derivations, an error representation
will be derived where we follow roughly \cite[Eq. (2.9)]{KaOf12}.
For higher index DAEs, any initial value problem is ill-posed. In
particular, this means that $\mathcal{T}$ is not surjective. Fix
an $[\bar{q},g]^{\T}\in\im\mathcal{T}$. Since (\ref{eq:G0bar})
is an accurate initial condition, the solution $x_{g}$ of $\mathcal{T}x=[\bar{q},g]^{\T}$
is unique. Later on, we will show that $\mathcal{T}_{\tau}\mathcal{P}_{\pi}$
restricted onto $\mathcal X_{\pi}$ is injective. Then it holds $(\mathcal{T}_{\tau}\mathcal{P}_{\pi})^{+}\mathcal{T}_{\tau}\mathcal{P}_{\pi}=I_{\mathcal X_{\pi}}$.
Finally, assume that $[q^{\delta},g^{\delta}]^{\T}\in\mathcal Y$ is
a perturbed right-hand side with $\lVert[q^{\delta}-\bar{q},g^{\delta}-g]^{\T}\rVert_{\mathcal Y}\leq\delta$
which does not necessarily belong to $\im\mathcal{T}$. The discrete
least-squares solution for the perturbed right-hand side becomes $x_{\tau,\pi}^{\delta}=(\mathcal{T}_{\tau}\mathcal{P}_{\pi})^{+}[q^{\delta},g^{\delta}]^{\T}$.
Then it holds
\begin{align*}
x_{\tau,\pi}^{\delta}-x_{g} & =x_{\tau,\pi}^{\delta}-\mathcal{P}_{\pi}x_{g}-(I_{X}-\mathcal{P}_{\pi})x_{g}\\
 & =(\mathcal{T}_{\tau}\mathcal{P}_{\pi})^{+}\mathcal{T}_{\tau}\mathcal{P}_{\pi}(x_{\tau,\pi}^{\delta}-\mathcal{P}_{\pi}x_{g})-(I_{X}-\mathcal{P}_{\pi})x_{g}\\
 & =(\mathcal{T}_{\tau}\mathcal{P}_{\pi})^{+}([q^{\delta},g^{\delta}]^{\T}-\mathcal{T}_{\tau}\mathcal{P}_{\pi}x_{g})-(I_{X}-\mathcal{P}_{\pi})x_{g}\\
 & =(\mathcal{T}_{\tau}\mathcal{P}_{\pi})^{+}[q^{\delta}-\bar{q},g^{\delta}-g]^{\T}+(\mathcal{T}_{\tau}\mathcal{P}_{\pi})^{+}([\bar{q},g]^{\T}-\mathcal{T}_{\tau}\mathcal{P}_{\pi}x_{g})-(I_{X}-\mathcal{P}_{\pi})x_{g}\\
 & =(\mathcal{T}_{\tau}\mathcal{P}_{\pi})^{+}[q^{\delta}-\bar{q},g^{\delta}-g]^{\T}+(\mathcal{T}_{\tau}\mathcal{P}_{\pi})^{+}(\mathcal{T}x_{g}-\mathcal{T}_{\tau}\mathcal{P}_{\pi}x_{g})-(I_{X}-\mathcal{P}_{\pi})x_{g}\\
 & =(\mathcal{T}_{\tau}\mathcal{P}_{\pi})^{+}[q^{\delta}-\bar{q},g^{\delta}-g]^{\T}\\
 & \phantom{(\mathcal{T}_{h}\mathcal{P}_{h})^{+}}+(\mathcal{T}_{\tau}\mathcal{P}_{\pi})^{+}\mathcal{T}(I_{X}-\mathcal{P}_{\pi})x_{g}+(\mathcal{T}_{\tau}\mathcal{P}_{\pi})^{+}(\mathcal{T}-\mathcal{T}_{\tau})\mathcal{P}_{\pi}x_{g}-(I_{X}-\mathcal{P}_{\pi})x_{g}.
\end{align*}
Define 
\begin{align*}
\alpha_{\pi} & =\lVert(I_{X}-\mathcal{P}_{\pi})x_{g}\rVert_{\mathcal X}\\
\beta_{\pi} & =\lVert\mathcal{T}(I_{X}-\mathcal{P}_{\pi})x_{g}\rVert_{\mathcal Y}\\
\gamma_{\tau,\pi} & =\inf_{x\in\mathcal X_{\pi},x\neq0}\frac{\lVert\mathcal{T}_{\tau}x\rVert_{\mathcal Y}}{\lVert x\rVert_{\mathcal X}}=\lVert(\mathcal{T}_{\tau}\mathcal{P}_{\pi})^{+}\rVert^{-1}\\
\eta_{\tau,\pi} & =\lVert(\mathcal{T}-\mathcal{T}_{\tau})\mathcal{P}_{\pi}x_{g}\rVert_{\mathcal Y}.
\end{align*}
If $\gamma_{\tau,\pi}>0$, it holds
\begin{equation}
\lVert x_{\tau,\pi}^{\delta}-x_{g}\rVert_{\mathcal X}\leq\frac{\delta}{\gamma_{\tau,\pi}}+\frac{\beta_{\pi}}{\gamma_{\tau,\pi}}+\frac{\eta_{\tau,\pi}}{\gamma_{\tau,\pi}}+\alpha_{\pi}.\label{eq:TP}
\end{equation}
Similar derivations can be done for $\hat{x}_{\tau,\pi}^{\delta}=(\mathcal{R}_{\pi,M_{c}}\mathcal{T}_{\tau}\mathcal{P}_{\pi})^{+}[q^{\delta},g^{\delta}]^{\T}$
such that
\begin{equation}
\lVert\hat{x}_{\tau,\pi}^{\delta}-x_{g}\rVert_{\mathcal X}\leq\frac{\delta}{\hat{\gamma}_{\tau,\pi}}+\frac{\beta_{\pi}}{\hat{\gamma}_{\tau,\pi}}+\frac{\hat{\eta}_{\tau,\pi}}{\hat{\gamma}_{\tau,\pi}}+\alpha_{\pi}\label{eq:RTP}
\end{equation}
where
\begin{align*}
\hat{\gamma}_{\tau,\pi} & =\inf_{x\in\mathcal X_{\pi},x\neq0}\frac{\lVert\mathcal{R}_{\pi,M_{c}}\mathcal{T}_{\tau}x\rVert_{\mathcal Y}}{\lVert x\rVert_{\mathcal X}}=\lVert(\mathcal{R}_{\pi,M_{c}}\mathcal{T}_{\tau}\mathcal{P}_{\pi})^{+}\rVert^{-1}\\
\hat{\eta}_{\tau,\pi} & =\lVert(\mathcal{T}-\mathcal{R}_{\pi,M_{c}}\mathcal{T}_{\tau})\mathcal{P}_{\pi}x_{g}\rVert_{\mathcal Y}.
\end{align*}

\begin{proposition}
\label{prop-alpha}If the General Assumption A is fulfilled and $x_g:=x_*$ then:
\begin{enumerate}
\item There exist constants $C_{\alpha}$ and $C_{\beta}$ such that, for
all sufficiently fine grids,
\[
\alpha_{\pi}\leq C_{\alpha}H^{1/2}h^{N_{c}}\lVert x_{*}^{(N_{c}+1)}\rVert_{C(\mathcal{I},\R^{m})},\quad\beta_{\pi}\leq C_{\beta}h^{N_{c}}\lVert x_{*}^{(N_{c}+1)}\rVert_{C(\mathcal{I},\R^{m})}.
\]
The constants $C_{\alpha},C_{\beta}$ are independent of $\mathcal{I}$.
\item There exists a constant $C_{\eta}$ such that, for all sufficiently
fine grids,
\[
\eta_{\tau,\pi}\leq C_{\eta}C_{H}\tau^{N_{d}+2-\mu}\lVert x_{*}\rVert_{H_{D}^{1}(\mathcal{I})}.
\]
The constant $C_{\eta}$ is independent of $\mathcal{I}$.
\item Let $M_{c}\geq N_{c}+1$. Then there exists a constant $\hat{C}_{\eta}$
such that, for all sufficiently fine grids,
\[
\hat{\eta}_{\tau,\pi}\leq\hat{C}_{\eta}(h^{2M_{c}-2N_{c}-1}+C_{H}^{2}\tau^{2(N_{d}+2-\mu)})^{1/2}\lVert x_{*}\rVert_{H_{D}^{1}(\mathcal{I})}.
\]
The constant $\hat{C}_{\eta}$ is independent of $\mathcal{I}$.
\end{enumerate}
\end{proposition}
\begin{svmultproof}
A careful analysis of the proof of \cite[Theorem 3.4]{HaMaPad121}
shows that assertion 1 holds. The assumptions of that theorem are
fulfilled by \cite[Theorem 4.3]{HaMa231}, \cite[Theorem 2.42]{LMT}, and Lemma~\ref{Cor:Gbound}.

In order to prove assertion 2 consider, for $x\in H_{D}^{1}(\mathcal{I})$,
\[
(\mathcal{T}-\mathcal{T}_{\tau})x=\begin{bmatrix}
0\\
(G_{0}(\bar{t})-G_{\tau}(\bar{t}))x(\bar{t})
\end{bmatrix}.
\]
Hence,
\[
\lVert(\mathcal{T}-\mathcal{T}_{\tau})x\rVert_{\mathcal Y}=\lvert(G_{0}(\bar{t})-G_{\tau}(\bar{t}))x(\bar{t})\rvert\leq C_{d}\tau^{N_{d}+2-\mu}\lvert Dx(\bar{t})\rvert\leq C_{d}C_{H}\tau^{N_{d}+2-\mu}\lVert x\rVert_{\mathcal X}
\]
by \cite[Lemma 3.2]{HaMaPad121}. It follows $\eta_{\tau,\pi}\leq C_{d}C_{H}\tau^{N_{d}+2-\mu}\lVert x_{*}\rVert_{H_{D}^{1}(\mathcal{I})}$,
that is, assertion 2. 

Finally, we estimate
\begin{multline*}
\lVert(\mathcal{T}-\mathcal{R}_{\pi,M_{c}}\mathcal{T}_{\tau})\mathcal{P}_{\pi}x_{*}\rVert_{\mathcal Y}^{2} \\ =\lVert A(D\mathcal{P}_{\pi}x_{*})'+B\mathcal{P}_{\pi}x_{*}-R_{\pi,M_{c}}(A(D\mathcal{P}_{\pi}x_{*})'+B\mathcal{P}_{\pi}x_{*})\rVert_{L^{2}(\mathcal{I},\R^{m})}^{2}\\
+\lvert(G_{0}(\bar{t})-G_{\tau}(\bar{t}))\mathcal{P}_{\pi}x_{*}(\bar{t})\rvert^{2}.
\end{multline*}
Invoking \cite[Lemma 3.3(2)]{HaMaPad121}, the first term can be estimated
by 
\begin{multline*}
\lVert A(D\mathcal{P}_{\pi}x_{*})'+B\mathcal{P}_{\pi}x_{*}-R_{\pi,M_{c}}(A(D\mathcal{P}_{\pi}x_{*})'+B\mathcal{P}_{\pi}x_{*})\rVert_{L^{2}(\mathcal{I},\R^{m})}\\
\leq C_{AB1}h^{M_{c}-N_{c}-1/2}\lVert x_{*}\rVert_{H_{D}^{1}(\mathcal{I})}
\end{multline*}
where $C_{AB1}$ depends only on the coefficients $A,B$ but not on
$\mathcal{I}$. The second term can be estimated as before. Summarizing,
we obtain
\[
\hat{\eta}_{\tau,\pi}\leq\hat{C}_{\eta}(h^{2M_{c}-2N_{c}-1}+C_{H}^{2}\tau^{2(N_{d}+2-\mu)})^{1/2}\lVert x_{*}\rVert_{H_{D}^{1}(\mathcal{I})}.
\]
\end{svmultproof}

\begin{proposition}
\label{Cor:gamma}Let General Assumption A be fulfilled. Let $N_{c}\geq\mu$.
Moreover, let $\tau$ and $N_{d}$ be chosen such that $\tau^{N_{d}+2-\mu}\leq h^{N_{c}}$.
Then it holds:
\begin{enumerate}
\item There exists a constant $C_{\gamma}$ and an
$h_{0}=h_{0}(H)$ such that, for all $h\leq h_{0}$, it holds
\[
\gamma_{\tau,\pi}\geq C_{\gamma}h^{\mu-1}
\]
with a constant $C_{\gamma}$ independent of $\mathcal{I}$.
\item Let, additionally, $M_{c}\geq N_{c}+\mu$. Then
there exists a constant $\hat{C}_{\gamma}$ and an $h_{0}=h_{0}(H)$
such that, for all $h\leq h_{0}$, it holds
\[
\hat{\gamma}_{\tau,\pi}\geq\hat{C}_{\gamma}h^{\mu-1}
\]
with a constant $\hat{C}_{\gamma}$independent of $\mathcal{I}$.
\end{enumerate}
\end{proposition}
\begin{svmultproof}
We apply \cite[Theorem 3.5(2)]{HaMaPad121}. The assumptions of that
theorem are fulfilled by \cite[Theorem 4.3]{HaMa231}, \cite[Theorem 2.42]{LMW},
and Corollary~\ref{Cor:Gbound}.\footnote{The assumption $M_{c}\geq N_{c}+1$ in \cite{HaMaPad121} is not necessary
in the present case since neither $\mathcal{T\mathcal{P}}_{\pi}$
nor $\mathcal{T}_{\tau}\mathcal{P}_{\pi}$ depend on $M_{c}$.} This provides us with 
\[
\lVert(\mathcal{T\mathcal{P}}_{\pi})^{+}\rVert^{-1}=:\tilde{\gamma}_{\pi}\geq\tilde{C}_{\gamma}h^{\mu-1}
\]
for sufficiently small $h$ with $\tilde{C}_{\gamma}$ being independent
of $\mathcal{I}$. Equivalently, $\tilde{\gamma}_{\pi}\lVert x\rVert_{\mathcal X}\leq\rVert\mathcal{T}x\rVert_{\mathcal Y}$
for all $x\in\mathcal X_{\pi}$. Then we have similar to the proof of Proposition~\ref{prop-alpha},
for $x\in\mathcal X_{\pi}$,
\begin{align*}
\tilde{\gamma}_{\pi}\lVert x\rVert_{\mathcal X}\leq\rVert\mathcal{T}x\rVert_{\mathcal Y} 
& \leq\lVert\mathcal{T}_{\tau}x\rVert_{\mathcal Y}+\lVert(\mathcal{T}-\mathcal{T}_{\tau})x\rVert_{\mathcal Y}.\\
 & \leq\lVert\mathcal{T}_{\tau}x\rVert_{\mathcal Y}+C_{\eta}C_{H}\tau^{N_{d}+2-\mu}\lVert x\rVert_{\mathcal X}\\
 & \leq\lVert\mathcal{T}_{\tau}x\rVert_{\mathcal Y}+C_{\eta}C_{H}h^{N_{c}}\lVert x\rVert_{\mathcal X}
\end{align*}
Hence,
\[
(\tilde{C}_{\gamma}-C_{\eta}C_{H}h^{N_{c}-\mu+1})h^{\mu-1}\lVert x\rVert_{\mathcal X}\leq\lVert\mathcal{T}_{\tau}x\rVert_{\mathcal Y}.
\]
If $h_{0}$ is chosen such that $\frac{C_{\eta}C_{H}}{\tilde{C}_{\gamma}}h_0^{N_{c}-\mu+1}\leq\frac{1}{2}$,
$C_{\gamma}=\tilde{C}_{\gamma}/2$ has the required property. In order
to prove the second statement it suffices to invoke \cite[Theorem 3.5(3)]{HaMaPad121}.
The derivations before can be repeated verbatim.
\end{svmultproof}

The next corollary is an immediate consequence.
\begin{corollary}
Under the assumptions of Proposition~\ref{Cor:gamma}, both $\mathcal{T}_{\tau}\mathcal{P}_{\pi}$
and $\mathcal{R}_{\pi,M_{c}}\mathcal{T}_{\tau}\mathcal{P}_{\pi}$
restricted onto $\mathcal X_{\pi}$ are injective for $h\leq h_{0}$.
\end{corollary}

\begin{theorem}
\label{theo:Conv}Let General Assumption A be fulfilled. Let $N_{c}\geq\mu$.
Moreover, let $\tau$ and $N_{d}$ be chosen such that $\tau^{N_{d}+2-\mu}\leq h^{N_{c}}$. 
\begin{enumerate}
\item Then there exists an $h_{0}=h_{0}(H)$ such that,
for all $h\leq h_{0}$, there is a unique solution $x_{\tau,\pi}$
of the least-squares problem
\[
\min\{\Phi_{\tau,\pi}(x)\vert x\in\mathcal X_{\pi}\}
\]
and there exists a constant $C_{*}$ independent of $\mathcal{I}$
such that
\[
\lVert x_{\tau,\pi}-x_{\ast}\rVert_{H_{D}^{1}(\mathcal{I})}\leq C_{\ast}C_{H}h^{N_{c}-\mu+1}.
\]
\item Let, additionally, $M_{c}\geq N_{c}+\mu$. Then there exists
an $h_{0}=h_{0}(H)$ such that, for all $h\leq h_{0}$, there is a
unique solution $x_{\tau,\pi,M_{c}}$ of the least-squares problem
\[
\min\{\Phi_{\tau,\pi,M_{c}}(x)\vert x\in\mathcal X_{\pi}\}
\]
and there is a constant $\tilde{C}_{*}$ such that
\[
\lVert x_{\tau,\pi,M_{c}}-x_{\ast}\rVert_{H_{D}^{1}(\mathcal{I})}\leq\tilde{C}_{\ast}C_{H}h^{N_{c}-\mu+1}.
\]
\end{enumerate}
\end{theorem}
\begin{svmultproof}
Apply the estimates (\ref{eq:TP}) and (\ref{eq:RTP}), respectively,
with $\delta=0$ and $x_g=x_*$.
\end{svmultproof}
\bigskip

The results of Theorems \ref{th-Scan} and \ref{theo:Conv} suggest
the following numerical procedure for solving initial value problems on the entire interval $[a,b]$ : 
\begin{enumerate}
\item Let the interval $[a,b]$ be decomposed into $L$ subintervals,
\[
a=w_{0}<w_{1}<\cdots<w_{L}=b,
\]
with the length $H_{\lambda}=w_{\lambda}-w_{\lambda-1}$, $\lambda=1,\ldots,L$.
Denote $\mathcal{I}_{\lambda}=[w_{\lambda-1},w_{\lambda}]$.
\item For $\lambda=1$, compute the approximation $x_{\tau,\pi}^{[1]}\in H_{D}^{1}(\mathcal{I}_{1})$
(or, $x_{\tau,\pi,M_{c}}^{[1]}\in H_{D}^{1}(\mathcal{I}_{1})$, respectively)
to the solution of the initial value problem (\ref{eq:DAE})-(\ref{eq:Gar})
on $\mathcal{I}_{1}$ by using the functional $\Phi_{\tau,\pi}$ (or,
$\Phi_{\tau,\pi,M_c}$, respectively). with $G_{\tau}(a)$ replaced
by $G^{a}$.
\item For $\lambda>1$, determine $G_{\tau}(w_{\lambda-1})$ by using Algorithm
(\ref{alg:Reduction-procedure-1}) at $w_{\lambda-1}$. Let $x_{\tau,\pi}^{[\lambda-1]}\in H_{D}^{1}(\mathcal{I}_{\lambda-1})$
(or, $x_{\tau,\pi,M_{c}}^{[\lambda-1]}\in H_{D}^{1}(\mathcal{I}_{\lambda-1})$,
respectively) be available. Then solve the initial value problem (\ref{eq:DAE})
on $\mathcal{I}_{\lambda}$ subject to
\begin{equation}
G_{\tau}(w_{\lambda-1})x(w_{\lambda-1})=g_{\lambda-1},\quad g_{\lambda-1}=G_{\tau}(w_{\lambda-1})x_{\tau,\pi}^{[\lambda-1]}(w_{\lambda-1})\label{eq:subIC}
\end{equation}
or
\begin{equation}
G_{\tau}(w_{\lambda-1})x(w_{\lambda-1})=g_{\lambda-1},\quad g_{\lambda-1}=G_{\tau}(w_{\lambda-1})x_{\tau,\pi,M_c}^{[\lambda-1]}(w_{\lambda-1})\label{eq:subICR}
\end{equation}
for $x_{\pi}^{[\lambda]}\in H_{D}^{1}(\mathcal{I}_{\lambda})$ (or,
$x_{\pi,M_c}^{[\lambda]}\in H_{D}^{1}(\mathcal{I}_{\lambda})$, respectively).
\end{enumerate}

In order to simplify the notations, we omit the index $\tau$ in
the following.

The following theorem will be proven for constant stepsizes $H_{\lambda}=H$. It is
easy to see that it can be generalized to the case of quasi-uniform
grids $\rho H\leq H_{\lambda}\leq H$ with a constant $\rho$ independent
of $H$.
\begin{theorem}
\label{Theo:StepConv}
Let General Assumption A be fulfilled. Let $N_c\geq\mu$ and 
let $\tau$ and $N_{d}$ be chosen such that $\tau^{N_{d}+2-\mu}\leq h^{N_c}$. 

Moreover, let $C_N$ denote a uniform bound for all DAE solutions corresponding to the given $q$ and belonging to the set $U_r$ with sufficiently large $r$ in the sense of Remark \ref{r.derivbound} such that 
\begin{align*}
 \lvert x^{(N_c+1)}(t)\rvert\leq C_{N},\quad t\in[a,b],\quad x\in U_r=:\mathcal U, 
\end{align*}
Let $H_{\lambda}\equiv H=(b-a)/L$ be constant and $G_{h}(w_{\lambda-1})$,
$\lambda=1,\ldots,L-1$ be determined by Algorithm (\ref{alg:Reduction-procedure-1}). 
\begin{enumerate}
\item Then, $x_{\pi}^{[1]},\ldots,x_{\pi}^{[L]}$ are uniquely
determined. There exists a constant $\tilde C>0$ such that the estimate
\[
\lVert x_{\pi}^{[\lambda]}-x_{\ast}\rVert_{H_{D}^{1}(w_{\lambda-1},w_{\lambda})}\leq C_{\lambda}h^{N_c-\mu+1},\quad C_{\lambda}=\tilde{C}C_{H}\frac{(\tilde{C}C_{H})^{\lambda}-1}{\tilde{C}C_{H}-1}
\]
holds true for sufficiently small $h>0$.\footnote{A careful analysis of the proof shows that a suitable constant $\tilde C>0$ is given by $\tilde C = \max(\bar C,C_pc_G,C_\ast)$ with
$\bar C = \max(C_\eta\lVert x_\ast\rVert_{H^1_D(a,b)}/C_\gamma,C_\beta C_N/C_\gamma,C_\alpha C_N)$. Similar expressions can be provided for $\hat{\tilde C}$ in Assertion 2.}
\item Let $M_{c}>2N_c$. Then, $x_{\pi,M}^{[1]},\ldots,x_{\pi,M}^{[L]}$ are uniquely
determined. There exists a constant $\hat{\tilde{C}}>0$ such that the estimate
\[
\lVert x_{\pi,M}^{[\lambda]}-x_{\ast}\rVert_{H_{D}^{1}(w_{\lambda-1},w_{\lambda})}\leq\hat{C}_{\lambda}h^{N_c-\mu+1},\quad\hat{C}_{\lambda}=\hat{\tilde{C}}C_{H}\frac{(\hat{\tilde{C}}C_{H})^{\lambda}-1}{\hat{\tilde{C}}C_{H}-1}
\]
holds true for sufficiently small $h>0$.
\end{enumerate}
\end{theorem}
\begin{svmultproof}
For $\lambda=1$, Theorem~(\ref{theo:Conv}) provides the estimate
\[
d_{1}:=\lVert x_{\pi}^{[1]}-x_{\ast}\rVert_{H_{D}^{1}(\mathcal{I}_{1})}\leq C_{\ast}C_{H}h^{N_c-\mu+1}.
\]
For $\lambda>1$, let $\gamma=G_{h}(w_{\lambda-1})x_{\pi}^{[\lambda-1]}(w_{\lambda-1})$
and $x_{\gamma}$ be the solution of (\ref{eq:DAE}), (\ref{eq:G0bar})
on $\mathcal{I}_{\lambda}$. Then we estimate using \cite[Theorem 3.5]{HaMaPad121}
(note that $C_{H}\geq\sqrt{2}$)
\begin{align*}
d_{\lambda}&:=\lVert x_{\pi}^{[\lambda]}-x_{\ast}\rVert_{H_{D}^{1}(\mathcal{I}_{\lambda})} \\
& \leq\lVert x_{\pi}^{[\lambda]}-x_{\gamma}\rVert_{H_{D}^{1}(\mathcal{I}_{\lambda})}+\lVert x_{\gamma}-x_{\ast}\rVert_{H_{D}^{1}(\mathcal{I}_{\lambda})}\\
 & \leq\left(\frac{\beta_{\pi}}{\gamma_{\pi}}+\frac{\eta_{\pi}}{\gamma_{\pi}}+\alpha_{\pi}\right)+C_{p}\lvert\gamma-G_{h}(w_{\lambda-1})x_{\ast}(w_{\lambda-1})\rvert\\
 & \leq\bar{C}\left[C_{H}h^{N_c-\mu+1}+H^{1/2}h^{N_c}\right]+C_{p}\lvert G_{h}(w_{\lambda-1})(x_{\pi}^{[\lambda-1]}(w_{\lambda-1})-x_{\ast}(w_{\lambda-1}))\rvert\\
 & \leq\bar{C}\left[C_{H}h^{N_c-\mu+1}+H^{1/2}h^{N_c}\right]+C_{p}c_{G}C_{H}\lVert x_{\pi}^{[\lambda-1]}-x_{\ast}\rVert_{H_{D}^{1}(\mathcal{I}_{\lambda-1})}\\
 & \leq\bar{C}\left[C_{H}h^{N_c-\mu+1}+H^{1/2}h^{N_c}\right]+C_{p}c_{G}C_{H}d_{\lambda-1}\\
 & \leq\tilde{C}C_{H}(h^{N_c-\mu+1}+d_{\lambda-1})
\end{align*}
for sufficiently small $h$. 
From this recursion we obtain
\[
d_{\lambda}\leq\tilde{C}C_{H}\frac{(\tilde{C}C_{H})^{\lambda}-1}{\tilde{C}C_{H}-1}h^{N_c-\mu+1}.
\]
This proves assertion 1. Assertion 2 can be proved similarly.
\end{svmultproof}

\begin{remark}
Note that, for sufficiently small $H$, $\tilde{C}C_{H}>1$ and $\hat{\tilde{C}}C_{H}>1$
as well as $C_{H}=\sqrt{2}H^{-1/2}$. Thus, for constant $h$, the
bound $d_{\lambda}$ for the error is increasing with decreasing $H$!\hfill\qed
\end{remark}

We observe that the computational expense for computing $G_{\tau}(w_{\lambda-1})$
increases with increasing $N_{d}$, and the computation becomes more
sensitive with increasing $N_{d}$ and decreasing $\tau$. Therefore,
it is appropriate to choose $N=N_{c}=N_{d}$.

\begin{corollary}
Let the assumptions of Theorem~\ref{Theo:StepConv} be fulfilled. However, the assumptions on $\tau$ and $h$
shall be replaced by 
$N=N_c=N_d\geq\mu$
and $\tau=h^{\mu/2}\leq1$. Then the assertions of Theorem~\ref{Theo:StepConv} hold true.
\end{corollary}
The proof is simply given by observing that
\[
\tau=h^{\mu/2}=h^{1+(\mu-2)/2}\leq h^{1+(\mu-2)/(N-\mu+2)}
\]
such that $\tau^{N+2-\mu}\leq h^{N}$.

In the special case of the availability of the exact derivative $C_{0}'$,
the condition on $\tau$ and $h$ can be replaced by $\tau^{N_d+3-\mu}\leq h^{N_c}$ which follows
from \eqref{eq:simplified-error}. Similarly as above, if $N=N_d=N_c$, $\tau=h^{\mu/3}\leq 1$ is sufficient for this condition to be fulfilled. In particular, 
for $\mu=3$, the choice $\tau=h$ becomes appropriate.

Consider the concatenated functions $x_{L},x_{M,L}\in L^{2}((a,b),\R^{m})$
given by $x_{L}(t)=x_{\pi}^{[\lambda]}(t)$, $t\in\mathcal{I}_{\lambda}$
and $x_{M,L}(t)=x_{\pi,M}^{[\lambda]}(t),t\in\mathcal{I}_{\lambda}$.
In general, it does not hold $x_{L}\in H_{D}^{1}(a,b)$ or $x_{M,L}\in H_{D}^{1}(a,b)$
since the continuity of the differentiated components cannot be guaranteed.
In the following we will use the space
\[
\hat{H}_{D,L}^{1}(a,b)=\{x\in L^{2}((a,b),\R^{m}):x\lvert_{\mathcal{I}_{\lambda}}\in H_{D}^{1}(\mathcal{I}_{\lambda})\}
\]
equipped with the norm given by
\[
\lVert x\rVert_{\hat{H}_{D,L}^{1}(a,b)}^{2}=\sum_{\lambda=1}^{L}\lVert x\lvert_{\mathcal{I}_{\lambda}}\rVert_{H_{D}^{1}(\mathcal{I}_{\lambda})}^{2}.
\]

\begin{corollary}
\label{Cor:StepConv}Let the assumptions of Theorem~\ref{Theo:StepConv}
be fulfilled. Then it holds
\[
\lVert x_{L}-x_{\ast}\rVert_{\hat{H}_{D,L}^{1}(a,b)}\leq\frac{\left(\tilde{C}C_{H}\right)^{2}}{\tilde{C}C_{H}-1}\left(\frac{(\tilde{C}C_{H})^{2L}-1}{(\tilde{C}C_{H})^{2}-1}\right)^{1/2}h^{N_c-\mu+1}
\]
and
\[
\lVert x_{M,L}-x_{\ast}\rVert_{\hat{H}_{D,L}^{1}(a,b)}\leq\frac{\left(\hat{\tilde{C}}C_{H}\right)^{2}}{\hat{\tilde{C}}C_{H}-1}\left(\frac{(\hat{\tilde{C}}C_{H})^{2L}-1}{(\hat{\tilde{C}}C_{H})^{2}-1}\right)^{1/2}h^{N_c-\mu+1}.
\]
In particular, for sufficiently small $H$, it holds 
\begin{align*}
\lVert x_{L}-x_{\ast}\rVert_{\hat{H}_{D,L}^{1}(a,b)} &= O(H^{-(L+1)})h^{N_c-\mu+1},\\
\lVert x_{M,L}-x_{\ast}\rVert_{\hat{H}_{D,L}^{1}(a,b)} &= O(H^{-(L+1)})h^{N_c-\mu+1}.
\end{align*}
\end{corollary}
\begin{svmultproof}
It holds
\begin{align*}
\lVert x_{L}-x_{\ast}\rVert_{\hat{H}_{D,L}^{1}(a,b)}^{2} & =\sum_{\lambda=1}^{L}\lVert x_{_{\pi}}^{[\lambda]}-x_{*}\lvert_{\mathcal{I}_{\lambda}}\rVert_{H_{D}^{1}(\mathcal{I}_{\lambda})}^{2}\\
 & \leq\sum_{\lambda=1}^{L}\left[\tilde{C}C_{H}\frac{(\tilde{C}C_{H})^{\lambda}-1}{\tilde{C}C_{H}-1}h^{N_c-\mu+1}\right]^{2}\\
 & \leq\left(\frac{\tilde{C}C_{H}}{\tilde{C}C_{H}-1}\right)^{2}\sum_{\lambda=1}^{L}(\tilde{C}C_{H})^{2\lambda}h^{2(N_c-\mu+1)}\\
 & =\left(\frac{\tilde{C}C_{H}}{\tilde{C}C_{H}-1}\right)^{2}\left(\frac{(\tilde{C}C_{H})^{2(L+1)}-1}{(\tilde{C}C_{H})^{2}-1}-1\right)h^{2(N_c-\mu+1)}.
\end{align*}
Hence,
\[
\lVert x_{L}-x_{\ast}\rVert_{\hat{H}_{D,L}^{1}(a,b)}\leq\frac{\left(\tilde{C}C_{H}\right)^{2}}{\tilde{C}C_{H}-1}\left(\frac{(\tilde{C}C_{H})^{2L}-1}{(\tilde{C}C_{H})^{2}-1}\right)^{1/2}h^{N_c-\mu+1}.
\]
\end{svmultproof}
These estimates give rise to some comments:
\begin{itemize}
\item If the number $L$ of subintervals in the stepwise procedure (and thus,
$H$) is fixed, the asymptotic behavior with respect to the stepsize
$h$ and the approximation order $N_c$ is retained from the global
collocation approach.
\item For larger $L$ and thus small $H$, it holds $C_{H}=O(H^{-1/2})=O(L^{1/2})$.
Hence, the constant in front of $h^{N_c-\mu+1}$ grows such that smaller
stepsizes $h$ may be necessary in order to obtain an accuracy which
is in par with that of the global approach. However, the stepwise
approach may be faster and requiring less memory since the $L$ subproblems
to be solved are considerably smaller than the global one.
\item The stepwise procedure requires an additional overhead for computing
$G_{h}(w_{\lambda})$.
\end{itemize}
It should be noted that Algorithm~\ref{alg:Reduction-procedure}
(with an inhomogeneous right hand side in (\ref{eq:DAE})) has been
proposed before to be used to construct a time stepping method in
\cite{RaRh94}. The basic idea there is to reduce the DAE to an index-1
DAE (or an ODE) by using the reduction procedure and then to apply
an appropriate finite difference discretization. An implementation
of this method requires that the reduction procedure must be carried
out for each $t\in[a,b]$ where the reduced DAE must be evaluated.
Given the results of Section~\ref{sec:AIC}, such a procedure becomes
implementable. However, it is much more expensive than the algorithm
substantiated above.

\subsection{A numerical example}

We solve the DAE (\ref{eq:DAE}) with the matrices from Section \ref{subsec:LinCaMo}
with $\rho=5$ on $[0,5]$ subject to the accurate initial condition
\[
\begin{bmatrix}
0 & -1 & 0 & 0 & 0 & 0 & 0\\
0 & 1 & 1 & 0 & 0 & 0 & 0\\
0 & 0 & 0 & 0 & -1 & 0 & 0\\
-1 & 0 & 0 & 0 & 1 & 1 & 0
\end{bmatrix}x(0)=\begin{bmatrix}
-1\\
3\\
0\\
0
\end{bmatrix}.
\]
The right hand side $q$ is chosen in such a way that the exact solution
becomes
\begin{align*}
x_{\ast,1}(t)=\sin t, & x_{\ast,4}(t)=\cos t,\\
x_{\ast,2}(t)=\cos t, & x_{\ast,5}(t)=-\sin t,\\
x_{\ast,3}(t)=2\cos^{2}t, & x_{\ast,6}(t)=-2\sin2t,\\
x_{\ast,7}(t)=-\rho^{-1}\sin t.
\end{align*}
This corresponds to the setting in \cite{HaMaTi19,HaMaPad121}.

For each polynomial degree $N=N_c=N_d$, the problem is solved with $M_{d}=M_{c}=N+1$
(in the case of spectral differentiation) and $M_{d}=M_{c}=N+2$ (in
the case of least-squares differentiation) collocation points. These
values are used both in the least-squares collocation method as well
as in the computation of accurate initial conditions. Moreover, the
stepsizes for both algorithms are equal. For the computation of accurate
initial conditions central differences with Chebyshev nodes of the
second kind are used. The collocation nodes for the least-squares collocation
are the Gauss-Legendre points. The results are presented in Tables~\ref{tab:ErrorSpec}
and \ref{tab:ErrorsLSQ} for different values of $L$ and $h=5/(Ln)$.
In Table~\ref{tab:ErrorSpec}, spectral differentiation is used while
Table~\ref{tab:ErrorsLSQ} presents the error for the case of differentiation
by least-squares approximation. Note that, for $N\leq2$, the method
is not guaranteed to converge. However, even in that case we observe
bounded errors. Otherwise, the expected order of convergence $h^{N-2}$
is observed for constant $L$. Moreover, we observe that, for constant
$h$ and $N$, the error decreases with increasing $L$ which is in
accordance with the estimates in Theorem~\ref{Theo:StepConv} and
Corollary~\ref{Cor:StepConv}. However, the influence is only moderate.
Note that Theorem \ref{Theo:StepConv} requires the condition $M_{c}\geq2N$
to hold for convergence. However, our experiments indicate that even
a much smaller number of collocation points does not disturb the order
of convergence.

It should be noted that the computing
time is much too short to be measurable.

\begin{table}
\caption{Errors for the initial value solver in the norm $\hat{H}_{D,L}^{1}(0,5)$
in the case of using spectral differentiation for determining accurate
transfer conditions. In case $n=Ln$, we have $L=1$ such that the
global least-squares method is used\label{tab:ErrorSpec}}

\centering{}%
\begin{tabular}{|c|c|c|c|c|c|c|c|}
\hline 
$L\times n$ & $L$ & $n\backslash N$ & 2 & 4 & 6 & 8 & 10\tabularnewline
\hline 
\hline 
10 & 10 & 1 & 5.06e-01 & 1.18e-02 & 7.60e-05 & 2.67e-07 & 5.39e-10\tabularnewline
\hline 
 & 5 & 2 & 5.19e-01 & 8.92e-03 & 6.27e-05 & 2.05e-07 & 4.05e-10\tabularnewline
\hline 
 & 2 & 5 & 5.19e-01 & 6.51e-03 & 4.62e-05 & 1.59e-07 & 3.30e-10\tabularnewline
\hline 
 & 1 & 10 & 5.89e-01 & 6.24e-03 & 4.28e-05 & 1.40e-07 & 2.89e-10\tabularnewline
\hline 
20 & 20 & 1 & 2.61e-01 & 2.46e-03 & 3.38e-06 & 2.42e-09 & 5.94e-12\tabularnewline
\hline 
 & 10 & 2 & 2.25e-01 & 1.90e-03 & 3.06e-06 & 1.95e-09 & 7.21e-12\tabularnewline
\hline 
 & 4 & 5 & 2.03e-01 & 1.26e-03 & 2.31e-06 & 1.52e-09 & 7.25e-12\tabularnewline
\hline 
 & 1 & 20 & 2.02e-01 & 9.35e-04 & 1.93e-06 & 1.33e-09 & 5.97e-12\tabularnewline
\hline 
40 & 40 & 1 & 2.03e-01 & 5.84e-04 & 1.85e-07 & 2.60e-11 & 1.63e-11\tabularnewline
\hline 
 & 20 & 2 & 1.11e-01 & 4.50e-04 & 1.77e-07 & 2.41e-11 & 2.06e-11\tabularnewline
\hline 
 & 8 & 5 & 9.84e-02 & 2.94e-04 & 1.34e-07 & 2.46e-11 & 2.95e-11\tabularnewline
\hline 
 & 1 & 40 & 9.37e-02 & 1.66e-04 & 9.85e-08 & 2.06e-11 & 2.71e-11\tabularnewline
\hline 
80 & 80 & 1 & 1.88e-01 & 1.44e-04 & 1.11e-08 & 3.41e-11 & 7.17e-11\tabularnewline
\hline 
 & 40 & 2 & 5.74e-02 & 1.11e-04 & 1.08e-08 & 4.92e-11 & 8.36e-11\tabularnewline
\hline 
 & 16 & 5 & 5.29e-02 & 7.30e-05 & 8.27e-09 & 1.53e-10 & 1.92e-10\tabularnewline
\hline 
 & 1 & 80 & 4.63e-02 & 3.41e-05 & 5.61e-09 & 1.10e-10 & 1.27e-10\tabularnewline
\hline 
160 & 160 & 1 & 1.84e-01 & 3.59e-05 & 6.90e-10 & 2.90e-10 & 3.05e-10\tabularnewline
\hline 
 & 80 & 2 & 3.22e-02 & 2.77e-05 & 6.83e-10 & 1.98e-10 & 3.08e-10\tabularnewline
\hline 
 & 32 & 5 & 3.38e-02 & 1.82e-05 & 5.82e-10 & 9.50e-10 & 1.42e-09\tabularnewline
\hline 
 & 1 & 160 & 2.33e-02 & 7.69e-06 & 5.26e-10 & 6.52e-10 & 8.68e-10\tabularnewline
\hline 
320 & 320 & 1 & 1.83e-01 & 8.97e-06 & 5.69e-10 & 1.45e-09 & 1.53e-09\tabularnewline
\hline 
 & 160 & 2 & 2.18e-02 & 6.91e-06 & 5.10e-10 & 6.94e-10 & 1.38e-09\tabularnewline
\hline 
 & 64 & 5 & 2.70e-02 & 5.45e-06 & 1.93e-09 & 6.66e-09 & 9.14e-09\tabularnewline
\hline 
 & 1 & 320 & 1.18e-02 & 1.82e-06 & 3.09e-09 & 1.04e-08 & 8.18e-09\tabularnewline
\hline 
\end{tabular}
\end{table}

\begin{table}
\caption{Errors for the initial value solver in the norm $\hat{H}_{D,L}^{1}(0,5)$
in the case of using least-squares differentiation for determining
accurate transfer conditions. In case $n=Ln$, we have $L=1$ such
that the global least-squares method is used\label{tab:ErrorsLSQ}}

\centering{}%
\begin{tabular}{|c|c|c|c|c|c|c|c|}
\hline 
$L\times n$ & $L$ & $n\backslash N$ & 1 & 3 & 5 & 7 & 9\tabularnewline
\hline 
\hline 
10 & 10 & 1 & 3.68e+00 & 8.25e-02 & 1.03e-03 & 5.24e-06 & 1.49e-08\tabularnewline
\hline 
 & 5 & 2 & 3.35e+00 & 8.10e-02 & 7.66e-04 & 3.04e-06 & 8.32e-09\tabularnewline
\hline 
 & 2 & 5 & 3.01e+00 & 6.79e-02 & 6.37e-04 & 2.40e-06 & 6.25e-09\tabularnewline
\hline 
 & 1 & 10 & 2.59e+00 & 6.29e-02 & 5.71e-04 & 1.84e-06 & 3.89e-09\tabularnewline
\hline 
20 & 20 & 1 & 2.45e+00 & 2.61e-02 & 8.84e-05 & 9.32e-08 & 6.28e-11\tabularnewline
\hline 
 & 10 & 2 & 2.34e+00 & 2.38e-02 & 7.38e-05 & 6.18e-08 & 3.64e-11\tabularnewline
\hline 
 & 4 & 5 & 2.37e+00 & 2.06e-02 & 6.65e-05 & 5.33e-08 & 2.72e-11\tabularnewline
\hline 
 & 1 & 20 & 1.59e+00 & 1.76e-02 & 6.12e-05 & 4.52e-08 & 1.70e-11\tabularnewline
\hline 
40 & 40 & 1 & 2.21e+00 & 1.09e-02 & 9.61e-06 & 2.02e-09 & 1.49e-11\tabularnewline
\hline 
 & 20 & 2 & 2.08e+00 & 9.08e-03 & 8.58e-06 & 1.62e-09 & 2.01e-11\tabularnewline
\hline 
 & 8 & 5 & 1.91e+00 & 7.65e-03 & 7.84e-06 & 1.45e-09 & 2.22e-11\tabularnewline
\hline 
 & 1 & 40 & 1.26e+00 & 6.42e-03 & 7.31e-06 & 1.32e-09 & 2.17e-11\tabularnewline
\hline 
80 & 80 & 1 & 2.17e+00 & 5.14e-03 & 1.14e-06 & 5.09e-11 & 7.57e-11\tabularnewline
\hline 
 & 40 & 2 & 2.13e+00 & 4.12e-03 & 1.05e-06 & 5.49e-11 & 8.82e-11\tabularnewline
\hline 
 & 16 & 5 & 1.58e+00 & 3.40e-03 & 9.63e-07 & 7.67e-11 & 1.18e-10\tabularnewline
\hline 
 & 1 & 80 & 1.09e+00 & 2.84e-03 & 9.02e-07 & 9.62e-11 & 1.03e-10\tabularnewline
\hline 
160 & 160 & 1 & 2.16e+00 & 2.53e-03 & 1.40e-07 & 1.44e-10 & 3.59e-10\tabularnewline
\hline 
 & 80 & 2 & 2.15e+00 & 2.00e-03 & 1.31e-07 & 1.44e-10 & 4.09e-10\tabularnewline
\hline 
 & 32 & 5 & 1.85e+00 & 1.64e-03 & 1.20e-07 & 5.88e-10 & 7.81e-10\tabularnewline
\hline 
 & 1 & 160 & 8.84e-01 & 1.36e-03 & 1.12e-07 & 4.05e-10 & 8.53e-10\tabularnewline
\hline 
320 & 320 & 1 & 2.16e+00 & 1.26e-03 & 1.75e-08 & 5.08e-10 & 1.15e-09\tabularnewline
\hline 
 & 160 & 2 & 2.16e+00 & 9.94e-04 & 1.63e-08 & 5.40e-10 & 1.52e-09\tabularnewline
\hline 
 & 64 & 5 & 2.06e+00 & 8.13e-04 & 1.50e-08 & 3.67e-09 & 5.40e-09\tabularnewline
\hline 
 & 1 & 320 & 6.51e-01 & 6.74e-04 & 1.41e-08 & 6.67e-09 & 7.00e-09\tabularnewline
\hline 
\end{tabular}
\end{table}

\section{Conclusions and outlook}

In the present paper, we proposed an implementation of the reduction
procedure of \cite[Section 12]{RaRh02}, \cite{HaMa231} in order
to provide practical means to formulate accurate initial and transfer conditions for regular linear higher-index DAEs. This is a problem that does not occur at all with regular ODEs, but only emerges in the context of DAEs. To the best of the authors' knowledge, this is the first attempt ever to replace the rather difficult hand provision as in \cite{HaMa231} by a numerical method.
Stable
and accurate algorithms for that are investigated. The efficiency
and accuracy of the proposal has been demonstrated in numerical
examples. 
The implementation has then be used in a time step method which is
based on the least-squares approach advocated in \cite{HaMaPad121}.
The latter requires transfer conditions which must be formulated in
terms of accurate initial conditions. We showed that the combination
of these two approaches gives rise to an efficient initial value solver
for higher-index DAEs. A numerical example provides some insight
into the behavior of the method demonstrating its efficiency and accuracy.

A crucial assumption for the convergence of the time stepping approach
is Conjecture~AIC which is based on heuristic estimations and our numerical experiments. However, a strict proof of the convergence in terms of the
approximation polynomials and the discretization stepsize is missing so far. A thorough
proof would eliminate this gap in our derivations. Additionally, a
careful further rounding error analysis of Algorithm~\ref{alg:Reduction-procedure-1}
may contribute  to a better understanding of the numerical experiments.
In that respect, Proposition~\ref{prop:ChebCond} may be of use.

We believe that this work represents an appreciable step towards a fully numerical method for IVPs in higher-index DAEs which would undoubtedly be a great benefit. So far, there are no fully numerical methods for general higher-index DAEs available. 
The previously known methods for general higher-index DAEs  use so-called derivative arrays. This is very  effort-intensive and there are no usable estimates of errors.
There are approaches 
 which evaluate derivative arrays generated by automatic differentiations as in \cite{EstLam21,EstLam21b} and 
those that work with a derivative array provided from the beginning and a so-called remodeling into an index-one DAE \cite{KuhMeh06} as adumbrated in Remark \ref{r.array} above. Here we quote the assessment from \cite[p. 234]{TMS14} : \textit{to perform this reformulation... , there is currently no way to avoid differentiating the  whole model and to use full derivative arrays.} 
This shows and emphasizes the need to develop appropriate numerical methods together with error analyses.

\appendix

\section{Appendix}

\subsection{Choices of the basis functions in the procedures of Chistyakov and
Jansen\label{subsec:Choices}}

Let as in Algorithm~\ref{alg:Reduction-procedure}, for a given pair
$\{E,F\}$, $Y$ and $Z$ be smooth matrix functions such that $Y(t)$
is a basis of $\im E(t)$, and $Z(t)$ is a basis of $(\im E(t))^{\perp}$.
From (\ref{eq:EF}), we have $Z^{\T}Fx=0$. Chistyakov proposes to
construct a matrix function $K:[a,b]\rightarrow\R^{m\times m}$ being
nonsingular for every $t\in[a,b]$ such that $Z^{\T}FK=[F_{1},F_{2}]$
and $F_{2}(t)$ is nonsingular. Then, the basis $C(t)$ of $S(t)$
is given by
\begin{equation}
C=K\begin{bmatrix}
I_{r}\\
-F_{2}^{-1}F_{1}
\end{bmatrix}.\label{eq:ChistBasis}
\end{equation}
The elimination procedure by Jansen is more complex. Its aim is to
find a decomposition (and, thus, a basis) which requires only minimal
smoothness of the solution components. In \cite{Jansen14}, this decomposition
is mainly used for theoretical investigations. Moreover, a great effort
is spent in order to generalize this procedure to nonlinear problems.
In the present context, the construction of the basis function $C$ takes
the following steps:
\begin{itemize}
\item Let $Y,Z$ be given as before. Additionally, let $T(t)$ and $T^{c}(t)$
be bases of $\ker E(t)$ and $(\ker E(t))^{\perp}$, respectively.
\item Form $Z^{\T}FT$ and define bases $W(t)$ of $(\im Z^{\T}FT(t))^{\perp}$
and $V(t)$ of $\im Z^{\T}FT(t)$.
\item Let $G(t),H(t)$ be bases of $\ker W^{\T}Z^{\T}FT^{c}(t)$ and $\ker V^{\T}Z^{\T}FT(t)$,
respectively. Furthermore, let $G^{c}(t),H^{c}(t)$ be bases of $(\ker W^{\T}Z^{\T}FT^{c}(t))^{\perp}$
and $(\ker V^{\T}Z^{\T}FT(t))^{\perp}$, respectively.
\item Set $\mathfrak{E}=-(V^{\T}Z^{\T}FTH^{c})^{-1}V^{\T}Z^{\T}FT^{c}G$.
\end{itemize}
Then a suitable basis is obtained by
\begin{equation}
C=[T^{c}G+TH^{c}\mathfrak{E}\;\vert\; TH].\label{eq:JansenBasis}
\end{equation}
More details of this procedure can found in \cite[Appendix A.3]{HaMa231}.

\subsection{Proof of Proposition \ref{prop:ChebCond}\label{subsec:ProofA2}}

First, we observe that $\lvert w_{j}/w_{i}\rvert\leq4$. So it remains
to bound the sums
\[
s_{i}=\sum_{j=1,j\neq i}^{M}\frac{1}{\lvert\tilde{\sigma}_{i}-\tilde{\sigma}_{j}\rvert}=\sum_{j=1,j\neq i}^{M}\frac{1}{\lvert\cos\alpha_{i}-\cos\alpha_{j}\vert}
\]
where $\alpha_{i}=\frac{M-i}{M-1}\pi$. Since the nodes satisfy $\cos\alpha_{i}=-\cos\alpha_{M-i}$,
$i=1,\ldots,M$, it is sufficient to consider the sums $s_{i}$ for
$i\leq(M+1)/2$.

Let first $i=1$. Then it holds
\begin{align*}
s_{1} & =\sum_{j=2}^{M}\frac{1}{\cos\alpha_{j}+1}=\sum_{j=2}^{M}\frac{1}{2\cos^{2}\frac{\alpha_{j}}{2}}\\
 & =\frac{M-1}{\pi}\sum_{j=2}^{M}\frac{1}{2\cos^{2}\frac{\alpha_{j}}{2}}\frac{\pi}{M-1}.
\end{align*}
The latter sum can be interpreted as a Riemann sum for the integral
$\frac{M-1}{\pi}\int_{0}^{\alpha_{2}}\frac{1}{2}\cos^{-2}\frac{\alpha}{2}\dt\alpha$.
Since the integrand is monotonically increasing, we obtain
\begin{align*}
s_{1} & \leq\frac{1}{2}\cos^{-2}\frac{\alpha_{2}}{2}+\frac{M-1}{2\pi}\int_{0}^{\alpha_{2}}\cos^{-2}\frac{\alpha}{2}\dt\alpha\\
 & =\frac{1}{2}\cos^{-2}\frac{\alpha_{2}}{2}+\frac{M-1}{\pi}\left[\tan\frac{\alpha}{2}\right]_{\alpha=0}^{\alpha=\alpha_{2}}\\
 & =\frac{1}{2}\cos^{-2}\left(\frac{M-2}{M-1}\frac{\pi}{2}\right)+\frac{M-1}{\pi}\tan\left(\frac{M-2}{M-1}\frac{\pi}{2}\right).
\end{align*}
Since the sine function is convex in $[0,\frac{\pi}{2}]$, it holds
$\sin x\geq\frac{2}{\pi}x$ for $x\in[0,\frac{\pi}{2}]$. Hence,
\begin{align*}
\frac{1}{2}\cos^{-2}\left(\frac{M-2}{M-1}\frac{\pi}{2}\right) & =\frac{1}{2}\sin^{-2}\left(\frac{\pi}{2}-\frac{M-2}{M-1}\frac{\pi}{2}\right)=\frac{1}{2}\sin^{-2}\left(\frac{1}{M-1}\frac{\pi}{2}\right)\\
 & \leq\frac{1}{2}(M-1)^{2}.
\end{align*}
Similarly, it holds $\tan x\geq x$ for $x\in[0,\frac{\pi}{2}]$.
Hence,
\begin{align*}
\frac{M-1}{\pi}\tan\left(\frac{M-2}{M-1}\frac{\pi}{2}\right) & =\frac{M-1}{\pi}\cot\left(\frac{\pi}{2}-\frac{M-2}{M-1}\frac{\pi}{2}\right)=\frac{M-1}{\pi}\tan^{-1}\left(\frac{1}{M-1}\frac{\pi}{2}\right)\\
 & \leq\frac{2}{\pi^{2}}(M-1)^{2}.
\end{align*}
Hence,
\[
s_{1}\leq\left(\frac{1}{2}+\frac{2}{\pi^{2}}\right)(M-1)^{2}\leq(M-1)^{2}.
\]

Let $1<i\leq(M+1)/2$ now. This means in particular that $\alpha_{i}\in[\frac{\pi}{2},\frac{M-2}{M-1}\pi]\subseteq[\frac{\pi}{2},\pi)$.
The basic idea is similarly to the derivations before. It holds
\[
s_{i}=\sum_{j=1}^{j=i-1}\frac{1}{\cos\alpha_{i}-\cos\alpha_{j}}+\sum_{j=i+1}^{M}\frac{1}{\cos\alpha_{j}-\cos\alpha_{i}}.
\]
The first term can be interpreted as the Riemann sum for $\frac{M-1}{\pi}\int_{\alpha_{i-1}}^{\pi}(\cos\alpha_{i}-\cos\alpha)^{-1}\dt\alpha$
while the second term corresponds to the integral $\frac{M-1}{\pi}\int_{0}^{\alpha_{i+1}}(\cos\alpha-\cos\alpha_{i})^{-1}\dt\alpha$.
This leads to
\begin{align*}
s_{i} & \leq\frac{1}{\cos\alpha_{i}-\cos\alpha_{i-1}}+\frac{M-1}{\pi}\int_{\alpha_{i-1}}^{\pi}\frac{\dt\alpha}{\cos\alpha_{i}-\cos\alpha}\\
 & \phantom{\leq}+\frac{1}{\cos\alpha_{i+1}-\cos\alpha_{i}}+\frac{M-1}{\pi}\int_{0}^{\alpha_{i+1}}\frac{\dt\alpha}{\cos\alpha-\cos\alpha_{i}}.
\end{align*}
Use the integrals evaluated further below to obtain
\begin{align*}
s_{i} & \leq\frac{1}{\cos\alpha_{i}-\cos\alpha_{i-1}}+\frac{M-1}{\pi\sin\alpha_{i}}\ln\frac{\tan(\alpha_{i-1}/2)+\tan(\alpha_{i}/2)}{\tan(\alpha_{i-1}/2)-\tan(\alpha_{i}/2)}\\
 & \phantom{\leq}+\frac{1}{\cos\alpha_{i+1}-\cos\alpha_{i}}+\frac{M-1}{\pi\sin\alpha_{i}}\ln\frac{\tan(\alpha_{i}/2)+\tan(\alpha_{i+1}/2)}{\tan(\alpha_{i}/2)-\tan(\alpha_{i+1}/2)}.
\end{align*}
The inequality $\ln x\leq x-1$ is well-known. This provides us with
\begin{align*}
\frac{1}{\sin\alpha_{i}}\ln\frac{\tan(\alpha_{i-1}/2)+\tan(\alpha_{i}/2)}{\tan(\alpha_{i-1}/2)-\tan(\alpha_{i}/2)} & \leq\frac{1}{\sin\alpha_{i}}\left(\frac{\tan(\alpha_{i-1}/2)+\tan(\alpha_{i}/2)}{\tan(\alpha_{i-1}/2)-\tan(\alpha_{i}/2)}-1\right)\\
 & =\frac{1}{2\sin(\alpha_{i}/2)\cos(\alpha_{i}/2)}\left(\frac{2\tan(\alpha_{i}/2)}{\tan(\alpha_{i-1}/2)-\tan(\alpha_{i}/2)}\right)\\
 & =\frac{\tan(\alpha_{i}/2)\cos(\alpha_{i-1}/2)\cos(\alpha_{i}/2)}{\sin(\alpha_{i}/2)\cos(\alpha_{i}/2)\sin(\alpha_{i-1}/2-\alpha_{i}/2)}\\
 & =\frac{\cos(\alpha_{i-1}/2)}{\cos(\alpha_{i}/2)}\frac{1}{\sin\left(\frac{\pi}{2(M-1)}\right)}.
\end{align*}
We have $0<\alpha_{i}/2<\alpha_{i-1}/2<\pi/2$ such that $0<\cos(\alpha_{i-1}/2)<\cos(\alpha_{i}/2)$.
Hence,
\[
\frac{1}{\sin\alpha_{i}}\ln\frac{\tan(\alpha_{i-1}/2)+\tan(\alpha_{i}/2)}{\tan(\alpha_{i-1}/2)-\tan(\alpha_{i}/2)}\leq M-1.
\]
A similar estimation yields
\[
\frac{1}{\sin\alpha_{i}}\ln\frac{\tan(\alpha_{i}/2)+\tan(\alpha_{i+1}/2)}{\tan(\alpha_{i}/2)-\tan(\alpha_{i+1}/2)}\leq\frac{\sin(\alpha_{i+1}/2)}{\sin(\alpha_{i}/2)}\frac{1}{\sin\left(\frac{\pi}{2(M-1)}\right)}.
\]
It holds $0\leq\alpha_{i+1}/2<\alpha_{i}/2<\pi/2$. Hence, $0\leq\sin(\alpha_{i+1}/2)<\sin(\alpha_{i}/2)<1$
such that
\[
\frac{1}{\sin\alpha_{i}}\ln\frac{\tan(\alpha_{i}/2)+\tan(\alpha_{i+1}/2)}{\tan(\alpha_{i}/2)-\tan(\alpha_{i+1}/2)}\leq M-1.
\]
Furthermore,
\[
\frac{1}{\cos\alpha_{i}-\cos\alpha_{i-1}}=\frac{1}{2\sin((\alpha_{i}+\alpha_{i-1})/2)\sin\left(\frac{\pi}{2(M-1)}\right)}\leq\frac{1}{2\sin((\alpha_{i}+\alpha_{i-1})/2)}(M-1).
\]
Taking into account the conditions on $i$ we obtain
\[
\frac{1}{\cos\alpha_{i}-\cos\alpha_{i-1}}\leq\frac{3}{2}(M-1)^{2}.
\]
A similar estimate holds for $(\cos\alpha_{i+1}-\cos\alpha_{i})^{-1}$.
Consequently,
\[
s_{i}\leq3(M-1)^{2}+\frac{2}{\pi}(M-1)^{2}\leq4(M-1)^{2}.
\]
Summarizing the estimates, we obtain
\[
\lvert D_{2}^{N}\rvert_{\infty}=4\max\{s_{1},\ldots,s_{M}\}\leq16(M-1)^{2}.
\]

\subsubsection*{Evaluation of the integrals appearing in the proof of Proposition
\ref{prop:ChebCond}}

Consider the integral
\[
I_{\epsilon}=\int_{\alpha}^{\pi-\epsilon}\frac{\dt s}{a-\cos s}
\]
with $1>a>\cos\alpha$, $0<\alpha<\pi$, and $0<\epsilon<\pi-\alpha$.
Since the integrand is continuous, the limit $\lim_{\epsilon\rightarrow0}I_{\epsilon}=I_{0}$
exists and corresponds to the integral with $\epsilon=0$. Introduce
the change of variables $t=\tan\frac{s}{2}$. It holds $\dt t/\dt s=\frac{1}{2}(1+t^{2})$
and $\cos\alpha=\frac{1-t^{2}}{1+t^{2}}$. Hence, for $t_{2}=\tan\frac{1}{2}(\pi-\epsilon)$
and $t_{1}=\tan\frac{\alpha}{2}$,
\begin{align*}
I_{\epsilon} & =\int_{t_{1}}^{t_{2}}\frac{2}{1+t^{2}}\frac{1+t^{2}}{a(1+t^{2})-(1-t^{2})}\dt t\\
 & =\frac{2}{a+1}\int_{t_{1}}^{t_{2}}\frac{\dt t}{t^{2}+\frac{a-1}{a+1}}.
\end{align*}
According to the assumptions on $a$, it holds $(a-1)/(a+1)<0$ such
that 
\begin{align*}
I_{\epsilon} & =\frac{2}{a+1}\int_{t_{1}}^{t_{2}}\frac{\dt t}{t^{2}-b^{2}},\quad b=\sqrt{\frac{1-a}{1+a}}\\
 & =\frac{2}{(a+1)b}\int_{t_{1}/b}^{t_{2}/b}\frac{\dt\tau}{\tau^{2}-1}\\
 & =\frac{-1}{(a+1)b}\left[\ln\frac{\tau+1}{\tau-1}\right]_{\tau=t_{1}/b}^{\tau=t_{2}/b}.
\end{align*}
In the latter integral we used the fact that the denominator in $I_{\epsilon}$
is always positive, that is $\tau^{2}-1>0$. Moreover, $t=\tan\frac{s}{2}$
is positive in the domain of integration, hence, $\tau>1$. Then, for
the upper limit, $t_{2}/b\rightarrow+\infty$ such that
\[
I_{0}=\frac{1}{(a+1)b}\ln\frac{t_{1}/b+1}{t_{1}/b-1}=\frac{1}{\sqrt{1-a^{2}}}\ln\frac{\tan(\alpha/2)+b}{\tan(\alpha/2)-b},\quad b=\sqrt{\frac{1-a}{1+a}}.
\]
Consider now the integral
\[
J=\int_{0}^{\alpha}\frac{\dt s}{\cos s-a}.
\]
We assume $\cos\alpha>a>-1$. The substitution $t=\tan\frac{s}{2}$
provides, as before,
\begin{align*}
J & =\int_{0}^{t_{1}}\frac{2}{1+t^{2}}\frac{1+t^{2}}{1-t^{2}-a(1+t^{2})}\dt t\\
 & =\frac{-2}{a+1}\int_{0}^{t_{1}}\frac{\dt t}{t^{2}+\frac{a-1}{a+1}}
\end{align*}
where $t_{1}=\tan\frac{\alpha}{2}$. Similarly as before, $(a-1)/(a+1)<0$.
Hence, 
\[
J=\frac{-2}{(a+1)b}\int_{0}^{t_{1}/b}\frac{\dt\tau}{\tau^{2}-1}.
\]
In the present case, the denominator is negative, $\tau^{2}-1<0$.
Since $t=\tan\frac{s}{2}$ is nonnegative in the domain of integration,
$0\leq\tau<1$. This provides us with
\begin{align*}
J & =\frac{1}{\sqrt{a^{2}-1}}\left[\ln\frac{1+\tau}{1-\tau}\right]_{\tau=0}^{\tau=t_{1}/b}\\
 & =\frac{1}{\sqrt{1-a^{2}}}\ln\frac{b+\tan(\alpha/2)}{b-\tan(\alpha/2)},\quad b=\sqrt{\frac{1-a}{1+a}}.
\end{align*}
Note: If $a=\cos\beta$ with $0\leq\beta<\pi$, $b=\tan\frac{\beta}{2}$.
Moreover, $\sqrt{1-a^{2}}=\sin\beta$.

\subsection{Proof of Remark \ref{rem:CondEquidistant}}

Consider the equidistant nodes
\[
-1=\tilde{\sigma}_{1}<\cdots<\tilde{\sigma}_{M}=1
\]
where $\tilde{\sigma}_{i}=(i-1)\frac{2}{M-1}-1$, $i=1,\ldots,M$.
Then it holds
\[
w_{i}=(-1)^{i}\left(\begin{array}{c}
M-1\\
i-1
\end{array}\right)
\]
where $\left(\begin{array}{c}
n\\
k
\end{array}\right)$ denotes the binomial coefficient \cite{BerTre04}. The absolute row
sum $s_{i}$ for $D_{2}^{N}$ becomes
\begin{align*}
s_{i} & =\sum_{j=1,j\neq i}^{M}\left|\frac{w_{j}}{w_{i}}\right|\left|\frac{1}{\tilde{\tau}_{i}-\tilde{\tau}_{j}}\right|=\frac{M-1}{2}\frac{1}{\lvert w_{i}\rvert}\sum_{j=1,j\neq i}^{M}\lvert w_{j}\rvert\frac{1}{\lvert i-j\rvert}\\
 & \geq\frac{1}{2}\frac{1}{\lvert w_{i}\rvert}\sum_{j=1,j\neq i}^{M}\lvert w_{j}\rvert\\
 & \geq\frac{1}{2}\frac{1}{\lvert w_{i}\rvert}\left[\sum_{j=1}^{M}\lvert w_{j}\rvert-\lvert w_{i}\rvert\right]\\
 & \geq\frac{1}{2}\left(\begin{array}{c}
M-1\\
i-1
\end{array}\right)^{-1}\sum_{j=1}^{M}\left(\begin{array}{c}
M-1\\
j-1
\end{array}\right)-\frac{1}{2}\\
 & =\frac{1}{2}\left[\left(\begin{array}{c}
M-1\\
i-1
\end{array}\right)^{-1}2^{M-1}-1\right].
\end{align*}
Hence,
\[
\lvert D_{2}^{N}\rvert_{\infty}\geq\max_{i=1,\ldots,M}s_{i}=\frac{1}{2}(2^{M-1}-1).
\]
Similarly,
\begin{align*}
s_{i} & \leq\frac{M-1}{2}\frac{1}{\lvert w_{i}\rvert}\sum_{j=1,j\neq i}^{M}\lvert w_{j}\rvert\\
 & \leq\frac{M-1}{2}\left[\left(\begin{array}{c}
M-1\\
i-1
\end{array}\right)^{-1}2^{M-1}-1\right].
\end{align*}
Hence,
\[
\lvert D_{2}^{N}\rvert_{\infty}\leq\frac{M-1}{2}(2^{M-1}-1).
\]

\subsubsection*{Acknowlegement}

This preprint has not undergone
peer review or any post-submission improvements or corrections. The Version of Record of this article is
published in Numerical Algorithms, and is available online at https://doi.org/10.1007/s11075-025-02116-7.

\bibliographystyle{plain}
\bibliography{referenzen}

\end{document}